\documentclass[12pt]{amsart}

\usepackage{amsfonts}
\usepackage{amsmath}
\usepackage{amsthm}
\usepackage{amssymb}
\usepackage{latexsym}

\newtheorem{defn}{Definition}[section]
\newtheorem{lem}[defn]{Lemma}
\newtheorem{slem}[defn]{Sublemma}
\newtheorem{prop}[defn]{Proposition}
\newtheorem{thm}[defn]{Theorem}
\newtheorem{cor}[defn]{Corollary}
\newtheorem{ex}[defn]{Example}
\newtheorem{rem}[defn]{Remark}

\newenvironment{pf}{\proof}{\endproof}

\newcommand{\nc}{\newcommand}
\newcommand{\rnc}{\renewcommand}

\nc{\bal}{\begin{align*}}
\nc{\ball}{\begin{align}}
\nc{\bit}{\begin{itemize}}
\nc{\bdefn}{\begin{defn}}
\nc{\blem}{\begin{lem}}
\nc{\bslem}{\begin{slem}}
\nc{\bprop}{\begin{prop}}
\nc{\bthm}{\begin{thm}}
\nc{\bcor}{\begin{cor}}
\nc{\bex}{\begin{ex}}
\nc{\brem}{\begin{rem}}
\nc{\bpf}{\begin{pf}}
\nc{\lab}{\label}

\nc{\eit}{\end{itemize}}
\nc{\edefn}{\end{defn}}
\nc{\elem}{\end{lem}}
\nc{\eslem}{\end{slem}}
\nc{\eprop}{\end{prop}}
\nc{\ethm}{\end{thm}}
\nc{\ecor}{\end{cor}}
\nc{\eex}{\end{ex}}
\nc{\erem}{\end{rem}}
\nc{\epf}{\end{pf}}

\nc{\A}{\ensuremath{\mathbf{A}}}
\nc{\Dr}{\ensuremath{\mathbf{D}}}
\nc{\Q}{\ensuremath{\mathbf{Q}}}
\nc{\U}{\ensuremath{\mathbf{U}}}
\nc{\Z}{\ensuremath{\mathbf{Z}}}

\rnc{\a}{\ensuremath{\alpha}}
\rnc{\b}{\ensuremath{\beta}}
\rnc{\d}{\ensuremath{\delta}}
\rnc{\th}{\ensuremath{\theta}}
\nc{\pa}{\ensuremath{\varpi_{1}}}
\nc{\ps}{\ensuremath{\tilde{\psi}}}
\nc{\psd}{\ensuremath{\psi}}
\nc{\D}{\ensuremath{\Delta}}

\nc{\h}{\operatorname{h}}
\nc{\ih}{\operatorname{i}}
\nc{\sym}{\operatorname{Sym}}

\begin{document}

\title[An Integral PBW Basis of Type A$_{2}^{(2)}$]
{An Integral PBW Basis of \\ the Quantum Affine Algebra of Type A$_{2}^{(2)}$}
\author{Tatsuya Akasaka} 
\begin{abstract}
We construct an integral PBW basis and an integral crystal basis 
of the quantum affine algebra of type A$_{2}^{(2)}$.
\end{abstract}
\maketitle

{\allowdisplaybreaks


\section{Introduction.}
For the study of the precise structure of a quantum affine algebra $\U$, 
the construction of a good basis is an important problem.

Damiani constructed a PBW basis in the $A_{1}^{(1)}$ case \cite{Da1}.
Generalizing it, Beck constructed one
in the untwisted case \cite{B2}. 
Their bases are considered over $\Q(q)$.

However, to study the representations of $\U$ 
when $q$ is specialized at the roots of unity, 
one needs a $\Z [q,q^{-1}]$-basis of a certain 
$\Z [q,q^{-1}]$-subalgebra $\U_{\Z}$ of $\U$ (see \cite[Part 5]{L}). 
We call such a basis an integral basis of $\U$.
Also, to deal with the structure of the Grothendieck group of 
level 0 representations of $\U$,
an integral basis is needed (see \cite{K2}).

Around 1990, Kashiwara and Lusztig
introduced the notion of a global crystal basis 
(associated with a crystal basis) and a canonical basis respectively 
and proved the existence and uniqueness \cite{K1}\cite{L2}.
It turned out that these two notions are equivalent \cite{GL}.

Let $\U^{+}$ be the positive part of $\U$ and 
let $\U_{\Z}^{+}=\U_{\Z}\cap\U^{+}$.
Let {\bf B} be the global crystal basis of $\U^{+}$ and 
let $L$ be the $\Z[q^{-1}]$-lattice generated by {\bf B}.
Then, {\bf B} is a $\Q(q)$-basis of $\U^{+}$, 
a $\Z[q,q^{-1}]$-basis of $\U_{\Z}^{+}$, 
a $\Z[q^{-1}]$-basis of $L$, 
a $\Z[q]$-basis of $\overline{L}$, and 
a $\Z$-basis of $L\cap\overline{L}$.
Here, $-$ is the $\Q$-linear involution of $\U^{+}$ given by 
$\overline{q}=q^{-1},\ \overline{e_{i}}=e_{i}$.

An integral crystal basis $B$ of $\U^{+}$ is, by definition, 
a $\Z [q^{-1}]$-basis of $L$ that coincides with {\bf B} modulo $q^{-1}L$. 
Let $T$ be the transformation matrix with coefficients in $\Z [q,q^{-1}]$
between $B$ and $\overline{B}$, i.e. $B=T\overline{B}$ in the matrix form. 
There is a unique matrix $A$ with coefficients in $\Z [q^{-1}]$ 
such that $T=A^{-1}\overline{A}$. 
Then we have ${\bf B}=AB$.
Thus, we can recover {\bf B} from $B$.

Using the results in \cite{CP}, Beck, Chari, and Pressley 
constructed an integral PBW basis and an integral crystal basis 
in the simply-laced case \cite{BCP}.

The purpose of this paper is to construct an integral PBW basis and 
an integral crystal basis in the $A_{2}^{(2)}$ case.

Let us take a closer look at the results of this paper.
Let $\a_{0}$ be the long simple root and $\a_{1}$ the short simple root 
of type $A_{2}^{(2)}$.
Let $\U$ be the quantum affine algebra of type $A_{2}^{(2)}$ and 
let $\U_{\Z}$ be its integral form (see Definition \ref{defn-z-form}). 
Here, $\U$ is an algebra over $\Q(q_{1})$ and $\U_{\Z}$ is an algebra 
over $\Z [q_{1},q_{1}^{-1}]$ (see Section 2 for $q_{1}=q^{1/2}$).
Let $\U^{+}$ be the positive part of $\U$ and let 
$\U_{\Z}^{+}=\U_{\Z}\cap\U^{+}$.

We divide the set of positive real roots into two parts:
\bal
R_{re}^{+}(>)=\{ n\d +\a_{1},\ (2n+2)\d -\a_{0}|\ n\geq 0 \}
\end{align*}
and 
\bal
R_{re}^{+}(<)=\{ (n+1)\d -\a_{1},\ 2n\d +\a_{0}|\ n\geq 0 \}.
\end{align*}
We then define the total orders on $R_{re}^{+}(>)$ and $R_{re}^{+}(<)$ by
\bal
n\d+\a_{1}< (2n+2)\d-\a_{0}< (n+1)\d+\a_{1}
\end{align*}
and
\bal
(2n+2)\d+\a_{0}< (n+1)\d-\a_{1}< 2n\d+\a_{0},
\end{align*}
respectively. We also set $R_{im}^{+}=\{ n\d|\ n\geq 1 \}$, 
the set of positive imaginary roots.

Using the braid group action on $\U$ introduced by Lusztig (see Section 2), 
we define the real root vector associated with a real root in $R_{re}^{+}(>)$
as follows:
\bal
E_{n\d +\a_{1}}&=(T_{1}^{-1}T_{0}^{-1})^{n}(e_{1}),\\
E_{(2n+2)\d -\a_{0}}&=(T_{1}^{-1}T_{0}^{-1})^{n}T_{1}^{-1}(e_{0}).
\end{align*}
We also define the real root vector associated with a real root in 
$R_{re}^{+}(<)$ as follows:
\bal
E_{2n\d +\a_{0}}&=(T_{0}T_{1})^{n}(e_{0}),\\
E_{(n+1)\d -\a_{1}}&=(T_{0}T_{1})^{n}T_{0}(e_{1}).
\end{align*}

\bdefn
In the following, each $\Z_{\geq 0}^{(i)}$ is a copy of $\Z_{\geq 0}$, 
and $E_{\a}^{(c)}$ for a real root $\a$ and $c\in\Z_{\geq 0}$ 
denotes the divided power of $E_{\a}$ 
$($see $\rm{Definition}$ $\ref{defn-divided})$,
which belong to $\U_{\Z}^{+}$.
\bit
\item[(1)]
For ${\bf c}=({\bf c}_{i})\in\oplus_{i\in R_{re}^{+}(>)}\Z_{\geq 0}^{(i)}$, 
we set 
${\bf E_{c}}=E_{\a_{1}}^{({\bf c}_{\a_{1}})}
E_{2\d-\a_{0}}^{({\bf c}_{2\d-\a_{0}})}
E_{\d+\a_{1}}^{({\bf c}_{\d+\a_{1}})}\cdots$.
\item[(2)]
We set
$B(>)=
\{ {\bf E_{c}}|\ {\bf c}\in\oplus_{i\in R_{re}^{+}(>)}\Z_{\geq 0}^{(i)}\}$.
\item[(3)]
For ${\bf c}=({\bf c}_{i})\in\oplus_{i\in R_{re}^{+}(<)}\Z_{\geq 0}^{(i)}$, 
we set
${\bf E_{c}}=\cdots E_{2\d+\a_{0}}^{({\bf c}_{2\d+\a_{0}})}
E_{\d-\a_{1}}^{({\bf c}_{\d-\a_{1}})}E_{\a_{0}}^{({\bf c}_{\a_{0}})}$.
\item[(4)]
We set 
$B(<)=
\{ {\bf E_{c}}|\ {\bf c}\in\oplus_{i\in R_{re}^{+}(<)}\Z_{\geq 0}^{(i)}\}$.
\eit
\edefn

It is known that both of $B(>)$ and $B(<)$ are 
linearly independent over $\Q(q_{1})$ \cite[40.2.1]{L}.
Let $\U^{+}(>)$ (resp. $\U^{+}(<)$) be the vector subspace of $\U^{+}$  
over $\Q(q_{1})$ with basis $B(>)$ (resp. $B(<)$): 
it is known that they are algebras over $\Q(q_{1})$.

We define the imaginary root vector $\ps_{n}$ associated with 
an imaginary root $n\d\in R_{im}^{+}$
as follows:
\bal
\ps_{n}&=[E_{\d -\a_{1}},E_{(n-1)\d +\a_{1}}]_{q^{-1}}\\
&:=E_{\d -\a_{1}}E_{(n-1)\d +\a_{1}}-q^{-1}E_{(n-1)\d +\a_{1}}E_{\d -\a_{1}}.
\end{align*}
We prove that they are mutually commutative and 
algebraically independent over $\Q (q_{1})$. 
Let $\U^{+}(0)$ be the subalgebra of $\U^{+}$ generated by them.

We prove that the $\Q(q_{1})$-linear map 
$\U^{+}(>)\otimes\U^{+}(0)\otimes\U^{+}(<)\longrightarrow\U^{+}$
given by multiplication is an isomorphism:
this solves the problem raised in Lusztig's book \cite[40.2.5]{L}
in the $A_{2}^{(2)}$-case. 
Namely, we prove the following.

\bthm
For ${\bf c}=({\bf c}_{i})\in\oplus_{i\in\Z_{\geq 1}}\Z_{\geq 0}^{(i)}$, 
we set
${\bf E'_{c}}=\ps_{1}^{{\bf c}_{1}}\ps_{2}^{{\bf c}_{2}}
\ps_{3}^{{\bf c}_{3}}\cdots$.
Then, the following is a $\Q(q_{1})$-basis of $\U^{+}:$
\bal
\{ {\bf E_{c_{+}}} {\bf E'_{c_{0}}} {\bf E_{c_{-}}}|
\ {\bf c_{+}}\in\oplus_{i\in R_{re}^{+}(>)}\Z_{\geq 0}^{(i)},\ 
{\bf c_{0}}\in\oplus_{n\geq 1}\Z_{\geq 0}^{(n)},\ 
{\bf c_{-}}\in\oplus_{i\in R_{re}^{+}(<)}\Z_{\geq 0}^{(i)}\}.
\end{align*}
\ethm 

Modifying this basis,  
we construct a $\Z [q_{1},q_{1}^{-1}]$-basis of $\U_{\Z}^{+}$.
Let $\U_{\Z}^{+}(>)$ (resp. $\U_{\Z}^{+}(<)$) 
be the free $\Z [q_{1},q_{1}^{-1}]$-submodule of 
$\U_{\Z}^{+}$ with basis $B(>)$ (resp. $B(<)$). 
It turns out that they are algebras over $\Z [q_{1},q_{1}^{-1}]$ 
\cite[Prop.2.3]{BCP}. 
Now, if we choose $\U_{\Z}^{+}(0)$ as 
the $\Z [q_{1},q_{1}^{-1}]$-subalgebra of $\U_{\Z}^{+}$ 
generated by the $\ps_{n}$'s, then 
the $\Z [q_{1},q_{1}^{-1}]$-linear map 
$\U_{\Z}^{+}(>)\otimes_{\Z [q_{1},q_{1}^{-1}]}
\U_{\Z}^{+}(0)\otimes_{\Z [q_{1},q_{1}^{-1}]}\U_{\Z}^{+}(<)
\longrightarrow\U_{\Z}^{+}$
given by multiplication is injective, 
but it is not surjective. 
Therefore, in order to construct 
a $\Z [q_{1},q_{1}^{-1}]$-basis of $\U_{\Z}^{+}$, 
we have to find an appropriate definition of 
$\U_{\Z}^{+}(0)\subset\U_{\Z}^{+}\cap\U^{+}(0)$ 
so that the above multiplication morphism is an isomorphism. 
Instead of the $\ps_{n}$, we introduce  
the new imaginary root vectors $P_{n}\in\U^{+}(0)$ as follows: 
we set $P_{0}=1$ and 
\bal 
P_{n}=[2n]_{1}^{-1}\sum_{k=0}^{n-1}P_{k}\ps_{n-k}q^{-k}
\ \ \ \text{for}\ \ n\geq 1.
\end{align*} 
We prove that the $P_{n}$ belong to $\U_{\Z}^{+}$. 
This statement is not at all evident, whereas 
$\ps_{n}\in\U_{\Z}^{+}$ is evident.
Then we define $\U_{\Z}^{+}(0)$ as 
the $\Z [q_{1},q_{1}^{-1}]$-subalgebra of $\U_{\Z}^{+}$ 
generated by the $P_{n}$'s and 
prove that the $\Z [q_{1},q_{1}^{-1}]$-linear map 
$\U_{\Z}^{+}(>)\otimes_{\Z [q_{1},q_{1}^{-1}]}
\U_{\Z}^{+}(0)\otimes_{\Z [q_{1},q_{1}^{-1}]}\U_{\Z}^{+}(<)
\longrightarrow\U_{\Z}^{+}$ is an isomorphism. 
In this way, we obtain

\bthm
For ${\bf c}=({\bf c}_{i})\in\oplus_{i\in\Z_{\geq 1}}\Z_{\geq 0}^{(i)}$, 
we set
${\bf E_{c}}=P_{1}^{{\bf c}_{1}}P_{2}^{{\bf c}_{2}}
P_{3}^{{\bf c}_{3}}\cdots$. 
Then, the following is a $\Z [q_{1},q_{1}^{-1}]$-basis of $\U_{\Z}^{+}:$
\bal
\{ {\bf E_{c_{+}}} {\bf E_{c_{0}}} {\bf E_{c_{-}}}|\ 
{\bf c_{+}}\in\oplus_{i\in R_{re}^{+}(>)}\Z_{\geq 0}^{(i)},\ 
{\bf c_{0}}\in\oplus_{n\geq 1}\Z_{\geq 0}^{(n)},\ 
{\bf c_{-}}\in\oplus_{i\in R_{re}^{+}(<)}\Z_{\geq 0}^{(i)}\}.
\end{align*}
\ethm

However, this basis does not give an integral crystal basis: 
we need a further modification.
For a partition $\lambda$, we define ${\bf S}_{\lambda}\in\U_{\Z}^{+}(0)$ 
from the $P_{n}$
in the same way as the Schur functions are defined from the complete symmetric 
functions: namely, we set 
\bal
{\bf S}_{\lambda}=\det (P_{\lambda_{i}-i+j})_{i,j\geq 1}
\ \ \text{for}\ \ \lambda=(\lambda_{1}\geq\lambda_{2}\geq\cdots)
\end{align*}
where we understand that $P_{n}=0$ for $n\leq -1$.
Then, it follows that the ${\bf S}_{\lambda}$ are 
quasi-orthonormal with respect to the inner product on $\U^{+}$ 
introduced by Drinfeld, 
that is,  
\bal
({\bf S}_{\lambda},{\bf S}_{\mu})\equiv \d_{\lambda,\mu}\mod q_{1}^{-1}\A.
\end{align*}
Here, $\A=\Q(q_{1})\cap\Z[[q_{1}^{-1}]]\subset\Q((q_{1}^{-1}))$.  
Therefore, the following is another 
$\Z [q_{1},q_{1}^{-1}]$-basis of $\U_{\Z}^{+}$:
\bal
\{ {\bf E_{c_{+}}} {\bf S}_{\lambda} {\bf E_{c_{-}}}|
\ {\bf c_{+}}\in\oplus_{i\in R_{re}^{+}(>)}\Z_{\geq 0}^{(i)},
\ \lambda\ \text{is a partition},
\ {\bf c_{-}}\in\oplus_{i\in R_{re}^{+}(<)}\Z_{\geq 0}^{(i)}\},
\end{align*}
which is denoted by $B$. In view of \cite[40.2.4]{L}, we see that 
$B$ is quasi-orthonormal with respect to the inner product. 
Hence, by the same argument in \cite{BCP}, we obtain

\bthm
$B$ is an integral crystal bases of $\U^{+}$.
\ethm

The contents of this paper are as follows.

In Section 2, we fix notations.
Automorphism $T_{\pa}=T_{0}T_{1}$ and anti-automorphism $T_{1}^{-1}\ast$ 
of $\U$ play an important role in this paper.

In Section 3, we introduce the root vectors as above and 
study their commutation relations. 
We prove that the imaginary root vectors mutually commute and 
are invariant under $T_{\pa}=T_{0}T_{1}$ and $T_{1}^{-1}\ast$. 
The key step in this section 
is to express the real root vectors recursively 
using brackets (Corollary \ref{c-1}), 
which the author learned from \cite[8.2]{KhT}.

In Section 4, 
we introduce the subspaces $\U^{+}(>),\U^{+}(<),\U^{+}(0)$ of $\U^{+}$ 
as above, and prove that the $\Q(q_{1})$-linear map 
$\U^{+}(>)\otimes\U^{+}(0)\otimes\U^{+}(<)\longrightarrow\U^{+}$
given by multiplication is surjective; 
the proof of its injectivity 
with the help of \cite[40.1.2]{L} is postponed until Section 6. 
We also introduce the $P_{n}$ in this section.

In Section 5, we study the coproducts of the real root vectors.

In Section 6,
we calculate the coproducts and 
the inner products of the imaginary root vectors and
introduce the ${\bf S}_{\lambda}$; 
the results in this section are used to construct 
an integral crystal basis in Section 8. 
As a by-product, the algebraically independence of the imaginary root vectors 
$($the $\ps_{n}$'s or the $P_{n}$'s$)$ is proved; 
thus, the monomials of the imaginary root vectors 
$($the $\ps_{n}$'s or the $P_{n}$'s$)$
form a basis of $\U^{+}(0)$ and we obtain PBW bases of $\U^{+}$.

Section 7 is the preparation for the next section. 

In Section 8, we give the commutation relation 
between $E_{\a_{0}}^{(s)}$ and $E_{\a_{1}}^{(r)}$. Using it, we 
prove that the $P_{n}$ belong to $\U_{\Z}^{+}$ 
and construct an integral PBW basis of $\U^{+}$.
Then we obtain an integral crystal basis of $\U^{+}$.

In Appendix A, we give some commutation relation
between the real root vectors that is used in Section 4.

In Appendix B, we discuss the connection between our root vectors and 
the Drinfeld generators.

After writing up the main part of this paper, the author
learned the existence of \cite{Da2}, 
in which (non-integral) PBW bases of twisted quantum affine algebras 
are constructed.

The author is grateful to his advisor Professor Masaki Kashiwara
for valuable comments and useful discussions on this work.
Thanks are also due to Professor Tetsuji Miwa for his support.


\section{Notation}
Let $X,\ Y$ be the finitely generated free $\Z$-modules with a perfect pairing 
$\langle\ ,\ \rangle: Y\times X\longrightarrow\Z$. 
Let $\a_{0}\in X$ be the long simple root and 
$\a_{1}\in X$ the short simple root of type $A_{2}^{(2)}$:
we assume that they are linearly independent.
We set $I=\{ 0,1\}$.
For $i\in I$, let $h_{i}\in Y$ be the simple coroots:
we assume that they are linearly independent.
Let $Q=\Z\a_{0}\oplus\Z\a_{1}\subset X$ be the root lattice
and let 
$Q^{+}=\Z_{\geq 0}\a_{0}\oplus\Z_{\geq 0}\a_{1},\ Q^{-}=-Q^{+}$.
Let $(\ ,\ )$ be the $\Z$-valued symmetric bilinear form on $X$ such that 
\bal
(\a_{0},\a_{0})=4,\ (\a_{0},\a_{1})=-2,\ (\a_{1},\a_{1})=1,
\end{align*}
so that the Cartan matrix is given by 
\bal
(\langle h_{i},\a_{j}\rangle)_{i,j}=
\left( 
\begin{array}{cc}
a_{00} &a_{01}\\ 
a_{10} &a_{11}
\end{array} 
\right) 
= \left( 
\begin{array}{cc}
2 &-1 \\ 
-4 &2
\end{array} 
\right).
\end{align*}
Let $\d =\a_{0} +2\a_{1}\in X$ be the smallest positive imaginary root. 
Note that $Q=\Z\d\oplus\Z\a_{1}$.
We set 
\bal
R_{re}^{+}(>)=\{ n\d +\a_{1},\ (2n+2)\d -\a_{0}|\ n\geq 0 \}
\end{align*}
and define the total order on it by
\bal
n\d+\a_{1}< (2n+2)\d-\a_{0}< (n+1)\d+\a_{1}
\end{align*}
for $n\geq 0$; we set 
\bal
R_{re}^{+}(<)=\{ (n+1)\d -\a_{1},\ 2n\d +\a_{0}|\ n\geq 0 \}
\end{align*}
and define the total order on it by
\bal
(2n+2)\d+\a_{0}< (n+1)\d-\a_{1}< 2n\d+\a_{0}
\end{align*}
for $n\geq 0$; 
we also set $R_{im}^{+}=\{ n\d |\ n\geq 1 \}$.
Then, $R_{re}^{+}=R_{re}^{+}(>)\sqcup R_{re}^{+}(<)$
is the set of positive real roots;
$R^{+}=R_{re}^{+}\sqcup R_{im}^{+}$ is the set of positive roots; and
$R=R^{+}\sqcup (-R^{+})$ is the set of roots.

We set $q_{0}=q^{2},\ q_{1}=q^{1/2}$. 
For $i\in I,\ k\in\Z,\ n\in\Z_{\geq 1}$, we set
$[k]_{i}=\frac{q_{i}^{k}-q_{i}^{-k}}{q_{i}-q_{i}^{-1}},
\ [n]_{i}!=\prod_{p=1}^{l}[p]_{i},\ [0]_{i}!=1,
\ [k]=\frac{q^{k}-q^{-k}}{q-q^{-1}},
\ [n]!=\prod_{p=1}^{l}[p],\ [0]!=1$.
For $i\in I,\ n,m\geq 0$, we set 
$\genfrac{[}{]}
{0pt}
{}{n+m}{m}_{i}
=\frac{[n+m]_{i}!}{[n]_{i}![m]_{i}!}$, 
which belong to $\Z [q_{i},q_{i}^{-1}]$.

\bdefn
Let $\U$ be the quantum affine algebra of type $A_{2}^{(2)}$, 
which is the $\Q (q_{1})$-algebra generated by
\bal
\{ q^{h},\ e_{i},\ f_{i}|\ h\in 2^{-1}Y,\ i\in I\}
\end{align*}
with the following defining relations$:$
\bal
&q^{0}=1,\ q^{h}q^{h'}=q^{h+h'}\ \ \text{for}\ h,h'\in 2^{-1}Y,\tag{U1}\\
&q^{h}e_{i}q^{-h}=q^{\langle h,\a_{i}\rangle}e_{i}\ \ 
\text{for}\ h\in 2^{-1}Y,\ i\in I,\tag{U2}\\
&q^{h}f_{i}q^{-h}=q^{-\langle h,\a_{i}\rangle}e_{i}\ \ 
\text{for}\ h\in 2^{-1}Y,\ i\in I,\tag{U3}\\
&[e_{i},f_{j}]=\d_{ij}\frac{k_{i}-k_{i}^{-1}}{q_{i}-q_{i}^{-1}}\ \ 
\text{for}\ i,j\in I,\tag{U4}\\
&\sum_{k=0}^{1-a_{ij}}(-1)^{k}e_{i}^{(k)}e_{j}e_{i}^{(1-a_{ij}-k)}=0\ \  
\text{for}\ i,j\in I\ \text{with}\ i\neq j,\tag{U5}\\
&\sum_{k=0}^{1-a_{ij}}(-1)^{k}f_{i}^{(k)}f_{j}f_{i}^{(1-a_{ij}-k)}=0\ \  
\text{for}\ i,j\in I\ \text{with}\ i\neq j\tag{U6}
\end{align*} 
where we set $k_{0}=q^{2h_{0}},\ k_{1}=q^{2^{-1}h_{1}}$ and 
$e_{i}^{(k)}=e_{i}^{k}/[k]_{i}!,\ f_{i}^{(k)}=f_{i}^{k}/[k]_{i}!$ 
for $i\in I,\ k\geq 0$.
\edefn

\brem
We have 
$k_{i}e_{j}k_{i}^{-1}=q^{(\a_{i},\a_{j})}e_{j},
\ k_{i}f_{j}k_{i}^{-1}=q^{-(\a_{i},\a_{j})}f_{j}$ for $i,j\in I$.
\erem

\bdefn
For $\mu =n\a_{0}+m\a_{1}\in Q$, we set $k_{\mu}=k_{0}^{n}k_{1}^{m}$.
We also set $c=k_{\d}=k_{0}k_{1}^{2}=q^{2h_{0}+h_{1}}$, 
which is a central element of $\U$.
\edefn

\bdefn
\bit
\item[(1)]
Let $\U^{0}$ be the $\Q (q_{1})$-algebra generated by $q^{h}$ 
for $h\in 2^{-1}Y$ with the defining relations \thetag{U1}.
\item[(2)]
Let $\U^{+}$ be the $\Q (q_{1})$-algebra generated by $e_{i}$ for $i\in I$ 
with the defining relations \thetag{U5}.
\item[(3)]
Let $\U^{-}$ be the $\Q (q_{1})$-algebra generated by $f_{i}$ for $i\in I$ 
with the defining relations \thetag{U6}.
\eit
\edefn
Then, $\{ q^{h}\in\U^{0}|\ h\in 2^{-1}Y\}$ 
is a $\Q (q_{1})$-basis of $\U^{0}$.

\bprop\lab{prop-triangular}\cite[3.2.5]{L}
The $\Q (q_{1})$-linear map 
$\U^{-}\otimes\U^{0}\otimes\U^{+}\rightarrow\U$ 
given by multiplication is an isomorphism.
\eprop

Hence, $\U^{0},\ \U^{+},\ \U^{-}$ can be considered as 
the subalgebras of $\U$.

If a nonzero element $x$ of $\U^{+}$ has a homogeneous expression 
in terms of $e_{i}$ ($i\in I$), then the indices $i_{1},\ldots ,i_{k}$ 
($i_{j}\in I,\ k\geq 0$) appearing in it 
are uniquely determined up to permutation, 
and we say that $x$ is homogeneous 
of weight $\sum_{j=1}^{k}\a_{i_{j}}$,
which is denoted by $|x|$. 
We apply the similar definition for $\U^{-}$.

\bdefn
For $\a\in Q^{+}$, we set
\bal
\U^{+}_{\a}=
\{ x\in \U^{+}|\ x=0\ \text{or}\ x\ \text{is homogeneous of weight}\ \a \}.
\end{align*}
We also set
\bal
\U^{+,h}=
\{ x\in \U^{+}|\ x=0\ \text{or}\ x\ \text{is homogeneous}\},
\end{align*}
which is closed under multiplication.
For $\a\in Q^{-}$, we define $\U^{-}_{\a}$ similarly.
\edefn

We have 
$\U^{+,h}=\sqcup_{\a\in Q^{+}}(\U^{+}_{\a}\backslash\{0\})\sqcup\{0\}$.

\bdefn
For $\nu\in Q$, we set 
$\U_{\nu}=\oplus_{\lambda\in Q^{-},\mu\in Q^{+};\lambda +\mu =\nu}
\U^{-}_{\lambda}\U^{0}\U^{+}_{\mu}$.
\edefn

\bdefn
For a subset $A$ of $\U^{+}$ and for $i\in\Z$, we write 
\bal
&A^{h}=A\cap\U^{+,h},\tag{1}\\
&A_{\leq i}=A\cap\oplus_{\a=n\d+r\a_{1}\in Q^{+};r\leq i}\U^{+}_{\a},\tag{2}\\
&A_{\geq i}=A\cap\oplus_{\a=n\d+r\a_{1}\in Q^{+};r\geq i}\U^{+}_{\a}.\tag{3}
\end{align*} 
\edefn

\bdefn
We define the function $\h$ from $\U^{+,h}\backslash\{0\}$ to $\Z$ by
$\h (x)=r$ for $x\in\U^{+}_{n\d+r\a_{1}}\backslash\{0\}$.
We call $\h (x)$ the height of $x$.
We also define the function $\ih$ from $\U^{+,h}\backslash\{0\}$ to 
$\Z_{\geq 0}$ by $\ih (x)=n$ 
for $x\in\U^{+}_{n\d+r\a_{1}}\backslash\{0\}$.
\edefn

We have $\h (xy)=\h (x)+\h (y)$ and $\ih (xy)=\ih (x)+\ih (y)$ 
for $x,y\in\U^{+,h}\backslash\{0\}$.

\bdefn
For elements $x,y$ of a $\Q (q_{1})$-algebra and for 
$v\in\Q (q_{1})^{\times}$, we set
\bal
[x,y]_{v}=xy-vyx.
\end{align*}
When $v=1$, we omit the suffix $v$.
\edefn

\blem\label{lem-bracket}
Let $x,y,z$ be elements of a $\Q (q_{1})$-algebra and 
let $\a,\b\in\Q (q_{1})^{\times}$. Then,
\bit
\item[(1)] 
$[x,yz]_{\a}=[x,y]_{\b}z+y[x,z]_{\a /\b}\b$,
\item[(2)]
$[xy,z]_{\a}=x[y,z]_{\b}+[x,z]_{\a /\b}y\b$.
\eit 
\elem
\bpf
This is clear.
\epf

\bdefn\label{defn-z-form}
Let $\U_{\Z}$ be the $\Z [q_{1},q_{1}^{-1}]$-subalgebra of $\U$ generated by
$e_{i}^{(r)},\ f_{i}^{(r)},\ k_{i},\ k_{i}^{-1}$\ for $i\in I,\ r\geq 0$.
\edefn

\bdefn\label{defn-z-k}
For $i\in I,\ m\in\Z,\ r\in\Z_{\geq 0}$, we set
\bal
\genfrac{[}{]}
{0pt}
{}{k_{i},m}{r}_{i}
= \prod_{s=1}^r \frac{k_{i} q_{i}^{m-s+1} - k_{i}^{-1}q_{i}^{-m+ s-1}}
{q_{i}^s - q_{i}^{-s}}.
\end{align*}
We understand that 
$\genfrac{[}{]}
{0pt}
{}{k_{i},m}{0}_{i}
=1$.
Then, they belong to $\U_{\Z}$, which follows from 
$\rm{Lemma}\ \ref{lem-com-ef}$ below.
\edefn

\blem\label{lem-com-ef}\cite[3.1.9]{L}
Let $n,m\geq 0,\ i\in I$. Then,
\bal
e_{i}^{(n)}f_{i}^{(m)}=\sum_{t=0}^{\min (n,m)}f_{i}^{(m-t)}
\genfrac{[}{]}
{0pt}
{}{k_{i},2t-n-m}{t}_{i} 
e_{i}^{(n-t)}
\end{align*}
\elem

\bcor\label{cor-com-ef}
Let $r\geq 1,\ i\in I$. Then,
\bal 
[e_{i}^{(r)},f_{i}]= e_{i}^{(r-1)}  
\frac{q_{i}^{r-1}k_{i}-q_{i}^{1-r}k_{i}^{-1}}{q_{i}-q_{i}^{-1}}.
\end{align*}
\ecor
\bpf
This follows from Lemma \ref{lem-com-ef} with $m=1$.
\epf

\blem\label{lem-com-ek,kf}
Let $i,j\in I,\ m\in\Z,\ r\in\Z_{\geq 0}$. Then,
\bit
\item[(1)]
$e_{j}
\genfrac{[}{]}
{0pt}
{}{k_{i},m}{t}_{i} 
=\genfrac{[}{]}
{0pt}
{}{k_{i},m-a_{ij}}{t}_{i} 
e_{j}$,
\item[(2)]
$\genfrac{[}{]}
{0pt}
{}{k_{i},m}{t}_{i} 
f_{j}=f_{j}
\genfrac{[}{]}
{0pt}
{}{k_{i},m-a_{ij}}{t}_{i}$.
\eit
\elem
\bpf
This is directly checked.
\epf

\bdefn
\bit
\item[(1)]
Let $\U_{\Z}^{0}$ be the $\Z [q_{1},q_{1}^{-1}]$-subalgebra of $\U$ 
generated by 
$k_{i},\ k_{i}^{-1},
\ \genfrac{[}{]}
{0pt}
{}{k_{i},m}{r}_{i}$ 
for $i\in I,\ m\in\Z,\ r\in\Z_{\geq 1}$.
\item[(2)]
Let $\U_{\Z}^{+}$ be the $\Z [q_{1},q_{1}^{-1}]$-subalgebra of $\U$ 
generated by $e_{i}^{(r)}$ for $i\in I,\ r\geq 0$.
\item[(3)]
Let $\U_{\Z}^{-}$ be the $\Z [q_{1},q_{1}^{-1}]$-subalgebra of $\U$ 
generated by $f_{i}^{(r)}$ for $i\in I,\ r\geq 0$.
\eit
\edefn

We have $\U_{\Z}^{+}=\oplus_{\a\in Q^{+}}(\U_{\Z}^{+}\cap\U_{\a}^{+})$ and  
$\U_{\Z}^{-}=\oplus_{\a\in Q^{-}}(\U_{\Z}^{-}\cap\U_{\a}^{-})$, 
since the generators of $\U_{\Z}^{+}$ and $\U_{\Z}^{-}$ are homogeneous.

\blem\label{lem-z-k-basis}\cite[4.5]{L1}
The following set $B_{0}$ is a $\Z [q_{1},q_{1}^{-1}]$-basis of $\U_{\Z}^{0}:$
\bal
\Big{\{} k_{0}^{a}k_{1}^{b}
\genfrac{[}{]}
{0pt}
{}{k_{0},0}{r}_{0}
\genfrac{[}{]}
{0pt}
{}{k_{1},0}{s}_{1}
\Big{|}\ a,b\in\{0,1\},\ r,s\geq 0\Big{\}}.
\end{align*}
\elem

\bprop\lab{prop-z-triangular}
The $\Z [q_{1},q_{1}^{-1}]$-linear map 
$\U_{\Z}^{-}\ \otimes_{\Z [q_{1},q_{1}^{-1}]}
\U_{\Z}^{0}\ \otimes_{\Z [q_{1},q_{1}^{-1}]}\U_{\Z}^{+}\rightarrow\U_{\Z}$ 
given by multiplication is an isomorphism.
\eprop

\bpf
Note that $\U_{\Z}^{+}$ and $\U_{\Z}^{-}$ are free 
$\Z [q_{1},q_{1}^{-1}]$-modules by 
the existence of the canonical bases:
thus, the injectivity follows from Proposition \ref{prop-triangular}.
The surjectivity follows from 
Lemma \ref{lem-com-ef} and Lemma \ref{lem-com-ek,kf}.
\epf

Noting that $\U_{\Z}^{0}$, $\U_{\Z}^{+}$, $\U_{\Z}^{-}$ are free 
$\Z [q_{1},q_{1}^{-1}]$-modules, by Proposition \ref{prop-z-triangular}, 
we have $\U_{\Z}^{+}=\U_{\Z}\cap\U^{+}$.

\bdefn\cite[3.4.1]{K1}\cite[1.2.13]{L} 
Let $i\in I$.
We define the $\Q (q_{1})$-linear maps $r_{i}$ and $\ _{i}r$
from $\U^{+}$ to itself by $r_{i}(1)=\ _{i}r(1)=0,
\ r_{i}(e_{j})=\ _{i}r(e_{j})=\d_{ij}$, and 
\bit
\item[(1)]
$r_{i}(xy)=q^{(\a_{i},|y|)}r_{i}(x)y+xr_{i}(y)$\ \ for $x,y\in\U^{+,h}$,
\item[(2)]
$_{i}r(xy)=\ _{i}r(x)y+q^{(\a_{i},|x|)}x\ _{i}r(y)$\ \ for $x,y\in\U^{+,h}$.
\eit
\edefn

\blem\lab{r-1}\cite[3.4.2]{K1}\cite[1.2.15]{L} 
Let $\a\in Q^{+}\backslash\{0\},\ x\in\U^{+,h}_{\a}$.
\bit
\item[(1)]
If $r_{i}(x)=0$ for all $i\in I$, then $x=0$.
\item[(2)]
If $_{i}r(x)=0$ for all $i\in I$, then $x=0$.
\eit
\elem

\blem\lab{r-2}\cite[3.4.1]{K1}\cite[3.1.6]{L}
Let $x\in\U^{+,h},\ i\in I$. Then,
\bal
[x,f_{i}]=\frac{r_{i}(x)k_{i}-k_{i}^{-1}\ _{i}r(x)}{q_{i}-q_{i}^{-1}}.
\end{align*}
\elem

Let us recall the braid group action on $\U$ introduced by Lusztig.
Note that the braid group of type $A_{2}^{(2)}$ is the free group generated by 
$T_{0}$ and $T_{1}$.

\bdefn\label{defn-braid}\cite[37.1.3]{L}
For $i\in I$, we define the automorphisms $T_{i}$ of $\U$ by 
\bal
&T_{i}(q^{h})=q^{h-\langle h,\a_{i}\rangle \a_{i}}
\ \ \text{for}\ \ h\in 2^{-1}Y,\\
&T_{i}(e_{i})=-f_{i}k_{i},\
T_{i}(f_{i})=-k_{i}^{-1}e_{i},\\
&T_{i}(e_{j})
=\sum_{r=0}^{-a_{ij}}(-1)^{r}q_{i}^{-r}e_{i}^{(-a_{ij}-r)}e_{j}e_{i}^{(r)}
\ \ \text{for}\ \ j\in I,\ j\neq i,\\
&T_{i}(f_{j})
=\sum_{r=0}^{-a_{ij}}(-1)^{r}q_{i}^{r}f_{i}^{(r)}f_{j}f_{i}^{(-a_{ij}-r)}
\ \ \text{for}\ \ j\in I,\ j\neq i.
\end{align*}
Then we have
\bal
&T_{i}^{-1}(q^{h})=q^{h-\langle h,\a_{i}\rangle \a_{i}}
\ \ \text{for}\ \ h\in 2^{-1}Y,\\
&T_{i}^{-1}(e_{i})=-k_{i}^{-1}f_{i},\
T_{i}^{-1}(f_{i})=-e_{i}k_{i},\\
&T_{i}^{-1}(e_{j})
=\sum_{r=0}^{-a_{ij}}(-1)^{r}q_{i}^{-r}e_{i}^{(r)}e_{j}e_{i}^{(-a_{ij}-r)}
\ \ \ \text{if}\ \ j\in I,\ j\neq i,\\
&T_{i}^{-1}(f_{j})
=\sum_{r=0}^{-a_{ij}}(-1)^{r}q_{i}^{r}f_{i}^{(-a_{ij}-r)}f_{j}f_{i}^{(r)}
\ \ \ \text{for}\ \ j\in I,\ j\neq i.
\end{align*}
Moreover, $\U_{\Z}$ is closed under the braid group action.
\edefn

We have
$T_{i}(k_{j})=T_{i}^{-1}(k_{j})=k_{i}^{-a_{ij}}k_{j}$ for $i,j\in I$.

\brem
In \cite{L}, $T_{i}$ and $T_{i}^{-1}$ are denoted 
by $T_{i,1}^{''}$ and $T_{i,-1}^{'}$ respectively.
\erem

\bdefn
We set $T_{\pa}=T_{0}T_{1}$.
\edefn

We have $T_{\pa}(c)=c,\ T_{\pa}(k_{1})=c^{-1}k_{1}$ and 
$T_{\pa}(k_{0})=c^{2}k_{0}$. 

\bdefn
Let $\D$ be the coproduct of $\U$ given by
\bal
\D (q^{h})&=q^{h}\otimes q^{h},\\
\D (e_{i})&=e_{i}\otimes 1 + k_{i}\otimes e_{i},\\
\D (f_{i})&=f_{i}\otimes k_{i}^{-1} + 1\otimes f_{i}
\end{align*}
for $h\in 2^{-1}Y$, $i\in I$.
\edefn

\bdefn 
Let $-$ be the $\Q$-linear involution of $\U$ given  by
\begin{align*}
\ \overline{q_{1}}=q_{1}^{-1}, 
\ \overline{q^{h}}=q^{-h}, 
\ \overline{e_{i}}=e_{i},
\ \overline{f_{i}}=f_{i}\ \ \text{for}\ \ h\in 2^{-1}Y,\ i\in I.
\end{align*}

\edefn

\bdefn 
Let $\ast$ be the $\Q (q_{1})$-linear anti-involution of $\U$ 
given by
\begin{align*}
\ast(q^{h})=q^{-h},\ \ast (e_{i})=e_{i},\ \ast (f_{i})=f_{i}
\ \ \text{for}\ \ h\in 2^{-1}Y,\ i\in I.
\end{align*}
\edefn

We have $\ast T_{i}=T_{i}^{-1}\ast$ for $i\in I$.

\bdefn 
Let $\Omega$ be the $\Q$-linear anti-involution of $\U$ 
given  by
\begin{align*}
\Omega (q_{1})=q_{1}^{-1},\ 
\Omega (q^{h})=q^{-h},\  
\Omega (e_{i})=f_{i},\ 
\Omega (f_{i})=e_{i}\ \ \text{for}\ \ h\in 2^{-1}Y,\ i\in I.
\end{align*}
\edefn

We have $\Omega T_{i}=T_{i}\Omega$ for $i\in I$.
We also have $\Omega (\U^{+})=\U^{-}$ and $\Omega (\U_{\Z}^{+})=\U_{\Z}^{-}$.

\brem\label{rem-basis}
Given a $\Q (q_{1})$-basis $B$ of $\U^{+}$, 
we obtain a $\Q (q_{1})$-basis of $\U:$ 
\bal
\{ \Omega (b)q^{h}b|\ b\in B,\ h\in 2^{-1}Y\}.
\end{align*}
Given a $\Z [q_{1},q_{1}^{-1}]$-basis $B$ of $\U^{+}_{\Z}$, 
by $\rm{Lemma}\ \ref{lem-z-k-basis}$, 
we obtain a $\Z [q_{1},q_{1}^{-1}]$-basis of $\U_{\Z}:$ 
\bal
\{ \Omega (b)xb|\ x\in B_{0},\ b\in B\}.
\end{align*}
\erem

\brem
Let us set the general notation used in this paper.
\bit
\item[(1)]
Algebras are always assumed to be associative with $1$
and any morphism of algebras is assumed to preserve $1$.
\item[(2)]
For a statement $P$, the symbol $\th (P)$ means $1$ if $P$ is true 
and $0$ otherwise.
\item[(3)]
For a set $S$ and for $i,j\in S$,
we write $\d_{ij}=\th(i=j)$.
\item[(4)]
Any summation considered in $\U$ is taken over a finite set.
\item[(5)]
For vector subspaces $A,B$ of $\U$, we denote by $A\otimes B$ 
the tensor product of $A$ and $B$ over $\Q(q_{1})$.
\item[(6)]
For subsets $A,B$ of $\U$, we denote by $AB$ the $\Z$-submodule of $\U$ 
generated by $\{ ab|\ a\in A,\ b\in B\}:$ thus, 
for a subring $C$ of $\Q(q_{1})$, if $A$ or $B$ is a $C$-module, then 
so is $AB$.
\eit
\erem


\section{Root Vectors}

\bdefn\label{defn-E}
For $n\geq 1$, we set 
\bal
&E_{n\d +\a_{1}}=(T_{1}^{-1}T_{0}^{-1})^{n}(e_{1}),\\
&E_{(2n+2)\d -\a_{0}}=(T_{1}^{-1}T_{0}^{-1})^{n}T_{1}^{-1}(e_{0}),\\
&E_{2n\d +\a_{0}}=(T_{0}T_{1})^{n}(e_{0}),\\
&E_{(n+1)\d -\a_{1}}=(T_{0}T_{1})^{n}T_{0}(e_{1}).
\end{align*}
They are homogeneous and the suffices indicate their weights.
We call them real root vectors. 
\edefn

\bdefn\label{defn-divided}
For $\a\in R^{+}_{re}$ such that $(\a,\a)=(\a_{i},\a_{i})$ and for $k\geq 0$,
we write $E_{\a}^{(k)}=E_{\a}^{k}/[k]_{i}!$.
They belong to $\U_{\Z}^{+}$ by \cite[37.1.3]{L} and \cite[40.2.3]{L}.
\edefn

\bdefn\label{defn-ps} 
For $n\geq 1$, we define the elements $\ps_{n}$ of $\U_{\Z}^{+}$ by
\bal
\ps_{n}&=[E_{\d -\a_{1}},E_{(n-1)\d +\a_{1}}]_{q^{-1}}.
\end{align*}
We call them imaginary root vectors.
\edefn

At the end of this section, we shall show that the imaginary root vectors are 
mutually commutative (Proposition \ref{p-3}).

We also set
\bal
\ps_{0}=(q_{1}-q_{1}^{-1})^{-1}.
\end{align*}

The following two lemmas follow from Definition \ref{defn-E} and 
Definition \ref{defn-braid}, 
and will be often used in this paper.
\blem\label{l-1}
Under the action of the automorphism $T_{\pa}$, we have
\bal
&T_{\pa}(E_{n\d -\a_{1}})=E_{(n+1)\d -\a_{1}} \ \ \text{for \ $n\geq 1$},\\
&T_{\pa}(E_{2n\d +\a_{0}})=E_{(2n+2)\d +\a_{0}}\ \ \text{for \ $n\geq 0$},\\
&T_{\pa}(E_{n\d +\a_{1}})=E_{(n-1)\d +\a_{1}}\ \ \text{for \ $n\geq 1$},\\
&T_{\pa}(E_{2n\d -\a_{0}})=E_{(2n-2)\d -\a_{0}}\ \ \text{for \ $n\geq 2$},\\
&T_{\pa}^{-1}(E_{\d -\a_{1}})=-k_{1}^{-1}f_{1},\\
&T_{\pa}(E_{2\d -\a_{0}})=-f_{0}k_{0}.
\end{align*}
\elem

\blem\label{l-2}
Under the action of the anti-involution $T_{1}^{-1}\ast$, we have
\bal
&(T_{1}^{-1}\ast )(E_{n\d -\a_{1}})=E_{n\d +\a_{1}} 
\ \ \text{for \ $n\geq 1$},\\
&(T_{1}^{-1}\ast )(E_{2n\d +\a_{0}})=E_{(2n+2)\d -\a_{0}}
\ \ \text{for \ $n\geq 0$},\\
&(T_{1}^{-1}\ast )(E_{\a_{1}})=-k_{1}^{-1}f_{1}.
\end{align*}
\elem

We are going to express the real root vectors 
without using the braid group action (Corollary \ref{c-1}).

\blem\label{l-3}
We have 
\bit
\item[(1)]
$E_{\d -\a_{1}}=[E_{\a_{0}},E_{\a_{1}}]_{q^{-2}}$,

\item[(2)]
$\ps_{1}=[2]_{1}\sum_{r=0}^{2}
(-1)^{r}q_{1}^{-3r}e_{1}^{(r)}e_{0}e_{1}^{(2-r)}$,

\item[(3)]
$E_{\d +\a_{1}}=\sum_{r=0}^{3}
(-1)^{r}q^{-r}e_{1}^{(r)}e_{0}e_{1}^{(3-r)}$,

\item[(4)]
$[E_{\d -\a_{1}},f_{1}]=[4]_{1}k_{1}E_{\a_{0}}$,

\item[(5)]
$[\ps_{1},f_{1}]=[3]_{1}!k_{1}E_{\d -\a_{1}}$,

\item[(6)]
$[E_{\d +\a_{1}},f_{1}]=k_{1}\ps_{1}$.
\eit
\elem

\bpf
(1)
By Definition \ref{defn-E}, we have 
\bal
E_{\d -\a_{1}}=T_{0}(e_{1})
=\sum_{r=0}^{1}(-1)^{r}q_{0}^{-r}e_{0}^{(1-r)}e_{1}e _{0}^{(r)}
=[E_{\a_{0}},E_{\a_{1}}]_{q^{-2}}.
\end{align*}

(2)
By Definition \ref{defn-ps} and (1), we have 
\bal
\ps_{1}&=[E_{\d -\a_{1}},E_{\a_{1}}]_{q^{-1}}
=[[E_{\a_{0}},E_{\a_{1}}]_{q^{-2}},E_{\a_{1}}]_{q^{-1}}\\
&=(e_{0}e_{1}-q^{-2}e_{1}e_{0})e_{1}-q^{-1}e_{1}(e_{0}e_{1}-q^{-2}e_{1}e_{0})\\
&=e_{0}e_{1}^{2}-(q^{-1}+q^{-2})e_{1}e_{0}e_{1}+q^{-3}e_{1}^{2}e_{0}\\
&=[2]_{1}\sum_{r=0}^{2}\ (-1)^{r}q_{1}^{-3r}e_{1}^{(r)}e_{0}e_{1}^{(2-r)}.
\end{align*}

(3)
By Definition \ref{defn-E}, we have 
\bal
E_{\d +\a_{1}}&=T_{1}^{-1}T_{0}^{-1}(e_{1})
=T_{1}^{-1}\Big(\sum_{r=0}^{1}(-1)^{r}q_{0}^{-r}e_{0}^{(r)}e_{1}e _{0}^{(1-r)}
\Big)\\
&=T_{1}^{-1}([e_{1},e_{0}]_{q^{-2}})
=\Big[-k_{1}^{-1}f_{1},
\sum_{r=0}^{4}(-1)^{r}q_{1}^{-r}e_{1}^{(r)}e_{0}e _{1}^{(4-r)}\Big]_{q^{-2}}\\
&=k_{1}^{-1}\Big[\sum_{r=0}^{4}(-1)^{r}q_{1}^{-r}e_{1}^{(r)}
e_{0}e _{1}^{(4-r)},f_{1}\Big].
\end{align*}
Applying Corollary \ref{cor-com-ef}, we have 
\bal
E_{\d +\a_{1}}&=k_{1}^{-1}\sum_{r=0}^{4}(-1)^{r}q_{1}^{-r}
\Big( e_{1}^{(r-1)}\frac{q_{1}^{r-1}k_{1}-q_{1}^{1-r}k_{1}^{-1}}
{q_{1}-q_{1}^{-1}} e_{0}e _{1}^{(4-r)}\\
&\ + e_{1}^{(r)}e_{0}e _{1}^{(3-r)}
\frac{q_{1}^{3-r}k_{1}-q_{1}^{r-3}k_{1}^{-1}}{q_{1}-q_{1}^{-1}}\Big)\\
&=k_{1}^{-1}\sum_{r=0}^{4}(-1)^{r}q_{1}^{-r}
\frac{q_{1}^{1-r}k_{1}-q_{1}^{r-1}k_{1}^{-1}}{q_{1}-q_{1}^{-1}}
(e_{1}^{(r-1)}e_{0}e _{1}^{(4-r)}+e_{1}^{(r)}e_{0}e _{1}^{(3-r)})\\
&=k_{1}^{-1}\sum_{r=0}^{3}\Big( (-1)^{r+1}q_{1}^{-r-1}
\frac{q_{1}^{-r}k_{1}-q_{1}^{r}k_{1}^{-1}}{q_{1}-q_{1}^{-1}}\\
&\ +(-1)^{r}q_{1}^{-r}\frac{q_{1}^{1-r}k_{1}-q_{1}^{r-1}k_{1}^{-1}}
{q_{1}-q_{1}^{-1}}\Big) e_{1}^{(r)}e_{0}e_{1}^{(3-r)}\\
&=\sum_{r=0}^{3}(-1)^{r}q^{-r}e_{1}^{(r)}e_{0}e_{1}^{(3-r)}.
\end{align*}

(4)
By (1), we have 
\bal
[E_{\d -\a_{1}},f_{1}]&=[[e_{0},e_{1}]_{q^{-2}},f_{1}]
=[e_{0},[e_{1},f_{1}]]_{q^{-2}}
=\Big[e_{0},\frac{k_{1}-k_{1}^{-1}}{q_{1}-q_{1}^{-1}}\Big]_{q^{-2}}\\
&=\frac{q^{2}k_{1}-q^{-2}k_{1}^{-1}}{q_{1}-q_{1}^{-1}} e_{0}
 -q^{-2}\frac{k_{1}-k_{1}^{-1}}{q_{1}-q_{1}^{-1}} e_{0}
=[4]_{1}k_{1}E_{\a_{0}}.
\end{align*}

(5)
By Definition \ref{defn-ps} and (4), we have 
\bal
[\ps_{1},f_{1}]
&=[[E_{\d -\a_{1}},E_{\a_{1}}]_{q^{-1}},f_{1}]\\
&=[[E_{\d -\a_{1}},f_{1}],E_{\a_{1}}]_{q^{-1}}
+[E_{\d -\a_{1}},[E_{\a_{1}},f_{1}]]_{q^{-1}}\\
&=[4]_{1}[k_{1}E_{\a_{0}},E_{\a_{1}}]_{q^{-1}}
+\Big[E_{\d -\a_{1}},\frac{k_{1}-k_{1}^{-1}}{q_{1}-q_{1}^{-1}}\Big]_{q^{-1}}\\
&=[4]_{1}k_{1}[E_{\a_{0}},E_{\a_{1}}]_{q^{-2}}
+\frac{qk_{1}-q^{-1}k_{1}^{-1}}{q_{1}-q_{1}^{-1}}E_{\d -\a_{1}} 
-q^{-1}\frac{k_{1}-k_{1}^{-1}}{q_{1}-q_{1}^{-1}}E_{\d -\a_{1}}\\
&=([4]_{1}+[2]_{1})k_{1}E_{\d -\a_{1}}\\
&=[3]_{1}!\ k_{1}E_{\d -\a_{1}}.
\end{align*}

(6)
By (3) and Corollary   \ref{cor-com-ef}, we have 
\bal
[E_{\d +\a_{1}},f_{1}]&=\sum_{r=0}^{3}(-1)^{r}q^{-r}
[e_{1}^{(r)}e_{0}e_{1}^{(3-r)},f_{1}]\\
&=\sum_{r=0}^{3}\ (-1)^{r}q^{-r}\Big( e_{1}^{(r-1)}
\frac{q_{1}^{r-1}k_{1}-q_{1}^{1-r}k_{1}^{-1}}{q_{1}-q_{1}^{-1}}
e_{0}e_{1}^{(3-r)}\\
&\ +e_{1}^{(r)}e_{0}e_{1}^{(2-r)}\frac{q_{1}^{2-r}k_{1}-q_{1}^{r-2}k_{1}^{-1}}
{q_{1}-q_{1}^{-1}}\Big) \\
&=\sum_{r=0}^{3}\ (-1)^{r}q^{-r}\Big(
\frac{q_{1}^{1-r}k_{1}-q_{1}^{r-1}k_{1}^{-1}}{q_{1}-q_{1}^{-1}}
e_{1}^{(r-1)}e_{0}e_{1}^{(3-r)}\\
&\ +\frac{q_{1}^{2-r}k_{1}-q_{1}^{r-2}k_{1}^{-1}}{q_{1}-q_{1}^{-1}}
e_{1}^{(r)}e_{0}e_{1}^{(2-r)}\Big)\\
&=\sum_{r=0}^{2}\Big( (-1)^{r+1}q^{-r-1}
\frac{q_{1}^{-r}k_{1}-q_{1}^{r}k_{1}^{-1}}{q_{1}-q_{1}^{-1}}\\
&\ +(-1)^{r}q^{-r}\frac{q_{1}^{2-r}k_{1}-q_{1}^{r-2}k_{1}^{-1}}
{q_{1}-q_{1}^{-1}}\Big) e_{1}^{(r)}e_{0}e_{1}^{(2-r)}\\
&=\sum_{r=0}^{2}\ (-1)^{r}q_{1}^{-3r}[2]_{1}k_{1}e_{1}^{(r)}e_{0}e_{1}^{(2-r)},
\end{align*}
which is equal to $k_{1}\ps_{1}$ by (2).
\epf

\blem\label{l-4}
Let $n\geq 1$. Then,
\bal
[E_{n\d +\a_{1}},f_{1}]=k_{1}T_{\pa}^{-1}(\ps_{n}).
\end{align*}
\elem
\bpf 
By Lemma \ref{l-1}, we have 
\bal
[E_{n\d +\a_{1}},f_{1}]
&=[E_{n\d +\a_{1}},-k_{1}T_{\pa}^{-1}(E_{\d -\a_{1}})]
=k_{1}[T_{\pa}^{-1}(E_{\d -\a_{1}}),E_{n\d +\a_{1}}]_{q^{-1}}\\
&=k_{1}T_{\pa}^{-1}([E_{\d -\a_{1}},E_{(n-1)\d +\a_{1}}]_{q^{-1}}).
\end{align*}
Hence, the lemma follows from Definition \ref{defn-ps}.
\epf

\blem{\label{l-5}}
We have 
\bit
\item[(1)]
$T_{\pa}(\ps_{1})=\ps_{1}$,
\item[(2)]
$(T_{1}^{-1}\ast )(\ps_{1})=\ps_{1}$.
\eit
\elem

\bpf
By Lemma \ref{l-4} (6) and Lemma \ref{l-4}, we have
$T_{\pa}^{-1}(\ps_{1})=\ps_{1}$, obtaining (1).
By Definition \ref{defn-ps}, we have
\bal
(T_{1}^{-1}\ast )(\ps_{1})
&=(T_{1}^{-1}\ast )([E_{\d -\a_{1}},E_{\a_{1}}]_{q^{-1}})
=[-k_{1}^{-1}f_{1},E_{\d +\a_{1}}]_{q^{-1}}\\
&=k_{1}^{-1}[E_{\d +\a_{1}},f_{1}].
\end{align*}
Applying Lemma \ref{l-4} and (1), we obtain (2).
\epf

\blem{\label{l-6}} 
We have 
\bit
\item[(1)]
$[E_{2\d -\a_{0}},E_{\a_{1}}]_{q^{2}}=0$,

\item[(2)]
$[E_{\a_{0}},E_{\d -\a_{1}}]_{q^{2}}=0$,

\item[(3)]
$[\ps_{1},E_{\a_{1}}]=[3]_{1}!\ E_{\d +\a_{1}}$,

\item[(4)]
$[E_{\a_{0}},\ps_{1}]=(q-1)[3]_{1}E_{\d -\a_{1}}^{2}$,

\item[(5)]
$[E_{\d +\a_{1}},E_{\a_{1}}]_{q}=[4]_{1}E_{2\d -\a_{0}}$.
\eit
\elem

\bpf
(1)
By Definition \ref{defn-E}, we have 
\bal
[E_{2\d -\a_{0}},E_{\a_{1}}]_{q^{2}}
&=[T_{1}^{-1}(E_{\a_{0}}),E_{\a_{1}}]_{q^{2}}
=T_{1}^{-1}([E_{\a_{0}},T_{1}E_{\a_{1}}]_{q^{2}})\\
&=T_{1}^{-1}([e_{0},-f_{1}k_{1}]_{q^{2}})
=-T_{1}^{-1}([e_{0},f_{1}]k_{1})=0.
\end{align*}

(2)
By Definition \ref{defn-E}, we have 
\bal
[E_{\a_{0}},E_{\d -\a_{1}}]_{q^{2}}
&=[E_{\a_{0}},T_{0}(E_{\a_{1}})]_{q^{2}}
=T_{0}([T_{0}^{-1}(E_{\a_{0}}),E_{\a_{1}}]_{q^{2}})\\
&=T_{0}([-k_{0}^{-1}f_{0},e_{1}]_{q^{2}})
=T_{0}(k_{0}^{-1}[e_{1},f_{0}])=0.
\end{align*}

(3)
By Lemma \ref{l-5} (2) and Lemma \ref{l-3} (5), we have 
\bal
(T_{1}^{-1}\ast )([\ps_{1},E_{\a_{1}}])
=[-k_{1}^{-1}f_{1},\ps_{1}]
=k_{1}^{-1}[\ps_{1},f_{1}]
=[3]_{1}!\ E_{\d -\a_{1}}.
\end{align*}
Applying Lemma \ref{l-2}, we obtain (3).

(4)
By Definition \ref{defn-ps} and Lemma \ref{lem-bracket}, we have 
\bal
[E_{\a_{0}},\ps_{1}]
&=[E_{\a_{0}},E_{\d -\a_{1}}E_{\a_{1}}-q^{-1}E_{\a_{1}}E_{\d -\a_{1}}] \\
&=[E_{\a_{0}},E_{\d -\a_{1}}]_{q^{2}}E_{\a_{1}}
+E_{\d -\a_{1}}[E_{\a_{0}},E_{\a_{1}}]_{q^{-2}}q^{2}\\
&\ -[E_{\a_{0}},E_{\a_{1}}]_{q^{-2}}E_{\d -\a_{1}}q^{-1}
-E_{\a_{1}}[E_{\a_{0}},E_{\d -\a_{1}}]_{q^{2}}q^{-3}.
\end{align*}
The first term and the last one vanish by (2). 
Thus, by Lemma \ref{l-3} (1), we obtain (4).

(5)
By Lemma \ref{l-3} (4), we have 
\bal
(T_{1}^{-1}\ast )([E_{\d +\a_{1}},E_{\a_{1}}]_{q})
&=[-k_{1}^{-1}f_{1},E_{\d -\a_{1}}]_{q}
=k_{1}^{-1}[E_{\d -\a_{1}},f_{1}]=[4]_{1}E_{\a_{0}}.
\end{align*}
Applying Lemma \ref{l-2}, we obtain (3). 
\epf

\bcor{\lab{c-1}}
We have 
\bit
\item[(1)] 
$[E_{(2n+2)\d -\a_{0}},E_{n\d +\a_{1}}]_{q^{2}}=0$ \ \ \ for \ $n\geq 0$,

\item[(2)]
$[E_{n\d -\a_{1}},E_{2n\d +\a_{0}}]_{q^{2}}=0$ \ \ \ for \ $n\geq 1$,

\item[(3)]
$[E_{2n\d +\a_{0}},E_{(n+1)\d -\a_{1}}]_{q^{2}}=0$ \ \ \ for \ $n\geq 0$,

\item[(4)]
$[E_{n\d +\a_{1}},E_{2n\d -\a_{0}}]_{q^{2}}=0$ \ \ \ for \ $n\geq 1$,

\item[(5)]
$[E_{(n+1)\d +\a_{1}},E_{n\d +\a_{1}}]_{q}=[4]_{1}E_{(2n+2)\d -\a_{0}}$
\ \ \ for \ $n\geq 0$,
\item[(6)]
$[E_{n\d -\a_{1}},E_{(n+1)\d -\a_{1}}]_{q}=[4]_{1}E_{2n\d +\a_{0}}$
\ \ \ for \ $n\geq 1$,
\item[(7)]
$[\ps_{1},E_{n\d +\a_{1}}]=[3]_{1}!\ E_{(n+1)\d +\a_{1}}$
\ \ \ for \ $n\geq 0$,
\item[(8)]
$[E_{n\d -\a_{1}},\ps_{1}]=[3]_{1}!\ E_{(n+1)\d -\a_{1}}$
\ \ \ for \ $n\geq 1$,
\item[(9)]
$[E_{2n\d +\a_{0}},\ps_{1}]=(q-1)[3]_{1}E_{(n+1)\d -\a_{1}}^{2}$
\ \ \ for \ $n\geq 0$,
\item[(10)]
$[\ps_{1},E_{2n\d -\a_{0}}]=(q-1)[3]_{1}E_{n\d +\a_{1}}^{2}$
\ \ \ for \ $n\geq 1$.
\eit
\ecor

\bpf 
We use Lemma \ref{l-1}, Lemma \ref{l-2}, and Lemma \ref{l-5}.

Applying $T_{\pa}^{n}$ to Lemma \ref{l-6} (1), we obtain (1).

Applying $T_{1}^{-1}\ast$ to (1) with $n\geq 1$, we obtain (2).

Applying $T_{\pa}^{n}$ to Lemma \ref{l-6} (2), we obtain (3).

Applying $T_{1}^{-1}\ast$ to (3) with $n$ replaced by $n-1$, we obtain (4).

Applying $T_{\pa}^{n}$ to Lemma \ref{l-6} (5), we obtain (5).

Applying $T_{1}^{-1}\ast$ to (5), we obtain (6).

Applying $T_{\pa}^{n}$ to Lemma \ref{l-6} (3), we obtain (7).

Applying $T_{1}^{-1}\ast$ to (7), we obtain (8).

Applying $T_{\pa}^{n}$ to Lemma \ref{l-6} (4), we obtain (9).

Applying $T_{1}^{-1}\ast$ to (9) with $n$ replaced by $n-1$, we obtain (10).
\epf

We introduce certain elements of $\Z[q,q^{-1}]$
that will be used for the commutation relations between the root vectors.

\bdefn\label{defn-b}
For $n\in\Z$, we define the elements $b_{n}$ of $\Z[q,q^{-1}]$ as follows$:$
for $i\in\Z$, we set
\bal
&b_{2i}=(1-q+q^{2})^{-1}(1-q)((-q)^{-i}-q^{2i}(q+q^{-1})),\\
&b_{2i+1}=(1-q+q^{2})^{-1}((-q)^{-i}+q^{2i+1}(q-1)).
\end{align*} 
Note that the roots $e^{\pm\pi i/3}$ of $1-q+q^{2}$
are roots of both of the last factors.
\edefn

\brem
It follows from $\rm{Definition}\ \ref{defn-b}$
that $b_{-n}=-q^{-1}\overline{b_{n}}$ for $n\in\Z$
where $-$ is the automorphism of $\Z[q,q^{-1}]$ given by $\overline{q}=q^{-1}$.
\erem
 
\bex  
We have 
\bal
&b_{-2}=1-q^{-3},\ b_{-1}=-q^{-1},\ b_{0}=1-q^{-1}\\
&b_{1}=1,\ b_{2}=q^{2}-q^{-1}=(q-1)[3]_{1},\\
&b_{3}=q^{2}-1-q^{-1}=b_{2}-1,\ b_{4}=q^{4}-q-1+q^{-2}.
\end{align*}
\eex

\blem\lab{b-1}
Let $i\in\Z$. Then, 
\bit
\item[(1)]
$b_{2i+1}+b_{2i-1}=b_{2i}$,
\item[(2)]
$q b_{2i}+(q+q^{-1})(q-1)\ b_{2i+1}=b_{2(i+1)}$,
\item[(3)]
$qb_{2i}+b_{2}b_{2(i+1)}=b_{2(i+2)}$.
\eit
\elem

\bpf 
We can directly check this by using Definition \ref{defn-b}.
\epf

\blem\lab{b-2}
Let $k,l\in\Z$. Then, 
\bit
\item[(1)]
$b_{2(k+1)}b_{2(l+1)}+qb_{2k}b_{2l}=b_{2(k+l+2)}-b_{2(k+l)}$,
\item[(2)]
$b_{2k}b_{2l+1}+qb_{2(k-1)}b_{2(l-1)+1}=b_{2(k+l)}-b_{2(k+l-1)}$,
\item[(3)]
$b_{2k}b_{2l}-(q+q^{-1})(q-1) b_{2k-1}b_{2l-1}=b_{2(k+l)}-b_{2(k+l-1)}$.
\eit
\elem

\bpf 
We can directly check (1) and (2) by using Definition \ref{defn-b}.
We check (3): by Lemma \ref{b-1} ((2), (1)), 
the left hand side is equal to 
\bal
b_{2k}b_{2l}-(b_{2k}-qb_{2k-2})b_{2l-1}
&=b_{2k}(b_{2l}-b_{2l-1})+qb_{2(k-1)}b_{2l-1}\\
&=b_{2k}b_{2l+1}+qb_{2(k-1)}b_{2l-1},
\end{align*}
which is equal to the right hand side by (2).
\epf

\blem\lab{b-3}
Let $s\in\Z$, $l\in\Z_{\geq{1}}$. Then,
\bal
q b_{2s}&+(1+q)\sum_{i=1}^{l-1}b_{2i}b_{2(i+s)}\th(l\geq 2)+b_{2l}b_{2(l+s)}
=b_{2(2l+s)},\tag{1}\\
q b_{2s}&+(1+q)\sum_{i=1}^{l-1}b_{2i}b_{2(i+s)}\th(l\geq 2)\tag{2}\\
&+(q+q^{-1})(q-1)b_{2l-1}b_{2(l+s)-1}=b_{2(2l+s-1)}.
\end{align*}
\elem

\bpf 
(1) is checked by the induction on $l$: 
the case where $l=1$ is nothing but Lemma \ref{b-1} (3), and 
the induction proceeds by Lemma \ref{b-2} (1).
(2) follows from (1) and Lemma \ref{b-2} (3).
\epf

\blem\lab{b-4}
Let $l\geq{1}$. Then,
\bal
&\sum_{i=1}^{l} b_{2i}[2l-2i+1]_{1}q^{i}=(q-1)q^{l}[l][2l+1]_{1},\tag{1}\\
&\sum_{i=1}^{l} b_{2i-1}[2l-2i+1]_{1}q^{i}=q^{l}[l]_{0}.\tag{2}
\end{align*}
\elem

\bpf 
The lemma is reduced to the following identities in $\Z[q,q^{-1}][[u]]$:
\bal
&\Big( \sum_{i\geq 1}b_{2i}q^{i}u^{i}\Big)
\Big( \sum_{k\geq 0}[2k+1]_{1}u^{k}\Big)
=(q-1)\sum_{l\geq 1}q^{l}[l][2l+1]_{1}u^{l},\\
&\Big( \sum_{i\geq 1}b_{2i-1}q^{i}u^{i}\Big)
\Big(\sum_{k\geq 0}[2k+1]_{1}u^{k}\Big)
=\sum_{l\geq 1}q^{l}[l]_{0}u^{l}.
\end{align*}
By Definition \ref{defn-b}, we can directly check  
\bal
&\sum_{i\geq 1}b_{2i}q^{i}u^{i}
=\frac{(q-1)u(q^{2}u+(q^{2}+q+1))}{(1+u)(1-q^{3}u)},\\
&\sum_{i\geq 1}b_{2i-1}q^{i}u^{i}
=\frac{qu(1-qu)}{(1+u)(1-q^{3}u)}.
\end{align*}
We can also directly check
\bal
&\sum_{k\geq 0}[2k+1]_{1}u^{k}=\frac{1+u}{(1-qu)(1-q^{-1}u)},\\
&\sum_{l\geq 1}q^{l}[l][2l+1]_{1}u^{l}
=\frac{u(q^{2}u+(q^{2}+q+1))}{(1-q^{3}u)(1-qu)(1-q^{-1}u)},\\
&\sum_{l\geq 1}q^{l}[l]_{0}u^{l}
=\frac{qu}{(1-q^{3}u)(1-q^{-1}u)}.
\end{align*}
Hence, the lemma follows.
\epf


We shall study the commutation relations 
between the real root vectors 
of height $1$ (Corollary \ref{c-2} (1))
and of height $-1$ (Corollary \ref{c-2} (2)). 

\blem\lab{le-1}
Let $n\geq{2}$. Then,
\bal
[E_{n\d +\a_{1}},E_{\a_{1}}]_{q}
&=[\ps_{1},[E_{(n-1)\d +\a_{1}},E_{\a_{1}}]_{q}]([3]_{1}!)^{-1}\\
&\ -T_{\pa}^{-1}([E_{(n-2)\d +\a_{1}},E_{\a_{1}}]_{q}).
\end{align*}
\elem

\bpf
By Corollary \ref{c-1} (7) and Lemma \ref{lem-bracket}, we have  
\bal
[E_{n\d +\a_{1}},E_{\a_{1}}]_{q}[3]_{1}!
&=[\ps_{1}E_{(n-1)\d +\a_{1}}-E_{(n-1)\d +\a_{1}}\ps_{1},E_{\a_{1}}]_{q}\\
&=\ps_{1}[E_{(n-1)\d +\a_{1}},E_{\a_{1}}]_{q}
+[\ps_{1},E_{\a_{1}}]E_{(n-1)\d +\a_{1}}q\\
&\ -E_{(n-1)\d +\a_{1}}[\ps_{1},E_{\a_{1}}]
-[E_{(n-1)\d +\a_{1}},E_{\a_{1}}]\ps_{1}\\
&=[\ps_{1},[E_{(n-1)\d +\a_{1}},E_{\a_{1}}]_{q}]
-[E_{(n-1)\d +\a_{1}},E_{\d +\a_{1}}]_{q}[3]_{1}!.
\end{align*}
Applying Lemma \ref{l-1}, we obtain the lemma.
\epf

\bprop\lab{p-1}
Let $n\geq{1}$. Then, 
\bal
[E_{n\d +\a_{1}},E_{\a_{1}}]_{q}&=
 (1+q)\sum_{i=1}^{[(n-1)/2]}b_{2i}E_{i\d +\a_{1}}E_{(n-i)\d +\a_{1}}\th (n\geq 3)\\
&\ +\begin{cases}
 \ b_{n}E_{\frac{n}{2}\d +\a_{1}}^{2}\ \ \ \text{if \ $n$\ is even},\\
 \ [4]_{1} b_{n}E_{(n+1)\d -\a_{0}}\ \ \ \text{if \ $n$\ is odd}.
 \end{cases}
 \end{align*}
\eprop

\bpf
We argue by the induction on $n$.
The case where $n=1$ is Lemma \ref{l-6} (5).
The case where $n=2$ follows 
from Lemma \ref{le-1} with $n=2$ and Lemma \ref{l-6} (5):
\begin{align*}
[E_{2\d +\a_{1}},E_{\a_{1}}]_{q}
&=[\ps_{1},[E_{\d +\a_{1}},E_{\a_{1}}]_{q}]([3]_{1}!)^{-1}
-(1-q)T_{\pa}^{-1}(E_{\a_{1}}^{2})\\
&=[\ps_{1},[4]_{1}E_{2\d -\a_{0}}]([3]_{1}!)^{-1}
+(q-1)E_{\d +\a_{1}}^{2}\\
&=(q-1)([4]_{1}[2]_{1}^{-1}+1)E_{\d +\a_{1}}^{2}
=b_{2}E_{\d +\a_{1}}^{2}.
\end{align*}
Now, assume that $n\geq{3}$. 
First, we consider the case  where $n=2m+1$ with $m\geq1$.
Using Lemma \ref{le-1}, the induction hypothesis, and Corollary \ref{c-1} (7), 
we have 
\bal
&[E_{(2m+1)\d +\a_{1}},E_{\a_{1}}]_{q}\\
&=[\ps_{1},[E_{2m\d +\a_{1}},E_{\a_{1}}]_{q}]([3]_{1}!)^{-1}
 -T_{\pa}^{-1}([E_{(2m-1)\d +\a_{1}},E_{\a_{1}}]_{q})\\
&=\Big[\ps_{1},(1+q)\sum_{i=1}^{m-1}b_{2i}E_{i\d +\a_{1}}
E_{(2m-i)\d +\a_{1}}+b_{2m}E_{m\d +\a_{1}}^{2}\Big]([3]_{1}!)^{-1}\\
&\ - T_{\pa}^{-1}\Big((1+q)\sum_{i=1}^{m-1}b_{2i}E_{i\d +\a_{1}}
E_{(2m-i-1)\d +\a_{1}}+[4]_{1}b_{2m-1}E_{2m\d -\a_{0}}\Big)\\
&=(1+q)\sum_{i=1}^{m-1}b_{2i}(E_{(i+1)\d +\a_{1}}E_{(2m-i)\d +\a_{1}}
+E_{i\d +\a_{1}}E_{(2m-i+1)\d +\a_{1}})\\
&\ +b_{2m}(E_{(m+1)\d +\a_{1}}E_{m\d +\a_{1}}+E_{m\d +\a_{1}}
E_{(m+1)\d +\a_{1}})\\
&\ -(1+q)\sum_{i=1}^{m-1}b_{2i}E_{(i+1)\d +\a_{1}}E_{(2m-i)\d +\a_{1}}
   -[4]_{1}b_{2m-1}E_{(2m+2)\d -\a_{0}}\\
&=(1+q)\sum_{i=1}^{m-1}b_{2i}E_{i\d +\a_{1}}E_{(2m-i+1)\d +\a_{1}}\\
&\ +b_{2m}((1+q)E_{m\d +\a_{1}}E_{(m+1)\d +\a_{1}}
 +[4]_{1}E_{(2m+2)\d -\a_{0}})\\
&\ -[4]_{1}b_{2m-1}E_{(2m+2)\d -\a_{0}}.
\end{align*}
By Corollary \ref{c-1} (5), this is equal to
\bal
&(1+q)\sum_{i=1}^{m}b_{2i}E_{i\d +\a_{1}}E_{(2m+1-i)\d +\a_{1}}
+[4]_{1}(b_{2m}-b_{2m-1})E_{(2m+2)\d -\a_{0}}.
\end{align*}
It remains to apply Lemma \ref{b-1} (1).
Next, we consider the case where $n=2m$ with $m\geq{2}$. 
Using Lemma \ref{le-1}, 
the induction hypothesis, and Corollary \ref{c-1} ((7), (10)), 
we have 
\bal
&[E_{2m\d +\a_{1}},E_{\a_{1}}]_{q}\\
&=[\ps_{1},[E_{(2m-1)\d +\a_{1}},E_{\a_{1}}]_{q}]([3]_{1}!)^{-1}
 -T_{\pa}^{-1}([E_{(2m-2)\d +\a_{1}},E_{\a_{1}}]_{q})\\
&=\Big[\ps_{1},(1+q)\sum_{i=1}^{m-1}b_{2i}E_{i\d +\a_{1}}E_{(2m-1-i)\d +\a_{1}}
+[4]_{1}b_{2m-1}E_{2m\d -\a_{0}}\Big]([3]_{1}!)^{-1}\\
&\ - T_{\pa}^{-1}\Big((1+q)\sum_{i=1}^{m-2}b_{2i}E_{i\d +\a_{1}}
E_{(2m-2-i)\d +\a_{1}}+b_{2m-2}E_{(m-1)\d +\a_{1}}^{2}\Big)\\
&=(1+q)\sum_{i=1}^{m-1}b_{2i}(E_{(i+1)\d +\a_{1}}E_{(2m-1-i)\d +\a_{1}}
+E_{i\d +\a_{1}}E_{(2m-i)\d +\a_{1}})\\
&\ +(q+q^{-1}) (q-1)b_{2m-1}E_{m\d +\a_{1}}^{2}\\
&\ -(1+q)\sum_{i=1}^{m-2}b_{2i}E_{(i+1)\d +\a_{1}}E_{(2m-1-i)\d +\a_{1}}
-b_{2m-2}E_{m\d +\a_{1}}^{2}\\
&=(1+q)\sum_{i=1}^{m-1}b_{2i}E_{i\d +\a_{1}}E_{(2m-i)\d +\a_{1}}\\
&\ +(q b_{2m-2}+(q+q^{-1})(q-1)b_{2m-1})E_{m\d +\a_{1}}^{2}.
\end{align*}
It remains to apply Lemma \ref{b-1} (2).
The proposition is proved.
\epf

\bcor\lab{c-2}
\bit
\item[(1)]
Let $k>l\geq{0}$. Then,
\bal
[E_{k\d +\a_{1}},E_{l\d +\a_{1}}]_{q}
&=(1+q)\sum_{i=1}^{[(k-l-1)/2]}b_{2i}E_{(l+i)\d +\a_{1}}E_{(k-i)\d +\a_{1}}
\th (k-l\geq 3)\\
&\ \ +\begin{cases}
 \ b_{k-l}E_{\frac{k+l}{2}\d +\a_{1}}^{2}\ \ \ \text{if \ $k-l$\ is even},\\
 \ [4]_{1} b_{k-l}E_{(k+l+1)\d -\a_{0}}\ \ \ \text{if \ $k-l$\ is odd}.\
 \end{cases}
\end{align*}
\item[(2)]
Let $l>k\geq{1}$. Then,
\begin{align*}
[E_{k\d -\a_{1}},E_{l\d -\a_{1}}]_{q}
&=(1+q)\sum_{i=1}^{[(l-k-1)/2]}b_{2i}E_{(l-i)\d -\a_{1}}E_{(k+i)\d -\a_{1}}
\th (l-k\geq 3)\\
&\ \ +\begin{cases}
 \ b_{l-k}E_{\frac{k+l}{2}\d -\a_{1}}^{2}\ \ \ \text{if \ $l-k$\ is even},\\
 \ [4]_{1} b_{l-k}E_{(k+l-1)\d +\a_{0}}\ \ \ \text{if \ $l-k$\ is odd}.\
 \end{cases}
 \end{align*}
\eit
\ecor
\bpf
Applying $T_{\pa}^{-(k-l)}$ to Proposition \ref{p-1}, we obtain (1).
Applying $T_{1}^{-1}\ast$ to (1) with $k$ and $l$ exchanged, we obtain (2).
\epf


We shall prove that the imaginary root vectors are invariant 
under $T_{\pa}$ and $T_{1}^{-1}\ast$ (Proposition \ref{p-2} ((1), (3)))
and are mutually commutative (Proposition \ref{p-3}).
As by-products,
we shall obtain the commutation relations 
between the imaginary root vectors
and the real root vectors 
of height $1$ or  $-1$ (Proposition \ref{p-2} ((4), (5))), 
and the ones 
between the real root vectors 
of height $1$ and $-1$ (Proposition \ref{p-2} (2)),
of height $1$ and $-2$ (Corollary \ref{c-3} (1)),
and of height $2$ and $-1$ (Corollary \ref{c-3} (2)); 
as for the ones between the real root vectors of height $2$ and $-2$, 
see Appendix A.
First, we prove several lemmas.

\blem\lab{le-2}
Let $n\geq 1$. Then,
\bal
[E_{\a_{0}},E_{n\d +\a_{1}}]_{q^{-2}}[3]_{1}!&
=(q_{1}-q_{1}^{-1})[3]_{1}!\ps_{n} E_{\d -\a_{1}}
+(q-1)[3]_{1}[E_{\d -\a_{1}},\ps_{n}]\\
&\ -[[E_{\a_{0}},E_{(n-1)\d +\a_{1}}]_{q^{-2}},\ps_{1}].
\end{align*}
\elem

\bpf
Using Corollary \ref{c-1} (7), Lemma \ref{lem-bracket}, 
and Lemma \ref{l-6} (4), we have 
\bal
&[E_{\a_{0}},E_{n\d +\a_{1}}]_{q^{-2}}[3]_{1}!\\
&=[E_{\a_{0}},\ps_{1}E_{(n-1)\d +\a_{1}}-E_{(n-1)\d +\a_{1}}\ps_{1}]_{q^{-2}}\\
&=[E_{\a_{0}},\ps_{1}]E_{(n-1)\d +\a_{1}}
+\ps_{1}[E_{\a_{0}},E_{(n-1)\d +\a_{1}}]_{q^{-2}}\\
&-[E_{\a_{0}},E_{(n-1)\d +\a_{1}}]_{q^{-2}}\ps_{1}
-E_{(n-1)\d +\a_{1}}[E_{\a_{0}},\ps_{1}]q^{-2}\\
&=[\ps_{1},[E_{\a_{0}},E_{(n-1)\d +\a_{1}}]_{q^{-2}}]
+(q-1)[3]_{1}[E_{\d -\a_{1}}^{2},E_{(n-1)\d +\a_{1}}]_{q^{-2}}.
\end{align*}
We rewrite the last term: using Definition \ref{defn-ps} and 
Lemma \ref{lem-bracket}, we have 
\bal
&[E_{\d -\a_{1}}^{2},E_{(n-1)\d +\a_{1}}]_{q^{-2}}\\
&=E_{\d -\a_{1}}[E_{\d -\a_{1}},E_{(n-1)\d +\a_{1}}]_{q^{-1}}
+[E_{\d -\a_{1}},E_{(n-1)\d +\a_{1}}]_{q^{-1}}E_{\d -\a_{1}}q^{-1}\\
&=E_{\d -\a_{1}}\ps_{n}+\ps_{n}E_{\d -\a_{1}}q^{-1}\\
&=(1+q^{-1})\ps_{n}E_{\d -\a_{1}}+[E_{\d -\a_{1}},\ps_{n}].
\end{align*}
Hence, we obtain the lemma.
\epf

\blem\lab{le-3}
Let $n\geq 2$. Then,
\bal
T_{\pa}(\ps_{n})=\ps_{n}+[\ps_{n-1},\ps_{1}]([3]_{1}!)^{-1}.
\end{align*}
\elem

\bpf  Using Definition \ref{defn-ps}, Corollary \ref{c-1} ((7), (8)), 
and Lemma \ref{lem-bracket}, we have 
\bal
\ps_{n}[3]_{1}!
&=[E_{\d -\a_{1}},E_{(n-1)\d +\a_{1}}]_{q^{-1}}[3]_{1}!\\
&=[E_{\d -\a_{1}},\ps_{1}E_{(n-2)\d +\a_{1}}
-q^{-1}E_{(n-2)\d +\a_{1}}\ps_{1}]_{q^{-1}}\\
&=[E_{\d -\a_{1}},\ps_{1}]E_{(n-2)\d +\a_{1}}
+\ps_{1}[E_{\d -\a_{1}},E_{(n-2)\d +\a_{1}}]_{q^{-1}}\\
&\ -[E_{\d -\a_{1}},E_{(n-2)\d +\a_{1}}]_{q^{-1}}\ps_{1}
-E_{(n-2)\d +\a_{1}}[E_{\d -\a_{1}},\ps_{1}]q^{-1}\\
&=[\ps_{1},[E_{\d -\a_{1}},E_{(n-2)\d +\a_{1}}]_{q^{-1}}]
+[E_{2\d -\a_{1}},E_{(n-2)\d +\a_{1}}]_{q^{-1}}[3]_{1}!\\
&=[\ps_{1},\ps_{n-1}]+T_{\pa}(\ps_{n})[3]_{1}!.
\end{align*}
Hence, we obtain the lemma.
\epf

\blem\lab{le-4}
Let $n\geq 1$. Then,
\bal
[\ps_{n},E_{\a_{1}}]&
=(1+q^{-1})qb_{2}E_{n\d +\a_{1}}\ps_{0}
+(1+q^{-1})qb_{0}E_{(n-1)\d +\a_{1}}\ps_{1}\\
&\ +[E_{\d -\a_{1}},[E_{(n-1)\d +\a_{1}},E_{\a_{1}}]_{q}]_{q^{-2}}.
\end{align*}
\elem

\bpf
Using Definition \ref{defn-ps}, Corollary \ref{c-1} (7), 
and Lemma \ref{lem-bracket}, 
we see that the left hand side is equal to
\bal
&[E_{\d -\a_{1}}E_{(n-1)\d +\a_{1}}
-q^{-1}E_{(n-1)\d +\a_{1}}E_{\d -\a_{1}},E_{\a_{1}}]\\
&=E_{\d -\a_{1}}[E_{(n-1)\d +\a_{1}},E_{\a_{1}}]_{q}
+[E_{\d -\a_{1}},E_{\a_{1}}]_{q^{-1}}E_{(n-1)\d +\a_{1}}q\\
&\ -E_{(n-1)\d +\a_{1}}[E_{\d -\a_{1}},E_{\a_{1}}]_{q^{-1}}q^{-1}
-[E_{(n-1)\d +\a_{1}},E_{\a_{1}}]_{q}E_{\d -\a_{1}}q^{-2}\\
&=[E_{\d -\a_{1}},[E_{(n-1)\d +\a_{1}},E_{\a_{1}}]_{q}]_{q^{-2}}
+(q-q^{-1})E_{(n-1)\d +\a_{1}}\ps_{1}+q[3]_{1}!E_{n\d +\a_{1}}.
\end{align*}
Hence, we obtain the lemma.
\epf

\blem\lab{le-5}
Let $k,l\geq 0$. Then,
\bal
[E_{\d -\a_{1}},E_{k\d +\a_{1}}E_{l\d +\a_{1}}]_{q^{-2}}
=q^{-1}E_{k\d +\a_{1}}\ps_{l+1} + E_{l\d +\a_{1}}\ps_{k+1} 
+[\ps_{k+1},E_{l\d +\a_{1}}].
\end{align*}
\elem

\bpf
By Lemma \ref{lem-bracket} and Definition \ref{defn-ps}, 
the left hand side is equal to
\bal
&[E_{\d -\a_{1}},E_{k\d +\a_{1}}]_{q^{-1}}E_{l\d +\a_{1}}
+E_{k\d +\a_{1}}[E_{\d -\a_{1}},E_{l\d +\a_{1}}]_{q^{-1}}q^{-1}\\
&=q^{-1}E_{k\d +\a_{1}}\ps_{l+1}+\ps_{k+1}E_{l\d +\a_{1}}.
\end{align*}
Hence, we obtain the lemma.
\epf

\blem\lab{le-6} 
We have 
\bit
\item[(1)]
$r_{0}(E_{\a})=0$\ for $\a\in R_{re}^{+}\backslash \{ \a_{0} \}$,
\item[(2)]
$r_{0}(\ps_{n})=0$\ for $n\geq 1$.
\eit
\elem

\bpf
Using Lemma \ref{l-3} (1), we have 
$r_{0}(E_{\d -\a_{1}})=0$.
Hence, using Definition \ref{defn-ps}, we have $r_{0}(\ps_{1})=0$.
Then, applying Corollary \ref{c-1}, we obtain (1). 
Using (1) and Definition \ref{defn-ps}, we obtain (2).
\epf

\bprop\lab{p-2}
Let $n\geq{1}$. Then,
\bit
\item[(1)]
$T_{\pa}(\ps_{n})=\ps_{n}$,
\item[(2)]
$[E_{k\d -\a_{1}},E_{l\d +\a_{1}}]_{q^{-1}}=\ps_{n}$\ \   
for $k\geq 1,\ l\geq 0$ such that $k+l=n$,
\item[(3)]
$(T_{1}^{-1}\ast )(\ps_{n})=\ps_{n}$,
\item[(4)]
$[\ps_{n},E_{k\d +\a_{1}}]=(1+q^{-1})\sum_{i=0}^{n-1}b_{2(n-i)}
E_{(k+n-i)\d +\a_{1}}\ps_{i}$\ \ for $k\geq 0$,
\item[(5)]
$[E_{k\d -\a_{1}},\ps_{n}]=(1+q^{-1})\sum_{i=0}^{n-1}b_{2(n-i)}
\ps_{i}E_{(k+n-i)\d -\a_{1}}$\ \ for $k\geq 1$,
\item[(6)]
$[\ps_{n},\ps_{1}]=0$,
\item[(7)]
$[E_{\a_{0}},E_{n\d +\a_{1}}]_{q^{-2}}=(q_{1}-q_{1}^{-1})
\sum_{i=0}^{n}b_{2(n-i)+1}\ps_{i}E_{(n-i+1)\d -\a_{1}}$.
\eit
\eprop

\bpf
We denote by (a)$_{r}$ the statement (a) for $n=1,\ldots,r$.
We prove (1)$_{n}$--(7)$_{n}$ at once by the induction on $n$.
If $n=1$, we have (1) and (3) by Lemma \ref{l-5}, 
(2) by Definition \ref{defn-ps},
(4) and (5) by Corollary \ref{c-1} ((7),\ (8)),
(6) trivially, and (7) by Lemma \ref{le-2} with $n=1$.
Now, assume that $n\geq 2$.
 
(1)$_{n}$\ \ 
By Lemma \ref{le-3} and $(6)_{n-1}$, we obtain (1) for $n$.

(2)$_{n}$\ \ 
By $(1)_{n}$ and by applying $T_{\pa}^{k-1}$ to Definition \ref{defn-ps},  
we obtain (2) for $n$.

(3)$_{n}$\ \ 
By Definition \ref{defn-ps}, we have 
\bal
(T_{1}^{-1}\ast )(\ps_{n})&
=(T_{1}^{-1}\ast )([E_{\d -\a_{1}},E_{(n-1)\d +\a_{1}}]_{q^{-1}})
=[E_{(n-1)\d -\a_{1}},E_{\d +\a_{1}}]_{q^{-1}}.
\end{align*}
Applying $(2)_{n}$, we obtain (3) for $n$.

(4)$_{n}$\ \ 
By  $(1)_{n}$, we can assume that $k=0$.
First, we consider the case where $n=2m$ with $m\geq 1$.
Note that we have $2m > m$. 
Using Proposition \ref{p-1} and Lemma \ref{le-5}, we have 
\bal
&[E_{\d -\a_{1}},[E_{(2m-1)\d +\a_{1}},E_{\a_{1}}]_{q}]_{q^{-2}}\\
&=[E_{\d -\a_{1}},(1+q)\sum_{i=1}^{m-1}b_{2i}E_{i\d +\a_{1}}
E_{(2m-i-1)\d +\a_{1}}\th(m\geq 2)+[4]_{1}b_{2m-1}E_{2m\d -\a_{0}}]_{q^{-2}}\\
&=(1+q)\sum_{i=1}^{m-1}b_{2i}\Big( 
q^{-1}E_{i\d +\a_{1}}\ps_{2m-i}+E_{(2m-i-1)\d +\a_{1}}\ps_{i+1}\\
&\ +[\ps_{i+1},E_{(2m-i-1)\d +\a_{1}}]\Big)\th(m\geq 2)
+[4]_{1}b_{2m-1}(T_{1}^{-1}\ast T_{\pa}^{m-1})
([E_{\a_{0}},E_{m\d +\a_{1}}]_{q^{-2}}).
\end{align*}
By $(4)_{m}$ and $(7)_{m}$, this is equal to
\bal
&(1+q^{-1})\sum_{i=m+1}^{2m-1}b_{2(2m-i)}E_{(2m-i)\d +\a_{1}}\ps_{i}
\th(m\geq 2)\\
&\ +(1+q^{-1})\sum_{i=1}^{m-1}qb_{2(i-1)}E_{(2m-i)\d +\a_{1}}\ps_{i}
\th(m\geq 2)\\
&\ +(1+q^{-1})(1+q)\sum_{i=1}^{m-1}\sum_{j=0}^{i}b_{2i}b_{2(i-j+1)}
E_{(2m-j)\d +\a_{1}}\ps_{j}\th(m\geq 2)\\
&\ +[4]_{1}b_{2m-1}(q_{1}-q_{1}^{-1})(T_{1}^{-1}\ast T_{\pa}^{m-1})
\Big(\sum_{i=0}^{m}b_{2(m-i)+1}\ps_{i}E_{(m-i+1)\d -\a_{1}}\Big).
\end{align*}
In the third term, putting $j=i',\ i-j+1=j'$,
we see that $i=i'+j'-1$ and that 
$1\leq i \leq m-1,\ 0\leq j \leq i$ is equivalent to 
$0\leq i' \leq m-1,\ \max(1,2-i')\leq j'\leq m-i'$. 
We can apply $(1)_{m}$ and $(3)_{m}$ to the fourth one. 
Thus, by Lemma \ref{le-4}, we have
\bal
[\ps_{2m},E_{\a_{1}}]
&=(1+q^{-1})qb_{2}E_{2m\d +\a_{1}}\ps_{0}
+(1+q^{-1})qb_{0}E_{(2m-1)\d +\a_{1}}\ps_{1}\\
&\ +[E_{\d -\a_{1}},[E_{(2m-1)\d +\a_{1}},E_{\a_{1}}]_{q}]_{q^{-2}}\\
&=(1+q^{-1})\sum_{i=m+1}^{2m-1}b_{2(2m-i)}E_{(2m-i)\d +\a_{1}}\ps_{i}
\th(m\geq 2)\\
&\ +(1+q^{-1})\sum_{i=0}^{m}\Big( qb_{2}\th(i=0)
+qb_{2(i-1)}\th(i\geq 1)\\
&\ +(1+q)\sum_{j=\max(1,2-i)}^{m-i}b_{2j}b_{2(j+i-1)}
\th(m\geq 2)\th(i\leq m-1)\\
&\ +(q+q^{-1})(q-1)b_{2m-1}b_{2(m-i+1)-1}\Big) 
E_{i\d +\a_{1}}\ps_{2m-i}.
\end{align*}
Applying Lemma \ref{b-3} (2) with $l=m-i+1,\ s=i-1$ together with 
Lemma \ref{b-1} (2), we obtain the desired result.
Next, we consider the case where $n=2m+1$ with $m\geq 1$.
Note that we have $2m+1 > m+1$. 
Using Proposition \ref{p-1}, Lemma \ref{le-5}, and $(4)_{m+1}$, 
we have 
\bal
&[E_{\d -\a_{1}},[E_{2m\d +\a_{1}},E_{\a_{1}}]_{q}]_{q^{-2}}\\
&=\Big[E_{\d -\a_{1}},(1+q)\sum_{i=1}^{m-1}b_{2i}E_{i\d +\a_{1}}
E_{(2m-i)\d +\a_{1}}\th(m\geq 2)+b_{2m}E_{m\d +\a_{1}}^{2}\Big]_{q^{-2}}\\
&=(1+q)\sum_{i=1}^{m-1}b_{2i}\Big( q^{-1}E_{i\d +\a_{1}}\ps_{2m-i+1}
+E_{(2m-i)\d +\a_{1}}\ps_{i+1}\\
&\ +[\ps_{i+1},E_{(2m-i)\d +\a_{1}}]\Big)\th(m\geq 2)\\
&\ +b_{2m}((1+q^{-1})E_{m\d +\a_{1}}\ps_{m+1}+[\ps_{m+1},E_{m\d +\a_{1}}])\\
&=(1+q^{-1})\sum_{i=m+1}^{2m}b_{2(2m-i+1)}E_{(2m-i+1)\d +\a_{1}}\ps_{i}\\
&\ +(1+q^{-1})\sum_{i=2}^{m}qb_{2(i-1)}E_{(2m-i+1)\d +\a_{1}}\ps_{i}
\th(m\geq 2)\\
&\ +(1+q^{-1})(1+q)\sum_{i=1}^{m-1}\sum_{j=0}^{i}b_{2i}b_{2(i-j+1)}
E_{(2m-j+1)\d +\a_{1}}\ps_{j}\th(m\geq 2)\\
&\ +(1+q^{-1})b_{2m}\sum_{i=0}^{m}b_{2(m-i+1)}E_{(2m-i+1)\d +\a_{1}}\ps_{i}.
\end{align*}
In the third term, putting $j=i',\ i-j+1=j'$, 
we see that $i=i'+j'-1$
and that $1\leq i \leq m-1,\ 0\leq j \leq i$ 
is equivalent to $0\leq i' \leq m-1,\ \max(1,2-i')\leq j'\leq m-i'$.
Thus, by Lemma \ref{le-4}, we have 
\bal
&[\ps_{2m+1},E_{\a_{1}}]\\
&=(1+q^{-1})qb_{2}\th(i=0)+(1+q^{-1})qb_{2(i-1)}\th(i\geq 1)\\
&\ +[E_{\d -\a_{1}},[E_{2m\d +\a_{1}},E_{\a_{1}}]_{q}]_{q^{-2}}\\
&=(1+q^{-1})\sum_{i=m+1}^{2m}b_{2(2m-i+1)}E_{(2m-i+1)\d +\a_{1}}\ps_{i}\\
&\ +(1+q^{-1})\sum_{i=0}^{m}E_{(2m-i+1)\d +\a_{1}}\ps_{i}
\Big(qb_{2}\th(i=0)+qb_{2(i-1)}\th(i\geq 1)\\
&\ +(1+q)\sum_{j=\max(1,2-i)}^{m-i}b_{2j}b_{2(i+j-1)}
\th(i\geq m-1)\th(m\geq 2)+b_{2m}b_{2(m-i+1)}\Big).
\end{align*}
Applying Lemma \ref{b-3} (1) with $l=m-i+1,\ s=i-1$ together with 
Lemma \ref{b-1} (1), we obtain the desired result.
We obtain (4) for $n$.

(5)$_{n}$\ \ 
By $(3)_{n}$ and by applying $T_{1}^{-1}\ast$ to $(4)_{n}$ 
with $k\geq 1$, we obtain (5) for $n$.

(6)$_{n}$\ \ 
By Lemma \ref{r-1} (1) and Lemma \ref{le-6} (2), 
it is enough to show that $r_{1}([\ps_{n},\ps_{1}])=0$. 
Using Lemma \ref{l-4} and $(1)_{n-1}$,
we have $[E_{(n-1)\d +\a_{1}},f_{1}]=k_{1}\ps_{n-1}$. 
Applying Lemma \ref{r-2} and Proposition \ref{prop-triangular}, we have 
\bal
(_{1}r(E_{(n-1)\d +\a_{1}})=0\ \ \text{and})\ \ 
r_{1}(E_{(n-1)\d +\a_{1}})=(q_{1}-q_{1}^{-1})\ps_{n-1}.
\end{align*}
Similarly, using Lemma \ref{l-3} ((4), (5)), Lemma \ref{r-2}, 
and Proposition \ref{prop-triangular}, we have 
\bal
(_{1}r(E_{\d -\a_{1}})=0\ \ \text{and})\ \ 
r_{1}(E_{\d -\a_{1}})=(q_{1}-q_{1}^{-1})q^{-2}[4]_{1}E_{\a_{0}},\\
(_{1}r(\ps_{1})=0\ \ \text{and})\ \ 
r_{1}(\ps_{1})=(q_{1}-q_{1}^{-1})q^{-1}[3]_{1}!E_{\d -\a_{1}}.
\end{align*}
Hence, using the definition of $\ps_{n}$, we have ($_{1}r(\ps_{n})=0$ and)
\bal
r_{1}(\ps_{n})
&=r_{1}(E_{\d -\a_{1}}E_{(n-1)\d +\a_{1}}-q^{-1}E_{(n-1)\d +\a_{1}}
E_{\d -\a_{1}})\\
&=(q_{1}-q_{1}^{-1})(q^{-1}[4]_{1}E_{\a_{0}}E_{(n-1)\d +\a_{1}}
+E_{\d -\a_{1}}\ps_{n-1})\\
&\ -q^{-1}(q_{1}-q_{1}^{-1})(q^{-1}\ps_{n-1}E_{\d -\a_{1}}
+E_{(n-1)\d +\a_{1}}q^{-2}[4]_{1}E_{\a_{0}})\\
&=(q_{1}-q_{1}^{-1})
\Big(q^{-1}[4]_{1}[E_{\a_{0}},E_{(n-1)\d +\a_{1}}]_{q^{-2}}\\
&\ +(1-q^{-2})\ps_{n-1}E_{\d -\a_{1}}+[E_{\d -\a_{1}},\ps_{n-1}]\Big).
\end{align*}
Applying $(5)_{n-1}$ and $(7)_{n-1}$, we have
\bal
&r_{1}(\ps_{n})(q_{1}-q_{1}^{-1})^{-1}\\
&=q^{-1}[4]_{1}(q_{1}-q_{1}^{-1})\sum_{i=0}^{n-1}
b_{2(n-i)-1}\ps_{i}E_{(n-i)\d -\a_{1}}\\
&\ +(1-q^{-2})\ps_{n-1}E_{\d -\a_{1}}
+(1+q^{-1})\sum_{i=0}^{n-2}b_{2i}\ps_{i}E_{(n-i)\d -\a_{1}}\\
&=(1+q^{-1})q^{-1}
\sum_{i=0}^{n-1}\big( qb_{2(n-i-1)}+(q+q^{-1})(q-1)b_{2(n-i)-1}\big)
\ps_{i}E_{(n-i)\d -\a_{1}}.
\end{align*}
Applying  Lemma \ref{b-1} (2), we have 
\bal
r_{1}(\ps_{n})=(q_{1}-q_{1}^{-1})q^{-1}(1+q^{-1})
\sum_{i=0}^{n-1}b_{2(n-i)}\ps_{i}E_{(n-i)\d -\a_{1}}.
\end{align*}
Hence,
\bal
&r_{1}([\ps_{n},\ps_{1}])(q_{1}-q_{1}^{-1})^{-1}q\\
&=([r_{1}(\ps_{n}),\ps_{1}]+[\ps_{n},r_{1}(\ps_{1})])(q_{1}-q_{1}^{-1})^{-1}q\\
&=(1+q^{-1})
\Big[\sum_{i=0}^{n-1}b_{2(n-i)}\ps_{i}E_{(n-i)\d -\a_{1}},\ps_{1}\Big]
-[3]_{1}![E_{\d -\a_{1}},\ps_{n}].
\end{align*}
It follows from $(6)_{n-1}$ and Corollary \ref{c-1} (7) that
the former term is equal to 
\bal
(1+q^{-1})\sum_{i=1}^{n}b_{2i}\ps_{n-i}E_{(i+1)\d -\a_{1}}[3]_{1}!,
\end{align*}
which cancels out with the latter by $(5)_{n}$. 
We obtain (6) for $n$.

(7)$_{n}$\ \ 
Using Lemma \ref{le-2}, $(5)_{n}$, and $(7)_{n-1}$, we have 
\bal
[E_{\a_{0}},E_{n\d +\a_{1}}]_{q^{-2}}[3]_{1}!
&=(q_{1}-q_{1}^{-1})[3]_{1}!\ps_{n}E_{\d -\a_{1}}\\
&\ +(q-1)[3]_{1}(1+q^{-1})
\sum_{i=0}^{n-1}b_{2(n-i)}\ps_{i}E_{(n-i+1)\d -\a_{1}}\\
&\ -(q_{1}-q_{1}^{-1})\Big[\sum_{i=0}^{n-1}b_{2(n-i)-1}\ps_{i}
E_{(n-i)\d -\a_{1}},\ps_{1}\Big].
\end{align*}
Applying $(6)_{n-1}$ and Corollary \ref{c-1} (8), we have 
\bal
[E_{\a_{0}},E_{n\d +\a_{1}}]_{q^{-2}}
&=(q_{1}-q_{1}^{-1})\ps_{n}E_{\d -\a_{1}}\\
&\ +(q_{1}-q_{1}^{-1})\sum_{i=0}^{n-1}
(b_{2(n-i)}-b_{2(n-i)-1})\ps_{i}E_{(n-i+1)\d -\a_{1}}.
\end{align*}
Applying Lemma \ref{b-1} (1), we obtain (7) for $n$.

The proposition is proved.
\epf

\bcor\lab{c-3}
\bit
\item[(1)]
Let $m,n\geq 0$. Then,
\bal
[E_{2m\d +\a_{0}},E_{n\d +\a_{1}}]_{q^{-2}}
&=(q_{1}-q_{1}^{-1})
\sum_{i=0}^{m+n}b_{2(m-i)+1}\ps_{i}E_{(2m+n-i+1)\d -\a_{1}}.
\end{align*}
\item[(2)]
Let $m,n\geq 1$. Then,
\bal
[E_{n\d -\a_{1}},E_{2m\d -\a_{0}}]_{q^{-2}}
&=(q_{1}-q_{1}^{-1})
\sum_{i=0}^{m+n-1}b_{2(m-i)-1}E_{(2m+n-i-1)\d +\a_{1}}\ps_{i}.
\end{align*}
\eit
\ecor

\bpf
If $m=n=0$, then (1) is Lemma \ref{l-3} (1). 
Otherwise, applying $T_{\pa}^{m}$ to Proposition \ref{p-1} (6) 
with $n$ replaced by $m+n$, we obtain (1). 
Applying $T_{1}^{-1}\ast$ to (1) 
with $m$ replaced by $m-1$ and with $n\geq 1$, we obtain (2).
\epf

\bcor\lab{c-4}
We have 
\bit
\item[(1)]
$_{1}r(E_{\a})=0$ for $\a\in R_{re}^{+}\backslash \{\a_{1}\}$,
\item[(2)]
$_{1}r(\ps_{n})=0$ for $n\geq 1$,
\item[(3)]
$r_{1}(\ps_{n})=(q_{1}-q_{1}^{-1})q^{-1}(1+q^{-1})
\sum_{i=1}^{n}b_{2i}\ps_{n-i}E_{i\d -\a_{1}}$ for $n\geq 1$,
\item[(4)]
$r_{1}(E_{n\d +\a_{1}})=(q_{1}-q_{1}^{-1})\ps_{n}$ for $n\geq 0$,
\item[(5)]
$r_{1}(E_{\d -\a_{1}})=(q_{1}-q_{1}^{-1})q^{-2}[4]_{1}E_{\a_{0}}$.
\eit
\ecor
\bpf 
In the proof of the induction step for Proposition \ref{p-1} (6), 
we already have 
(1) for $\a= n\d+\a_{1}$ with $n\geq 0$ and for $\a=\d-\a_{1}$, and (2)--(5).
By Corollary \ref{c-1}, we see that it remains to check 
(1) for $\a=2\d-\a_{0}$, but it follows from Lemma \ref{l-6} (5).
\epf

\bprop\lab{p-3}
Let $k,l\geq{0}$. Then, $[\ps_{k},\ps_{l}]=0$. 
\eprop

\bpf
We can assume that $1\leq k\leq l$.
We argue by the induction on $k$.
By Definition \ref{defn-ps} and 
Proposition \ref{p-2} ((4), (5)), we have 
\bal
[\ps_{k},\ps_{l}]
&=[\ps_{k},[E_{\d -\a_{1}},E_{(l-1)\d +\a_{1}}]_{q^{-1}}]\\
&=[[\ps_{k},E_{\d -\a_{1}}],E_{(l-1)\d +\a_{1}}]_{q^{-1}}
+[E_{\d -\a_{1}},[\ps_{k},E_{(l-1)\d +\a_{1}}]]_{q^{-1}}\\
&=-(1+q^{-1})\Big[\sum_{i=0}^{k-1}b_{2(k-i)}\ps_{i}E_{(k-i+1)\d -\a_{1}}
,E_{(l-1)\d +\a_{1}}\Big]_{q^{-1}}\\
&\ +(1+q^{-1})\Big[E_{\d -\a_{1}},\sum_{i=0}^{k-1}b_{2(k-i)}
E_{(k+l-i-1)\d +\a_{1}}\ps_{i}\Big]_{q^{-1}}.
\end{align*}
Applying Proposition \ref{p-2} (2), we have 
\bal
&[\ps_{k},\ps_{l}](1+q^{-1})^{-1}\\
&=-\sum_{i=0}^{k-1}b_{2(k-i)}\Big(\ps_{i}
(q^{-1}E_{(l-1)\d +\a_{1}}E_{(k-i+1)\d -\a_{1}}+\ps_{k+l-i})\\
&\ -q^{-1}E_{(l-1)\d +\a_{1}}\ps_{i}E_{(k-i+1)\d -\a_{1}}\Big)\\
&\ +\sum_{i=0}^{k-1}b_{2(k-i)}
\Big( (q^{-1}E_{(k+l-i-1)\d +\a_{1}}E_{\d -\a_{1}}+\ps_{k+l-i})\ps_{i}\\
&\ -q^{-1}E_{(k+l-i-1)\d +\a_{1}}\ps_{i}E_{\d -\a_{1}}\Big)\\
&=-q^{-1}\sum_{i=1}^{k-1}
b_{2(k-i)}[\ps_{i},E_{(l-1)\d +\a_{1}}]E_{(k-i+1)\d -\a_{1}}\th(k\geq 2)\\
&\ +q^{-1}\sum_{i=1}^{k-1}b_{2(k-i)}
E_{(k+l-i-1)\d +\a_{1}}[E_{\d -\a_{1}},\ps_{i}]\th(k\geq 2)\\
&\ +\sum_{i=0}^{k-1}b_{2(k-i)}[\ps_{k+l-i},\ps_{i}].
\end{align*}
By the induction hypothesis, the last term vanishes. 
We can assume that $k\geq 2$.
Applying Proposition \ref{p-2} ((4), (5)), we have 
\bal
[\ps_{k},\ps_{l}](1+q^{-1})^{-2}q
&=-\sum_{i=1}^{k-1}\sum_{j=0}^{i-1}
b_{2(k-i)}b_{2(i-j)}E_{(l+i-j-1)\d +\a_{1}}\ps_{j}E_{(k-i+1)\d -\a_{1}}\\
&\ +\sum_{i=1}^{k-1}\sum_{j=0}^{i-1}
b_{2(k-i)}b_{2(i-j)}E_{(k+l-i-1)\d +\a_{1}}\ps_{j}E_{(i-j+1)\d -\a_{1}}.
\end{align*}
In the former term, putting $k-i=i'-j$, we see that $i=k+j-i'$ and that 
$1\leq i\leq k-1,\ 0\leq j\leq i-1$ is equivalent to 
$1\leq i'\leq k-1,\ 0\leq j\leq i'-1$. Thus, 
it cancels out with the latter 
and the induction proceeds.
The proposition is proved.
\epf


\section{Construction of Basis I}

\bdefn
In the following, each $\Z_{\geq 0}^{(i)}$ is a copy of $\Z_{\geq 0}$.
\bit
\item[(1)]
For ${\bf c}=({\bf c}_{i})\in\oplus_{i\in R_{re}^{+}(>)}\Z_{\geq 0}^{(i)}$, 
we set 
${\bf E_{c}}=E_{\a_{1}}^{({\bf c}_{\a_{1}})}
E_{2\d-\a_{0}}^{({\bf c}_{2\d-\a_{0}})}
E_{\d+\a_{1}}^{({\bf c}_{\d+\a_{1}})}\cdots$.
\item[(2)]
For ${\bf c}=({\bf c}_{i})\in\oplus_{i\in R_{re}^{+}(<)}\Z_{\geq 0}^{(i)}$, 
we set
${\bf E_{c}}=\cdots E_{2\d+\a_{0}}^{({\bf c}_{2\d+\a_{0}})}
E_{\d-\a_{1}}^{({\bf c}_{\d-\a_{1}})} E_{\a_{0}}^{({\bf c}_{\a_{0}})}$.
\item[(3)]
For ${\bf c}=({\bf c}_{i})\in\oplus_{i\in\Z_{\geq 1}}\Z_{\geq 0}^{(i)}$, 
we set
${\bf E_{c}'}=\ps_{1}^{{\bf c}_{1}}\ps_{2}^{{\bf c}_{2}}
\ps_{3}^{{\bf c}_{3}}\cdots$.
\item[(4)]
For ${\bf c}=({\bf c}_{i})\in\oplus_{i\in\Z_{\geq 1}}\Z_{\geq 0}^{(i)}$, 
we set
${\bf E_{c}}=P_{1}^{{\bf c}_{1}}P_{2}^{{\bf c}_{2}}P_{3}^{{\bf c}_{3}}\cdots$.
\eit
\edefn

\bdefn
Let $\a,\a'\in R_{re}^{+}(>)$ and $\a\leq\a'$.
Let $\b,\b'\in R_{re}^{+}(<)$ and $\b\leq\b'$. We set 
\bit
\item[(1)]
$B(>;\a,\a')=
\{ {\bf E_{c}}|\ {\bf c}\in\oplus_{i\in R_{re}^{+}(>)}\Z_{\geq 0}^{(i)},
\ {\bf c}_{i}=0\ \ \text{if}\ \ i<\a \ \text{or}\ i>\a'\}$,
\item[(2)]
$B(>;\a)=
\{ {\bf E_{c}}|\ {\bf c}\in\oplus_{i\in R_{re}^{+}(>)}\Z_{\geq 0}^{(i)},
\ {\bf c}_{i}=0\ \ \text{if}\ \ i<\a \}$,
\item[(3)]
$B(<;\b,\b')=
\{ {\bf E_{c}}|\ {\bf c}\in\oplus_{i\in R_{re}^{+}(<)}\Z_{\geq 0}^{(i)},
\ {\bf c}_{i}=0\ \ \text{if}\ \ i<\b \ \text{or}\ i>\b'\}$,
\item[(4)]
$B(<;\b)=
\{ {\bf E_{c}}|\ {\bf c}\in\oplus_{i\in R_{re}^{+}(<)}\Z_{\geq 0}^{(i)},
\  {\bf c}_{i}=0\ \ \text{if}\ \ i>\b \}$.
\eit
Then, each set is contained in $\U_{\Z}^{+}$ and 
linearly independent over $\Q(q_{1})$ \cite[40.2.1]{L}. 
We also set $B(>)=B(>;\a_{1})$ and $B(<)=B(<;\a_{0})$.
\edefn

\bdefn
Let $\a,\a'\in R_{re}^{+}(>)$ and $\a\leq\a'$.
\bit
\item[(1)]
Let $\U^{+}(>;\a,\a')$ be the $\Q(q_{1})$-subspace of $\U^{+}$
with basis $B(>;\a,\a')$. 
\item[(2)]
Let $\U^{+}(>;\a)$ be the $\Q(q_{1})$-subspace of $\U^{+}$
with basis $B(>;\a)$. We also set $\U^{+}(>)=\U^{+}(>;\a_{1})$
\item[(3)]
Let $\U^{+}_{\Z}(>;\a,\a')$ be the free $\Z [q_{1},q_{1}^{-1}]$-submodule 
of $\U_{\Z}^{+}$ with basis $B(>;\a,\a')$.
\item[(4)]
Let $\U^{+}_{\Z}(>;\a)$ be 
the free $\Z [q_{1},q_{1}^{-1}]$-submodule of $\U_{\Z}^{+})$ 
with basis $B(>;\a)$. We also set $\U_{\Z}^{+}(>)=\U_{\Z}^{+}(>;\a_{1})$.
\eit
\edefn

\bdefn
Let $\b,\b'\in R_{re}^{+}(<)$ and $\b\leq\b'$.
\bit
\item[(1)]
Let $\U^{+}(<;\b,\b')$ be the $\Q(q_{1})$-subspace of $\U^{+}$
with basis $B(<;\b,\b')$.
\item[(2)]
Let $\U^{+}(<;\b)$ be the $\Q(q_{1})$-subspace of $\U^{+}$
with basis $B(<;\b)$. We also set $\U^{+}(<)=\U^{+}(<;\a_{0})$.
\item[(3)]
Let $\U^{+}_{\Z}(<;\b,\b')$ be 
the free $\Z [q_{1},q_{1}^{-1}]$-submodule of $\U_{\Z}^{+}$ 
with basis $B(<;\b,\b')$. 
\item[(4)]
Let $\U_{\Z}^{+}(<;\b)$ be 
the free $\Z [q_{1},q_{1}^{-1}]$-submodule of $\U_{\Z}^{+}$ 
with basis $B(<;\b)$.
We also set $\U_{\Z}^{+}(<)=\U_{\Z}^{+}(<;\a_{0})$.
\eit
\edefn

\bprop\label{prop-BCP}\cite[Prop.2.3]{BCP}
Let $\a\in R_{re}^{+}(>)$ and let $\b\in R_{re}^{+}(<)$.
Then, both of $\U_{\Z}^{+}(>;\a_{1},\a)$ and $\U_{\Z}^{+}(<;\b,\a_{0})$ are 
closed under multiplication.
In fact, for $n\geq 0$, we have
\bit
\item[(1)]
$\U_{\Z}^{+}(>;\a_{1},n\d+\a_{1})
=\{x\in\U_{\Z}^{+}|\ T_{1}(T_{0}T_{1})^{n}(x)\in\U_{\Z}^{-}\U_{\Z}^{0}\}$,
\item[(2)]
$\U_{\Z}^{+}(>;\a_{1},(2n+2)\d-\a_{0})
=\{x\in\U_{\Z}^{+}|\ (T_{0}T_{1})^{n+1}(x)\in\U_{\Z}^{-}\U_{\Z}^{0}\}$,
\item[(3)]
$\U_{\Z}^{+}(<;2n\d+\a_{0},\a_{0})
=\{x\in\U_{\Z}^{+}|\ T_{0}^{-1}(T_{1}^{-1}T_{0}^{-1})^{n}(x)
\in\U_{\Z}^{-}\U_{\Z}^{0}\}$,
\item[(4)]
$\U_{\Z}^{+}(<;(n+1)\d-\a_{1},\a_{0})
=\{x\in\U_{\Z}^{+}|\ (T_{1}^{-1}T_{0}^{-1})^{n+1}(x)
\in\U_{\Z}^{-}\U_{\Z}^{0}\}$.
\eit
\eprop
\bpf
Beck et al.treated the simply-laced case, but
their proof of [Prop.2.3] is applicable to our case.
\epf

\bcor\label{cor-z-convex-1} 
Let $\a,\a'\in R_{re}^{+}(>)$ and $\a\leq\a'$.
Let $\b,\b'\in R_{re}^{+}(<)$ and $\b\leq\b'$.
Then, each of the following is closed under multiplication$:$
\bit
\item[(1)]
$\U_{\Z}^{+}(>;\a,\a')$,
\item[(2)]
$\U_{\Z}^{+}(>;\a)$,
\item[(3)]
$\U_{\Z}^{+}(<;\b,\b')$,
\item[(4)]
$\U_{\Z}^{+}(<;\b)$.
\eit
\ecor
\bpf
For $n,m\geq 0$, we have the following isomorphisms of
$\Z[q_{1},q_{1}^{-1}]$-modules, which commute with the
multiplication in $\U$:
\bal
\U_{\Z}^{+}(>;\a_{1},n\d+\a_{1})&\overset{T_{\pa}^{m+1}}{\longleftarrow}
\U_{\Z}^{+}(>;(m+1)\d+\a_{1},(n+m+1)\d+\a_{1})\\
&\overset{T_{1}^{-1}\ast}{\longleftarrow}
\U_{\Z}^{+}(<;(n+m+1)\d-\a_{1},(m+1)\d-\a_{1}),
\end{align*}
\bal
\U_{\Z}^{+}(>;\a_{1},2\d-\a_{0})
&\overset{T_{\pa}^{m+1}}{\longleftarrow}
\U_{\Z}^{+}(>;(m+1)\d+\a_{1},(2n+2m+4)\d-\a_{0})\\
&\overset{T_{1}^{-1}\ast}{\longleftarrow}
\U_{\Z}^{+}(<;(2n+2m+2)\d+\a_{0},(m+1)\d-\a_{1}),
\end{align*}
\bal
\U_{\Z}^{+}(<;2n\d+\a_{0},\a_{0})&\overset{T_{\pa}^{-m}}{\longleftarrow}
\U_{\Z}^{+}(<;(2n+2m)\d+\a_{0},2m\d+\a_{0})\\
&\overset{T_{1}^{-1}\ast}{\longleftarrow}
\U_{\Z}^{+}(>;(2m+2)\d-\a_{0},(2n+2m+2)\d-\a_{0}),
\end{align*}
\bal
\U_{\Z}^{+}(<;(n+1)\d-\a_{1},\a_{0})&\overset{T_{\pa}^{-m}}{\longleftarrow}
\U_{\Z}^{+}(<;(n+m+1)\d-\a_{1},2m\d+\a_{0})\\
&\overset{T_{1}^{-1}\ast}{\longleftarrow}
\U_{\Z}^{+}(>;(2m+2)\d-\a_{0},(n+m+1)\d+\a_{1}).
\end{align*}
(1) and (3) follow from Proposition \ref{prop-BCP}
by using the above morphisms.
(2) and (4) follow from (1) and (3) respectively.
\epf

\bcor\label{cor-convex-1} 
Let $\a,\a'\in R_{re}^{+}(>)$ and let $\b,\b'\in R_{re}^{+}(<)$.
Then, each of the following is closed under multiplication$:$
\bit
\item[(1)]
$\U^{+}(>;\a,\a')$,
\item[(2)]
$\U^{+}(>;\a)$,
\item[(3)]
$\U^{+}(<;\b,\b')$,
\item[(4)]
$\U^{+}(<;\b)$.
\eit
\ecor
\bpf
This follows from Corollary \ref{cor-z-convex-1}.
\epf


\bdefn 
Let $n\geq 0$. 
\bit
\item[(1)]
Let $\U^{+}(0)$ be the $\Q(q_{1})$-subalgebra of $\U^{+}$ 
generated by $\ps_{i}$ for $i\geq 1$.
\item[(2)]
Let $\U^{+}(0;n)$ be the $\Q(q_{1})$-subalgebra of $\U^{+}$ 
generated by $\ps_{i}$ for $1\leq i\leq n$.
We understand that $\U^{+}(0;0)=\Q(q_{1})$.
\eit
\edefn

We introduce new imaginary root vectors.

\bdefn\label{defn-P}
For $n\geq 0$, we define the elements $P_{n}$ of $\U^{+}(0;n)$ by the 
induction on $n$ as follows$:$ we set $P_{0}=1$ and
\bal 
P_{n}=[2n]_{1}^{-1}\sum_{k=0}^{n-1}P_{k}\ps_{n-k}q^{-k}\ \ \text{for}\ n\geq 1.
\end{align*} 
\edefn

\blem\label{lem-P}
Let $n,m\geq 0$. Then,
\bit
\item[(1)]
$[P_{n},P_{m}]=0$,
\item[(2)]
$T_{\pa}(P_{n})=P_{n}$,
\item[(3)]
$T_{1}^{-1}\ast(P_{n})=P_{n}$.
\eit
\elem  
\bpf
This follows from Proposition \ref{p-3} and Proposition \ref{p-2}.
\epf

\blem\label{lem-det}
Let $n\geq 1$. Then,
\bal
\ps_{n}=(-1)^{n-1}q^{-n}\left| 
\begin{array}{cccc}
\a_{1} &1 & 0 &0 \\
\a_{2} &P_{1} &\ddots  &0 \\
\vdots &\vdots &\ddots &1 \\ 
\a_{n} &P_{n-1} &\cdots  &P_{1}
\end{array} 
\right| 
\end{align*}
where we set $\a_{k}=P_{k}[2k]_{1}q^{k}$ for $1\leq k\leq n$.
\elem

\bpf
We argue by the induction on $n$.
The case where $n=1$ is clear: $\ps_{1}=P_{1}[2]_{1}$.
Assume that $n\geq 2$. Then, 
\begin{align*}
\left|
\begin{array}{cccc}
\a_{1} &1 & 0 &0 \\
\a_{2} &P_{1} &\ddots  &0 \\
\vdots &\vdots &\ddots &1 \\ 
\a_{n} &P_{n-1} &\cdots  &P_{1}
\end{array} 
\right| 
&= \sum_{k=1}^{n-1}(-1)^{k+1}P_{k}
\left|
\begin{array}{cccc}
\a_{1} &1 & 0 &0 \\
\a_{2} &P_{1} &\ddots  &0 \\
\vdots &\vdots &\ddots &1 \\ 
\a_{n-k} &P_{n-k-1} &\cdots  &P_{1}
\end{array}
\right| \\
&\ -(-q)^{n}[2n]_{1}P_{n},
\end{align*}
which is checked by developing the left hand side 
with respect to the last row.
Using the induction hypothesis, we see that this is equal to
\bal
(-q)^{n}\Big(\sum_{k=1}^{n-1}P_{k}\ps_{n-k}q^{-k}-[2n]_{1}P_{n}\Big).
\end{align*}
On the other hand, it follows from Definition \ref{defn-P} that
\bal
\ps_{n}=[2n]_{1}P_{n}-\sum_{k=1}^{n-1}P_{k}\ps_{n-k}q^{-k}.
\end{align*}
Thus, the induction proceeds.
\epf

\bcor
Let $n\geq 0$. Then, $\U^{+}(0;n)$ is generated 
by $P_{k}$ for $0\leq k\leq n$. 
\ecor
\bpf
This follows from Lemma \ref{lem-det}.
\epf

\brem 
$\rm{Lemma}\ \ref{lem-det}$ suggests that $P_{n}$ and $\ps_{n}$ 
are q-analogues of the complete symmetric function and 
two times of the power sum respectively \cite[1.2.ex.8]{M}.
\erem

\brem
If we define $E_{n\d}\in\U^{+}(0;n)$ for $n\geq 1$ as in $\rm{Appendix\ B}$, 
then we have 
\bal 
\sum_{k\geq 0}P_{k}u^{k}=\exp\Big(\sum_{n\geq 1}[2n]_{1}^{-1}E_{n\d}u^{n}\Big).
\end{align*} 
\erem

Let us study the commutation relations between $P_{n}$ and 
the real root vectors.

\blem\label{lem-com-PE}
Let $n\geq 0$. Then, 
\bal
P_{n}E_{m\d +\a_{1}}=\sum_{i=0}^{n}E_{(m+n-i)\d+\a_{1}}P_{i}[2n-2i+1]_{1}.
\end{align*}
\elem

\bpf
We can assume that $n\geq 1$. We argue by the induction on $n$.
By virtue of Lemma \ref{lem-P} (2), we can assume that $m=0$.
By Definition \ref{defn-P} and Proposition \ref{p-2} (5), we have
\bal 
&P_{n}E_{\a_{1}}[2n]_{1}\\
&=\sum_{k=0}^{n-1}P_{k}\ps_{n-k}E_{\a_{1}}q^{-k}\\
&=\sum_{k=0}^{n-1}P_{k}\Big(E_{\a_{1}}\ps_{n-k}
+(1+q^{-1})\sum_{i=0}^{n-k-1}b_{2(n-k-i)}
E_{(n-k-i)\d +\a_{1}}\ps_{i}\Big)q^{-k}.
\end{align*}
Using the induction hypothesis, we see that this is equal to
\bal
&\sum_{k=0}^{n-1}\sum_{i=0}^{k}E_{(k-i)\d +\a_{1}}P_{i}\ps_{n-k}
[2k-2i+1]_{1}q^{-k}\\
&\ +(1+q^{-1})\sum_{k=0}^{n-1}\sum_{i=0}^{n-k-1}\sum_{j=0}^{k}
E_{(n-i-j)\d +\a_{1}}P_{j}\ps_{i}b_{2(n-k-i)}[2k-2j+1]_{1}q^{-k}.
\end{align*}
In the former term, putting $k-i=n-i',\ i=j'$, we see that 
$k=n-i'+j'$ and that 
$0\leq k\leq n-1,\ 0\leq i\leq k$ is equivalent to 
$1\leq i'\leq n,\ 0\leq j'\leq i'-1$.
In the latter, putting $i+j=i',\ n-k-i=k'$, we see that
$i=i'-j,\ k=n-i'+j-k'$ and that 
$0\leq k\leq n-1,\ 0\leq i\leq n-k-1,\ 0\leq j\leq k$ is equivalent to 
$0\leq i'\leq n-1,\ 0\leq j\leq i',\ 1\leq k'\leq n-i'$.
Hence, using Lemma \ref{b-4}, we have
\bal
&P_{n}E_{\a_{1}}[2n]_{1}\\
&=\sum_{i=1}^{n}\sum_{j=0}^{i-1}E_{(n-i)\d +\a_{1}}P_{j}\ps_{i-j}
[2n-2i+1]_{1}q^{-n+i-j}\\
&\ +(1+q^{-1})\sum_{i=0}^{n-1}\sum_{j=0}^{i}\sum_{k=1}^{n-i}
E_{(n-i)\d +\a_{1}}P_{j}\ps_{i-j}b_{2k}[2n-2i-2k+1]_{1}q^{-n+i-j+k}\\
&=\sum_{i=1}^{n}\sum_{j=0}^{i-1}E_{(n-i)\d +\a_{1}}P_{j}\ps_{i-j}
[2n-2i+1]_{1}q^{-n+i-j}\\
&\ +\sum_{i=0}^{n-1}\sum_{j=0}^{i}
E_{(n-i)\d +\a_{1}}P_{j}\ps_{i-j}[2n-2i+1]_{1}(q^{n-i}-q^{-n+i})q^{-j}.
\end{align*}
In the latter term, on account of the factor $(q^{n-i}-q^{i-n})$, 
we can include the case for $i=n$. 
Thus, using Definition \ref{defn-P}, we have
\bal
P_{n}E_{\a_{1}}[2n]_{1}
&=\sum_{i=1}^{n}\sum_{j=0}^{i-1}E_{(n-i)\d +\a_{1}}P_{j}\ps_{i-j}
[2n-2i+1]_{1}q^{n-i-j}\\
&\ +\sum_{i=0}^{n}
E_{(n-i)\d +\a_{1}}P_{i}\ps_{0}[2n-2i+1]_{1}(q^{n-i}-q^{-n+i})q^{-i}\\
&=\sum_{i=0}^{n}E_{(n-i)\d +\a_{1}}P_{i}[2n-2i+1]_{1}
([2i]_{1}q^{n-i}+[2n-2i]_{1}q^{-i})\\
&=\sum_{i=0}^{n}E_{(n-i)\d +\a_{1}}P_{i}[2n-2i+1]_{1}[2n]_{1}.
\end{align*}
Thus, the induction proceeds.
\epf

\blem\label{lem-com-P,EE}
Let $n\geq 1,\ m\geq 0$. Then,
\bit
\item[(1)]
$[P_{n},E_{(2m+2)\d-\a_{0}}]=\sum_{i=0}^{n-1}x_{i}P_{i}$ 
for some $x_{i}\in\U^{+,h}(>;(m+1)\d+\a_{1})$ with  $\h(x_{i})=2$,
\item[(2)]
$[E_{2m\d+\a_{0}},P_{n}]=\sum_{i=0}^{n-1}P_{i}z_{i}$ 
for some $z_{i}\in\U^{+,h}(<;(m+1)\d-\a_{1})$ with $\h(z_{i})=-2$.
\eit
\elem
\bpf
By using $T_{1}^{-1}\ast$, (2) is reduced to (1).
By using $T_{\pa}^{m}$, (1) is reduced to the case where $m=0$,
which we shall check now.
By Corollary \ref{c-1} (5) and Lemma \ref{lem-com-PE}, we have
\bal
&P_{n}E_{2\d-\a_{0}}[4]_{1}\\
&=P_{n}[E_{\d+\a_{1}},E_{\a_{1}}]_{q}\\
&=\sum_{i=0}^{n}E_{(n-i+1)\d+\a_{1}}P_{i}[2n-2i+1]_{1}E_{\a_{1}}\\
&\ -q\sum_{i=0}^{n}E_{(n-i)\d+\a_{1}}P_{i}[2n-2i+1]_{1}E_{\d+\a_{1}}\\
&=\sum_{i=0}^{n}\sum_{j=0}^{i}E_{(n-i+1)\d+\a_{1}}
E_{(i-j)\d+\a_{1}}P_{j}[2n-2i+1]_{1}[2i-2j+1]_{1}\\
&\ -q\sum_{i=0}^{n}\sum_{j=0}^{i}E_{(n-i)\d+\a_{1}}
E_{(i-j+1)\d+\a_{1}}P_{j}[2n-2i+1]_{1}[2i-2j+1]_{1}
\end{align*}
The $j\leq i-1$ part in the former term plus 
the $i\leq n-1$ part in the latter
can be written in the desired form
(i.e. as in the right hand side of the statement) by Corollary \ref{c-2}.
The remainder is equal to
\bal
&\sum_{i=0}^{n}E_{(n-i+1)\d+\a_{1}}E_{\a_{1}}P_{i}[2n-2i+1]_{1}\\
&\ -q\sum_{j=0}^{n}E_{\a_{1}}E_{(n-j+1)\d+\a_{1}}P_{j}[2n-2j+1]_{1}\\
&=\sum_{i=0}^{n}[E_{(n-i+1)\d+\a_{1}},E_{\a_{1}}]_{q}P_{i}[2n-2i+1]_{1}.
\end{align*}
The $i=n$ part is equal to $E_{2\d-\a_{0}}P_{n}[4]_{1}$ 
by Corollary \ref{c-1} (5).
The $i\leq n-1$ part can be written in the desired form
by Proposition \ref{p-1}. Thus, we obtain the lemma.
\epf

\bprop\label{prop-convex-2}
Let $n\geq 0$ and let $\a\in R_{re}^{+}(>)$, $\b\in R_{re}^{+}(<)$. 
\bit
\item[(1)]
Let $x\in\U^{+,h}(>;\a),\ y\in\U^{+,h}(0;n)$. Then, 
$yx=\sum_{i}x_{i}y_{i}$ for some
$x_{i}\in\U^{+,h}(>;\a),\ y_{i}\in\U^{+,h}(0;n)$
with $\h(x_{i})=\h(x),\ \ih(y_{i})\leq\ih(y)$.
\item[(2)]
$\U^{+}(>;\a)\U^{+}(0;n)$ is closed under multiplication.
\item[(3)]
Let $y\in\U^{+,h}(0;n),\ z\in\U^{+,h}(<;\b)$. Then, 
$zy=\sum_{i}y_{i}z_{i}$ for some
$y_{i}\in\U^{+,h}(0;n),\ z_{i}\in\U^{+,h}(<;\b)$
with $\ih(y_{i})\leq\ih(y),\ \h(z_{i})=\h(z)$.
\item[(4)]
$\U^{+}(0;n)\U^{+}(<;\b)$ is closed under multiplication. 
\eit
\eprop

\bpf
We prove (1) by the induction on $\h(x)+\ih(y)$.
We can assume that $\h(x)\geq 1$ and $\ih(y)\geq 1$.
First, assume that $x=x'x''$ for some $x',x''\in\U^{+,h}(>)$ 
with $\h(x')<\h(x),\ \h(x'')<\h(x)$. 
Since $\h(x')+\ih(y)<\h(x)+\ih(y)$, by the induction hypothesis, we have
$yx=yx'x''=\sum x_{1}y_{1}x''$ for some 
$x_{1}\in\U^{+,h}(>;\a),\ y_{1}\in\U^{+,h}(0;n)$ with 
$\h(x_{1})=\h(x'),\ \ih(y_{1})\leq\ih(y)$.
Since $\h(x'')+\ih(y_{1})<\h(x)+\ih(y)$, by the induction hypothesis, we have
$yx=\sum x_{1}x_{2}y_{2}$ for some 
$x_{2}\in\U^{+,h}(>;\a),\ y_{2}\in\U^{+,h}(0;n)$ with
$\h(x_{2})=\h(x''),\ \ih(y_{2})\leq\ih(y_{1})$.
By Corollary \ref{cor-convex-1},
we have $x_{1}x_{2}\in\U^{+,h}(>;\a)$.
We also have
$\h(x_{1}x_{2})=\h(x')+\h(x'')=\h(x)$ and $\ih(y_{2})\leq\ih(y)$.
Thus, we are reduced to the case where
$x=E_{k\d+\a_{1}}$ with $k\d+\a_{1}\geq\a$ or 
$E_{(2k+2)\d-\a_{0}}$ with $(2k+2)\d-\a_{0}\geq\a$.
Similarly, we can also assume that $y=P_{i}$ for $1\leq i\leq n$.
Thus, (1) follows from Lemma \ref{lem-com-PE} and Lemma \ref{lem-com-P,EE}.

We prove (2). 
Let $x,x'\in\U^{+,h}(>;\a),\ y,y'\in\U^{+,h}(0;n)$. 
By (1), we have 
$xyx'y'=xx_{1}y_{1}y'$ for some 
$x_{1}\in\U^{+,h}(>;\a),\ y_{1}\in\U^{+,h}(0;n)$. 
Thus, (2) follows from Corollary \ref{cor-convex-1}.

Applying $T_{1}^{-1}\ast$ to (1) and (2), we obtain (3) and (4) 
respectively.
\epf

\brem
In the above proof of $(1)$, 
we understand that each suffix of $x,y$ indicates the number of 
real suffices$:$
$x_{1}=x_{i},\ y_{1}=y_{i},\ x_{2}=x_{ij},\ y_{2}=y_{ij}$
where $i,j$ belong to some finite sets $I_{1},I_{2}$ respectively.
We shall sometimes use such a notation.
\erem

\bprop\label{prop-convex-3}
Let $\a\in R_{re}^{+}(>),\ \b\in R_{re}^{+}(<)$. 
\bit
\item[(1)]
Let $x\in\U^{+,h}(>;\a),\ z\in\U^{+,h}(<;\b)$. Then,
$zx=\sum_{i}x_{i}y_{i}z_{i}$ for some 
$x_{i}\in\U^{+,h}(>;\a)$, $y_{i}\in\U^{+,h}(0),\ z_{i}\in\U^{+,h}(<;\b)$
with $\h(x_{i})\leq\h(x),\ \h(z_{i})\geq\h(z)$.
\item[(2)]
$\U^{+}(>;\a)\U^{+}(0)\U^{+}(<;\b)$ is closed under 
multiplication.
\eit
\eprop

\bpf
(1)
We argue by the induction on $\h(x)-\h(z)$.
If $\h(x)=0$ or $\h(z)=0$, the statement is clear.
We assume that $\h(x)\geq 1$ and $\h(z)\leq -1$.
First, assume that $x=x'x''$ with $x',x''\in\U^{+,h}(>)$ 
and $\h(x')<\h(x),\ \h(x'')<\h(x)$. 
Since $\h(x')-\h(z)<\h(x)-\h(z)$, by the induction hypothesis, we have
$zx=zx'x''=\sum x_{1}y_{1}z_{1}x''$ for some 
$x_{1}\in\U^{+,h}(>;\a),\ y_{1}\in\U^{+,h}(0)$, $z_{1}\in\U^{+,h}(<;\b)$
with  $\h(x_{1})\leq\h(x'),\ \h(z_{1})\geq\h(z)$.
Since $\h(x'')-\h(z_{1})<\h(x)-\h(z)$, by the induction hypothesis, we have
$zx=\sum x_{1}y_{1}x_{2}y_{2}z_{2}$ for some 
$x_{2}\in\U^{+,h}(>;\a),\ y_{2}\in\U^{+,h}(0),\ z_{2}\in\U^{+,h}(<;\b)$
with  $\h(x_{2})\leq\h(x''),\ \h(z_{2})\geq\h(z_{1})$.
By Proposition \ref{prop-convex-2} (1), we have
$zx=\sum x_{1}x_{3}y_{3}y_{2}z_{2}$ for some  
$x_{3}\in\U^{+,h}(>;\a),\ y_{3}\in\U^{+,h}(0)$
with  $\h(x_{3})=\h(x_{2})$.
Then, $x_{1}x_{3}\in\U^{+,h}(>;\a),\ y_{3}y_{2}\in\U^{+,h}(0)$ with
$\h(x_{1}x_{3})\leq\h(x')+\h(x'')=\h(x),\ \h(z_{2})\geq\h(z)$.
Thus, (1) is reduced to the case where
$x=E_{k\d+\a_{1}}$ with $k\d+\a_{1}\geq\a$ or 
$x=E_{2k\d-\a_{0}}$ with $2k\d-\a_{0}\geq\a$.
Similarly, we can also assume that
$z=E_{k\d-\a_{1}}$ with $k\d-\a_{1}\leq\b$ or 
$z=E_{2k\d+\a_{0}}$ with $2k\d+\a_{0}\leq\b$.
Now, (1) follows from 
Proposition \ref{p-2} (2), Corollary \ref{c-3}, and Corollary \ref{cor-appA},
which will be shown in Appendix A.

(2)
Let $x,x'\in\U_{\Z}^{+,h}(>;\a),
\ y,y'\in\U_{\Z}^{+,h}(0),
\ z,z'\in\U_{\Z}^{+,h}(<;\b)$. 
By (1) and Proposition \ref{prop-convex-2}, we have 
$xyzx'y'z'=xyx_{1}y_{1}z_{1}y'z'=xx_{2}y_{2}y_{1}y_{2}'z_{2}z'$ 
for some  $x_{1},x_{2}\in\U_{\Z}^{+,h}(>;\a),
\ y_{1},y_{2},y_{2}'\in\U_{\Z}^{+,h}(0),\ z_{1},z_{2}\in\U_{\Z}^{+,h}(<;\b)$.
Hence, (2) follows from Corollary \ref{cor-convex-1}.
\epf

\bcor\label{cor-surj}
The $\Q(q_{1})$-linear map 
$\U^{+}(>)\otimes\U^{+}(0)\otimes\U^{+}(<)\longrightarrow\U^{+}$
given by multiplication is surjective.
\ecor
\bpf
Let $V=\U^{+}(>)\U^{+}(0)\U^{+}(<)$ be the image of the above morphism.
Since $V$ contains the generators of the $\Q(q_{1})$-algebra $\U^{+}$,
it follows from Proposition \ref{prop-convex-3} (2) that $V=\U^{+}$.
\epf

In Section 6, we shall prove that 
the above morphism is an isomorphism 
and that both of $\{P_{n}|\ n\geq 1\}$ and $\{\ps_{n}|\ n\geq 1\}$
are algebraically independent over $\Q(q_{1})$; 
thus, we shall obtain $\Q(q_{1})$-bases of $\U^{+}$.


\section{Lemmas on Coproduct}

In this section, $\equiv$ means the congruence 
modulo $\U^{0}\U^{+}(<)_{\leq -2}\otimes\U^{+}$, 
unless otherwise stated.

\blem\label{lem-copro-Eps}
We have
\bal
\D(E_{\d -\a_{1}})
&=ck_{1}^{-1}\otimes E_{\d -\a_{1}}+E_{\d -\a_{1}}\otimes 1
 +(q^{2}-q^{-2})k_{1}E_{\a_{0}}\otimes E_{\a_{1}},\tag{1}\\
\D(\ps_{1})
&=c\otimes \ps_{1} +\ps_{1}\otimes 1
 +(q-q^{-1})[3]_{1}k_{1}E_{\d -\a_{1}}\otimes E_{\a_{1}}\tag{2}\\
&\ +(1-q^{-3})(q^{4}-1)k_{1}^{2}E_{\a_{0}}\otimes E_{\a_{1}}^{2}.
\end{align*}
\elem
\bpf
(1)
By Lemma \ref{l-3} (1), we have
\bal
\D(E_{\d -\a_{1}})
&=\D([E_{\a_{0}},E_{\a_{1}}]_{q^{-2}})\\
&=[E_{\a_{0}}\otimes 1 +k_{0}\otimes E_{\a_{0}},
E_{\a_{1}}\otimes 1+k_{1}\otimes E_{\a_{1}}]_{q^{-2}}\\
&=E_{\d -\a_{1}}\otimes 1
+(k_{0}E_{\a_{1}}-q^{-2}E_{\a_{1}}k_{0})\otimes E_{\a_{0}}\\
&\ +(E_{\a_{0}}k_{1}-q^{-2}k_{1}E_{\a_{0}})\otimes E_{\a_{1}}
+ck_{1}^{-1}\otimes E_{\d -\a_{1}}.
\end{align*}
The second term vanishes and we obtain (1).

(2)
By Definition \ref{defn-ps} and (1), we have 
\bal
\D(\ps_{1})
&=\D([E_{\d-\a_{1}},E_{\a_{1}}]_{q^{-1}})\\
&=\big[ ck_{1}^{-1}\otimes E_{\d -\a_{1}}+E_{\d -\a_{1}}\otimes 1
 +(q^{2}-q^{-2})k_{1}E_{\a_{0}}\otimes E_{\a_{1}},\\
&\ k_{1}\otimes E_{\a_{1}}+E_{\a_{1}}\otimes 1\big]_{q^{-1}}\\
&=c\otimes \ps_{1}+(E_{\d -\a_{1}}k_{1}-q^{-1}k_{1}E_{\d -\a_{1}})\otimes 
E_{\a_{1}}\\
&\ +(q^{2}-q^{-2})(k_{1}E_{\a_{0}}k_{1}-q^{-1}k_{1}^{2}E_{\a_{0}})\otimes
E_{\a_{1}}^{2}\\
&\ +c(k_{1}^{-1}E_{\a_{1}}-q^{-1}E_{\a_{1}}k_{1}^{-1})\otimes E_{\d -\a_{1}}
+\ps_{1}\otimes 1\\
&\ +(q^{2}-q^{-2})(k_{1}E_{\a_{0}}E_{\a_{1}}-q^{-1}E_{\a_{1}}k_{1}E_{\a_{0}})
\otimes E_{\a_{1}}\\
&=c\otimes \ps_{1} +\ps_{1}\otimes 1
+((q-q^{-1})+(q^{2}-q^{-2}))k_{1}E_{\d -\a_{1}}\otimes E_{\a_{1}}\\
&\ +(q^{2}-q^{-2})(q^{2}-q^{-1})k_{1}^{2}E_{\a_{0}}\otimes E_{\a_{1}}^{2}.
\end{align*}
Thus, we obtain (2).
\epf

\blem\label{lem-convex-ps}
Let $n,r\geq 1$,\ let $\b\in R^{+}_{re}(<)$, 
and let $x\in\U^{+}(<;\b)\cap\U^{+}_{n\d-r\a_{1}}$. Then, 
$[x,\ps_{1}]\in\U^{+}(<;\b)$.
\elem
\bpf
We argue by the induction on $r$. 
First, assume that $x=x_{1}x_{2}$ 
for some $x_{i}\in\U^{+}(<;\b)\cap\U^{+}_{n_{i}\d-r_{i}\a_{1}}$ 
with $r_{i}\geq 1$. Then we have 
\bal
[x,\ps_{1}]=[x_{1}x_{2},\ps_{1}]=x_{1}[x_{2},\ps_{1}]+[x_{1},\ps_{1}]x_{2}
\end{align*}
and this belongs to $\U^{+}(<;\b)$ 
by the induction hypothesis and Corollary \ref{cor-convex-1}.
Thus, the lemma is reduced to the case where 
$x=E_{n\d-\a_{1}}$, or $x=E_{2m\d+\a_{0}}$ with $m\geq 0$, 
but it follows from Corollary \ref{c-1}.
\epf

\blem\label{lem-copro-E}
Let $n\geq 1$. Then,
\bal
\D(E_{n\d -\a_{1}})
&\equiv c^{n}k_{1}^{-1}\otimes E_{n\d -\a_{1}} 
+(q_{1}-q_{1}^{-1})\sum_{i=1}^{n}c^{n-i}E_{i\d -\a_{1}}\otimes \ps_{n-i}.
\end{align*}
\elem

\bpf
We argue by the induction on $n$. 
The case where $n=1$ follows from Lemma \ref{lem-copro-Eps} (1).
Assuming the case for $n$, we shall prove the case for $n+1$.
Using Corollary \ref{c-1}, the induction hypothesis, and 
Lemma \ref{lem-copro-Eps} (2), in view of 
Proposition \ref{cor-convex-1}, we have
\bal
\D(E_{(n+1)\d -\a_{1}})[3]_{1}!
&=\D([E_{n\d -\a_{1}},\ps_{1}])\\
&\equiv 
\Big[ c^{n}k_{1}^{-1}\otimes E_{n\d -\a_{1}} 
+(q_{1}-q_{1}^{-1})\sum_{i=1}^{n}c^{n-i}E_{i\d -\a_{1}}\otimes \ps_{n-i},\\
&\ c\otimes\ps_{1}+\ps_{1}\otimes 1
 +(q-q^{-1})[3]_{1}k_{1}E_{\d -\a_{1}}\otimes E_{\a_{1}}\Big] \\
&\equiv 
c^{n+1}k_{1}^{-1}\otimes E_{(n+1)\d -\a_{1}}[3]_{1}!\\
&\ +(q_{1}-q_{1}^{-1})\sum_{i=1}^{n}c^{n-i}E_{(i+1)\d -\a_{1}}[3]_{1}!
 \otimes \ps_{n-i}\\
&\ +(q-q^{-1})[3]_{1}c^{n}E_{\d -\a_{1}}
 \otimes [E_{n\d -\a_{1}},E_{\a_{1}}]_{q^{-1}}.
\end{align*}
Applying Proposition \ref{p-2} to the last term, 
we obtain the case for $n+1$. 
\epf

\blem\label{lem-copro-EE}
Let $n\geq 0$. Then, 
\bal
\D(E_{2n\d +\a_{0}})
&\equiv 
c^{2n+1}k_{1}^{-2}\otimes E_{2n\d +\a_{0}}\\
&\ +(q_{1}-q_{1}^{-1})^{2}\sum_{i=1}^{n}c^{2n-i+1}
k_{1}^{-1}E_{i\d -\a_{1}}\\
&\ \otimes\sum_{j=0}^{n-i}b_{2(n-i-j)+1}
\ps_{j}E_{(2n-i-j+1)\d -\a_{1}}\th(n\geq 1).
\end{align*}
\elem
\bpf
We can assume that $n\geq 1$.
Using Corollary \ref{c-1} (6) and Lemma \ref{lem-copro-E}, in view of 
Proposition \ref{cor-convex-1}, we have
\bal
&\D(E_{2n\d +\a_{0}})[4]_{1}\\
&=\D([E_{n\d -\a_{1}},E_{(n+1)\d -\a_{1}}]_{q})\\
&\equiv 
\Big[ c^{n}k_{1}^{-1}\otimes E_{n\d -\a_{1}} 
+(q_{1}-q_{1}^{-1})\sum_{i=1}^{n}c^{n-i}E_{i\d -\a_{1}}\otimes \ps_{n-i},\\
&\ \ c^{n+1}k_{1}^{-1}\otimes E_{(n+1)\d -\a_{1}} 
+(q_{1}-q_{1}^{-1})\sum_{i=1}^{n+1}c^{n-i+1}E_{i\d -\a_{1}}
\otimes \ps_{n-i+1}\Big]_{q}\\
&\equiv
c^{2n+1}k_{1}^{-2}\otimes E_{2n\d +\a_{0}}[4]_{1}
+(q_{1}-q_{1}^{-1})\sum_{i=1}^{n}c^{2n-i+1}k_{1}^{-1}E_{i\d -\a_{1}}
\otimes A_{i}
\end{align*}
where we set 
\bal
A_{i}&=q^{-1}\ps_{n-i}E_{(n+1)\d -\a_{1}}-qE_{(n+1)\d -\a_{1}}\ps_{n-i}
+[E_{n\d -\a_{1}},\ps_{n-i+1}].
\end{align*}
Using Proposition \ref{p-2} (5), we have
\begin{align*}
A_{i}
&=(q^{-1}-q)\ps_{n-i}E_{(n+1)\d -\a_{1}}\\
&\ -q(1+q^{-1})\sum_{j=0}^{n-i-1}b_{2(n-i-j)}\ps_{j}E_{(2n-i-j+1)\d -\a_{1}}\\
&\ +(1+q^{-1})\sum_{j=0}^{n-i}b_{2(n-i-j+1)}\ps_{j}E_{(2n-i-j+1)\d -\a_{1}}\\ 
&=(1+q^{-1})\sum_{j=0}^{n-i}(b_{2(n-i-j+1)}-qb_{2(n-i-j)})
\ps_{j}E_{(2n-i-j+1)\d -\a_{1}}.
\end{align*}
Applying Lemma \ref{b-1} (2), we obtain the lemma.
\epf

\blem\label{cor-copro-E(<)}
Let $x\in\U^{+}(<)$. Then, $\D(x)\in\U^{0}\U^{+}(<)\otimes\U^{+}$.
\elem
\bpf
This follows from Lemma \ref{lem-copro-E} and Lemma \ref{lem-copro-EE} 
together with Corollary \ref{cor-convex-1}.
\epf

\blem\label{lem-copro-braid}
Let $x\in\U^{+}(<)$. Then,
\bal 
\D(T_{i}(x))&\equiv 
(T_{\pa}\otimes T_{\pa})(\D(x))\\
&\ +(q_{1}-q_{1}^{-1})\big[ (T_{\pa}\otimes T_{\pa})(\D(x)),
c^{-1}k_{1}E_{\d -\a_{1}}\otimes T_{\pa}(E_{\a_{1}})\big]\\
&\ \mod\U^{0}\U^{+}(<)_{\leq -2}\otimes\U.
\end{align*}
\elem
\bpf
It is known \cite[37.3.2]{L} that for $i\in I,\ x\in\U$, we have
the following equality in 
$\sum_{\lambda,\mu\in Q}\prod_{\xi\in Q^{+}}
\U_{\lambda +\xi}\otimes\U_{\mu-\xi}\subset\prod_{\lambda,\mu\in Q}
\U_{\lambda}\otimes\U_{\mu}$:
\bal
\D(T_{\pa}(x))=R_{i}(T_{i}\otimes T_{i})(\D(x))R_{i}'
\end{align*}
where we set
\bal
R_{i}&=\sum_{k\geq 0}q_{i}^{k(k-1)/2}(q_{i}-q_{i}^{-1})^{k}[k]_{i}!
\ T_{i}(f_{i}^{(k)})\otimes T_{i}(e_{i}^{(k)}),\\
R_{i}'&=\sum_{k\geq 0}(-1)^{k}q_{i}^{-k(k-1)/2}(q_{i}-q_{i}^{-1})^{k}[k]_{i}!
\ T_{i}(f_{i}^{(k)})\otimes T_{i}(e_{i}^{(k)}).
\end{align*}
Noting that $T_{i}(f_{i})=-k_{i}^{-1}e_{i}$ and 
$T_{i}(e_{i})=-f_{i}k_{i}$, we have 
\bal
\D(T_{\pa}(x))
&=\D(T_{0}(T_{1}(x)))\\
&=R_{0}(T_{0}\otimes T_{0})\big(\D(T_{1}(x))\big) R_{0}'\\
&=R_{0}(T_{0}\otimes T_{0})\big(
R_{1}(T_{1}\otimes T_{1})(\D(x))R_{1}'\big) R_{0}'.
\end{align*}
Noting that $T_{\pa}(f_{1})=-c^{-1}k_{1}E_{\d-\a_{1}}$ and 
$T_{\pa}(e_{1})=-c^{-1}\Omega (E_{\d-\a_{1}})k_{1}^{-1}$, we have 
\bal
\D(T_{\pa}(x))
&=R_{0}
R_{\pa}(T_{\pa}\otimes T_{\pa})(\D(x))R_{\pa}'R_{0}'
\end{align*}
where we set
\bal
R_{\pa}&=\sum_{k\geq 0}q_{1}^{k(k-1)/2}(q_{1}-q_{1}^{-1})^{k}[k]_{1}!
\ T_{\pa}(f_{1}^{(k)})\otimes T_{\pa}(e_{1}^{(k)}),\\
R_{\pa}'&=\sum_{k\geq 0}(-1)^{k}q_{1}^{-k(k-1)/2}(q_{1}-q_{1}^{-1})^{k}[k]_{1}!
\ T_{\pa}(f_{1}^{(k)})\otimes T_{\pa}(e_{1}^{(k)}).
\end{align*}
Hence, for $x\in\U^{+}(<)$, in view of Corollary \ref{cor-copro-E(<)} 
and Corollary \ref{cor-convex-1}, we have
\bal
\D(T_{\pa}(x))
&\equiv 
\big( 1-(q_{1}-q_{1}^{-1})c^{-1}k_{1}E_{\d-\a_{1}}
\otimes T_{\pa}(E_{\a_{1}})\big)
(T_{\pa}\otimes T_{\pa})(\D(x))\\
&\ \times\big( 1+(q_{1}-q_{1}^{-1})c^{-1}k_{1}E_{\d-\a_{1}}
\otimes T_{\pa}(E_{\a_{1}})\big)\\
&\equiv 
(T_{\pa}\otimes T_{\pa})(\D(x))\\
&\ +(q_{1}-q_{1}^{-1})\big[ (T_{\pa}\otimes T_{\pa})(\D(x)),
c^{-1}k_{1}E_{\d -\a_{1}}\otimes T_{\pa}(E_{\a_{1}})\big]
\end{align*}
where $\equiv$ means the congruence modulo 
$\sum_{\lambda,\mu\in Q}\prod_{\xi\in Q^{+}}
(\U^{0}\U^{+}(<)_{\leq -2}\cap\U_{\lambda +\xi})
\otimes\U_{\mu-\xi}$. 
Since both sides belong to $\U\otimes\U$,
we obtain the lemma.
\epf

\blem\label{lem-copro-m}
Let $m\geq 1,\ s\geq 0,\ t\geq 0$. Then we have 
\bal
&\D (E_{2m\d +\a_{0}}^{(s)}E_{m\d -\a_{1}}^{(t)})\\
&\equiv 
c^{(2m+1)s+mt}k_{1}^{-2s-t}\otimes 
E_{2m\d +\a_{0}}^{(s)}E_{m\d -\a_{1}}^{(t)}\\
&\ +c^{(2m+1)s+mt-m}k_{1}^{-2s-t+1}E_{m\d -\a_{1}} \\
&\ \otimes\Big( (q_{1}-q_{1}^{-1})q^{-t}
E_{(m+1)\d -\a_{1}}E_{2m\d +\a_{0}}^{(s-1)}E_{m\d -\a_{1}}^{(t)}
+q_{1}^{-t+1}E_{2m\d +\a_{0}}^{(s)}E_{m\d -\a_{1}}^{(t-1)}\Big) \\
&\ +\sum_{i=1}^{m-1}c^{(2m+1)s+mt-i}k_{1}^{-2s-t+1}E_{i\d -\a_{1}}\otimes
x_{i}\th(m\geq 2)
\end{align*}
for some $x_{i}\in\U^{+}$.
\elem
\bpf
Using Lemma \ref{lem-copro-E} and Lemma \ref{lem-copro-EE},  
in view of Proposition \ref{cor-convex-1}, we have
\bal
&\D (E_{2m\d +\a_{0}}^{s}E_{m\d -\a_{1}}^{t})\\
&\equiv
\Big( c^{2m+1}k_{1}^{-2}\otimes E_{2m\d +\a_{0}}
+(q_{1}-q_{1}^{-1})c^{m+1}k_{1}^{-1}E_{m\d -\a_{1}}\otimes 
E_{(m+1)\d -\a_{1}}\\
&\  +\sum_{i=1}^{m-1} c^{2m-i+1}k_{1}^{-1}E_{i\d -\a_{1}}
\otimes y_{i}\th(m\geq 2)\Big) ^{s}\\
&\ \times 
\Big( c^{m}k_{1}^{-1}\otimes E_{m\d -\a_{1}}+E_{m\d -\a_{1}}\otimes 1
+\sum_{i=1}^{m-1}c^{m-i}E_{i\d -\a_{1}}\otimes y_{i}^{'}\th(m\geq 2)
\Big) ^{t}
\end{align*}
for some $y_{i},y_{i}^{'}\in\U^{+}$.
By Corollary \ref{c-1} (3) and Proposition \ref{prop-convex-2}, we see that
\bal
\D (E_{2m\d +\a_{0}}^{s}E_{m\d -\a_{1}}^{t})
&\equiv
c^{(2m+1)s+mt}k_{1}^{-2s-t}\otimes E_{2m\d +\a_{0}}^{s}E_{m\d -\a_{1}}^{t}\\
&\ +c^{(2m+1)s+mt-m}k_{1}^{-2s-t+1}E_{m\d -\a_{1}} \\
&\ \otimes\Big( (q_{1}-q_{1}^{-1})
(q^{2(s-1)}+q^{2(s-2)}q^{-2}+\cdots +q^{-2(s-1)})q^{-t}\\ 
&\ \times E_{(m+1)\d -\a_{1}}E_{2m\d +\a_{0}}^{s-1}E_{m\d -\a_{1}}^{t}\\
&\ +(1+q^{-1}+\cdots +q^{-(t-1)})
E_{2m\d +\a_{0}}^{s}E_{m\d -\a_{1}}^{t-1}\Big) \\
&\ +\sum_{i=1}^{m-1}c^{(2m+1)s+mt-i}k_{1}^{-2s-t+1}E_{i\d -\a_{1}}
\otimes y_{i}^{''}
\end{align*}
for some $y_{i}^{''}\in\U^{+}$.
Thus, the lemma follows.
\epf

\bcor\label{cor-copro-m}
Let $m\geq 1,\ n\d +r\a\in Q^{+}$ and let 
$x\in\U^{+}(<;2m\d+\a_{0},\d-\a_{1})\cap\U^{+}_{n\d +r\a}$. Then we have
\bal
\D (x)\equiv c^{n}k_{1}^{r}\otimes x +\sum_{i=1}^{m}c^{n-i}k_{1}^{r+1}
E_{i\d -\a_{1}}\otimes z_{i}
\end{align*}
for some $z_{i}\in \U^{0}\U^{+}(<;\d-\a_{1})$.
\ecor
\bpf 
This follows from Lemma \ref{lem-copro-m} 
together with Proposition \ref{cor-convex-1}
by the induction on $m$.
\epf

\blem\label{lem-copro}
Let $m\geq 1,\ n\d +r\a_{1}\in Q^{+}\backslash\{ 0\}$ and let
$x\in\U^{+}(<;2m\d+\a_{0},\d-\a_{1})\cap\U^{+}_{n\d +r\a_{1}}$.
If $\D(x)\equiv c^{n}k_{1}^{r}\otimes x$,
then $x=0$.
\elem
\bpf
We argue by the induction on $m$.
First, we consider the case where $m=1$.
Then, $x=c_{s,t}E_{2\d +\a_{0}}^{(s)}E_{\d -\a_{1}}^{(t)}$
for some $c_{s,t}\in \Q(q_{1})$ where $s=n+r,\ t=-2n-3r,\ s,t\geq 0$,
and $s\geq 1$ or $t\geq1$.
By Lemma \ref{lem-copro-m} with $m=1$, we have
\bal
\D (x)&\equiv c^{n}k_{1}^{r}\otimes x\\
&\ +c^{n-1}k_{1}^{r+1}E_{\d -\a_{1}}\otimes 
c_{s,t}\Big( (q_{1}-q_{1}^{-1})q^{-t}
E_{2\d -\a_{1}}E_{2\d +\a_{0}}^{(s-1)}E_{\d -\a_{1}}^{(t)}\\
&\ +q_{1}^{-t+1}E_{2\d +\a_{0}}^{(s)}E_{\d -\a_{1}}^{(t-1)}\Big).
\end{align*}
Since $s\geq 1$ or $t\geq 1$, we have $c_{s,t}=0$; thus, $x=0$.
Now, assume that $m\geq 2$. Then, $x=\sum_{s,t\geq 0}c_{s,t}
E_{2m\d +\a_{0}}^{(s)}E_{m\d -\a_{1}}^{(t)}x_{s,t}$
where 
$c_{s,t}\in \Q(q_{1}),\ x_{s,t}\in \U^{+,h}(<;2(m-1)\d+\a_{0},\d-\a_{1})$.
By Lemma \ref{lem-copro-m} and by Corollary \ref{cor-copro-m} for $m-1$ 
with $x=x_{s,t}$, in view of Proposition \ref{cor-convex-1}, we have
\bal
\D (x)&\equiv c^{n}k_{1}^{r}\otimes x
+\sum_{i=1}^{m-1}c^{n-i}k_{1}^{r+1}E_{i\d -\a_{1}}\otimes w_{i}\\
&\ +c^{n-m}k_{1}^{r+1}E_{m\d -\a_{1}}\otimes\sum_{s,t\geq 0}c_{s,t}
\Big( (q_{1}-q_{1}^{-1})q^{-t}E_{(m+1)\d -\a_{1}}E_{2m\d +\a_{0}}^{(s-1)}
E_{m\d -\a_{1}}^{(t)}\\
&\ +q_{1}^{-t+1}E_{2m\d +\a_{0}}^{(s)}E_{m\d -\a_{1}}^{(t-1)}\Big) x_{s,t}
\end{align*}
for some $w_{i}\in\U^{+}$.
Since $\D (x)\equiv c^{n}k_{1}\otimes x$, we have 
($w_{i}=0$ for $1\leq i\leq m-1$, and)
$c_{s,t}x_{s,t}=0$ unless $s=t=0$.
Thus, $x\in \U^{+,h}(<;2(m-1)\d+\a_{0},\d-\a_{1})$.
Applying the induction hypothesis, we obtain $x=0$. 
Thus, the induction proceeds.
\epf

\bcor\label{cor-copro}
Let $n\d +r\a_{1}\in Q^{+}\backslash\{ 0\}$ 
and let $x\in\U^{+}(<;\d-\a_{1})\cap \U^{+}_{n\d +r\a_{1}}$. 
If $\D (x)\equiv c^{n}k_{1}^{r}\otimes x$, then $x=0$.
\ecor
\bpf 
This follows from Lemma \ref{lem-copro}.
\epf


\section{Inner Product}

We set $\A=\Q(q_{1})\cap\Z[[q_{1}^{-1}]]\subset\Q((q_{1}^{-1}))$.

\bdefn
\cite[1.2.2,\ 3.1.5]{L}
Let $r$ be the $\Q(q_{1})$-linear map from $\U^{+}$ to $\U^{+}\otimes\U^{+}$
defined as the composition of $\D:\U^{+}\rightarrow
\oplus_{\xi,\mu\in Q^{+}}(\U^{+}_{\xi}k_{\mu}\otimes\U^{+}_{\mu})$ 
and the vector space isomorphism 
from $\oplus_{\xi,\mu\in Q^{+}}(\U^{+}_{\xi}k_{\mu}\otimes\U^{+}_{\mu})$
to $\oplus_{\xi,\mu\in Q^{+}}(\U^{+}_{\xi}\otimes\U^{+}_{\mu})
=\U^{+}\otimes\U^{+}$ that sends $xk_{\mu}\otimes x'$ to $x\otimes x'$ 
for $x\in\U^{+}_{\xi}$ and $x'\in\U^{+}_{\mu}$ .
\edefn

\blem\label{lem-r-E}
Let $n\geq 1$. Then,
\bal
r(E_{n\d -\a_{1}})
&\equiv 1\otimes E_{n\d -\a_{1}} 
+(q_{1}-q_{1}^{-1})\sum_{i=1}^{n}E_{i\d -\a_{1}}\otimes \ps_{n-i}
\end{align*}
where $\equiv$ means the congruence modulo 
$\U^{+}(<)_{\leq -2}\otimes\U^{+}$.
\elem
\bpf
This follows from Lemma \ref{lem-copro-E}.
\epf

Let us recall the inner product introduced by Drinfeld.
\bdefn\label{defn-inner}
\cite[3.4.4,\ 3.4.8]{K1}\ \cite[1.2.3,\ 1.2.5,\ 1.2.13]{L}
Let $(\cdot,\cdot)$ be the $\Q(q_{1})$-valued symmetric bilinear form 
on $\U^{+}$ such that
\bit
\item[(1)]
$(1,1)=1$,
\item[(2)]
$(E_{i},E_{j})=\d_{ij}(1-q_{i}^{-2})^{-1}$ for $i,j\in I$,
\item[(3)]
$(x,yy')=(r(x),y\otimes y')$ for $x,y,y'\in\U^{+}$,
\item[(4)]
$(E_{i}x,y)=(E_{i},E_{i})(x,\ _{i}r(y))$ for $x,y\in\U^{+},\ i\in I$,
\item[(5)]
$(xE_{i},y)=(E_{i},E_{i})(x,r_{i}(y))$ for $x,y\in\U^{+},\ i\in I$,
\item[(6)]
for $\xi,\mu\in Q^{+}$ such that $\xi\neq\mu$ and for 
$x\in\U^{+}_{\xi},\ y\in\U^{+}_{\mu}$, we have $(x,y)=0$.
\item[(7)]
for $\nu\in Q^{+}$, the restriction of $(\cdot,\cdot)$ on $\U^{+}_{\nu}$
is non-degenerate,
\eit
In $(3)$, the bilinear form on $\U^{+}\otimes\U^{+}$ that sends a pair 
$(x_{1}\otimes x_{2},x_{1}'\otimes x_{2}')$ 
to $(x_{1},x_{2})(x_{1}',x_{2}')\in\Q(q_{1})$
is also denoted by $(\cdot,\cdot)$.
\edefn

\bprop\label{prop-inner-L}\cite[40.2.4]{L}
Let ${\bf c_{+}}=({\bf c_{+}}_{i}),\ {\bf c_{+}'}=({\bf c_{+}'}_{i})
\in\oplus_{i\in R_{re}^{+}(>)}\Z_{\geq 0}^{i}$, 
let ${\bf c_{-}}=({\bf c_{-}}_{i}),\ {\bf c_{-}'}=({\bf c_{-}'}_{i})
\in\oplus_{i\in R_{re}^{+}(>)}\Z_{\geq 0}^{i}$,
and let $x,x'\in \U^{+}(0)$. Then,
\bit
\item[(1)]
$({\bf E_{c_{+}}}x{\bf E_{c_{-}}},{\bf E_{c_{+}'}}x'{\bf E_{c_{-}'}})
=\d_{{\bf c_{+}},{\bf c_{+}'}}\d_{{\bf c_{-}},{\bf c_{-}'}}(x,x')
({\bf E_{c_{+}}},{\bf E_{c_{+}}})({\bf E_{c_{-}}},{\bf E_{c_{-}}})$,
\item[(2)]
$({\bf E_{c_{+}}},{\bf E_{c_{+}}})\equiv 1\mod q_{1}^{-1}\A$,
\item[(3)]
$({\bf E_{c_{-}}},{\bf E_{c_{-}}})\equiv 1\mod q_{1}^{-1}\A$.
\eit
\eprop
\bpf
If we set ${\bf h}=(i_{k})_{k\in\Z}$ where
$i_{k}$ is $0$ if $k$ is even and $1$ otherwise, and $p=0$, then 
our symbols $\U^{+}(>)$ and $\U^{+}(<)$ coincide with those
in \cite[40.2.5]{L}, and
$\U^{+}(0)$ is contained in the space 
${\bf P}$ defined in \cite[40.2.3]{L} by Lemma \ref{lem-P}
and the fact that $\U^{+}$ is invariant under $\ast$.
\epf

\brem
It follows from \rm{Corollary} $\ref{cor-surj}$ that $\U^{+}(0)$ 
coincides with 
the above ${\bf P}$.
\erem

\bprop\label{prop-isom}
The $\Q(q_{1})$-linear map 
$\U^{+}(>)\otimes\U^{+}(0)\otimes\U^{+}(<)\longrightarrow\U^{+}$
given by multiplication is an isomorphism.
\eprop
\bpf
By Corollary \ref{cor-surj}, it is enough to show the injectivity. 
Let
$\sum
{\bf E_{c_{+}}}x_{{\bf c_{+}},{\bf c_{-}}}{\bf E_{c_{-}}}=0$ 
where the sum is taken over 
${\bf c_{+}}\in\oplus_{i\in R_{re}^{+}(>)}\Z_{\geq 0}^{i}$ and 
${\bf c_{-}}\in\oplus_{i\in R_{re}^{+}(<)}\Z_{\geq 0}^{i}$; and 
$x_{{\bf c_{+}},{\bf c_{-}}}\in\U^{+}(0)$.
We fix ${\bf c_{+}}\in\oplus_{i\in R_{re}^{+}(>)}\Z_{\geq 0}^{i}$ and 
${\bf c_{-}}\in\oplus_{i\in R_{re}^{+}(<)}\Z_{\geq 0}^{i}$.
Then, for any $x\in\U^{+}(0)$, by Proposition \ref{prop-inner-L}, we have 
$(x_{{\bf c_{+}},{\bf c_{-}}},x)
({\bf E_{c_{+}}},{\bf E_{c_{+}}})({\bf E_{c_{-}}},{\bf E_{c_{-}}})=0$;
thus, $(x_{{\bf c_{+}},{\bf c_{-}}},x)=0$.
Hence, by Proposition \ref{prop-inner-L} together with 
Corollary \ref{cor-surj}, we have 
$(x_{{\bf c_{+}},{\bf c_{-}}},x')=0$ for any $x'\in\U^{+}$.
In view of Definition \ref{defn-inner} (7), we conclude that 
$x_{{\bf c_{+}},{\bf c_{-}}}=0$. The proposition is proved.
\epf

\blem\label{lem-inner-ps}
Let $n\geq 1$. Then, 
\bal
(\ps_{n},\ps_{n})=(1+q^{-1})(q_{1}-q_{1}^{-1})^{-2}b_{2n}.
\end{align*}
\elem
\bpf
By Definition \ref{defn-ps} and \ref{defn-inner} ((4), (5)),  we have
\bal
(\ps_{n},\ps_{n})
&=(E_{n\d -\a_{1}}E_{\a_{1}}-q^{-1}E_{\a_{1}}E_{n\d -\a_{1}},\ps_{n})\\
&=(E_{\a_{1}},E_{\a_{1}})
(E_{n\d -\a_{1}},r_{1}(\ps_{n})-_{1}r(\ps_{n})q^{-1}).
\end{align*}
Applying Corollary \ref{c-4}, we have
\bal
(\ps_{n},\ps_{n})
&=(E_{\a_{1}},E_{\a_{1}})(q_{1}-q_{1}^{-1})q^{-1}(1+q^{-1})
\sum_{i=1}^{n}b_{2i}(E_{n\d -\a_{1}},\ps_{n-i}E_{i\d -\a_{1}}).
\end{align*}
By Definition \ref{defn-inner} ((3), (6)) and Lemma \ref{lem-r-E}, 
we have
\bal
(\ps_{n},\ps_{n})
&=(E_{\a_{1}},E_{\a_{1}})(E_{n\d -\a_{1}},E_{n\d -\a_{1}})
q^{-1}(1+q^{-1})b_{2n}. 
\end{align*}
Since 
$(E_{n\d -\a_{1}},E_{n\d -\a_{1}})=(E_{\a_{1}},E_{\a_{1}})=(1-q^{-1})^{-1}$
by \cite[40.2.4]{L}, we obtain the desired result.
\epf

\blem\label{lem-convex-E}
Let $n,r\geq 1$, let $\b\in R_{re}^{+}(<)$,
and let $x\in\U^{+}(<;\b)\cap\U^{+}_{n\d-r\a_{1}}$. Then, 
$[x,E_{\a_{1}}]_{q^{-r}}\in\U^{+}(0)\U^{+}(<;\b)$.
\elem
\bpf
We argue by the induction on $r$.
First, assume that $x=x_{1}x_{2}$ for some 
$x_{i}\in\U^{+}(0)\cap\U^{+}_{n_{i}\d -r_{i}\a_{1}}$ with $r_{i}\geq 1$.
By Lemma \ref{lem-bracket}, we have 
\bal
[x,E_{\a_{1}}]_{q^{-r}}
&=[x_{1}x_{2},E_{\a_{1}}]_{q^{-r_{1}-r_{2}}}\\
&=x_{1}[x_{2},E_{\a_{1}}]_{q^{-r_{2}}}
+[x_{1},E_{\a_{1}}]_{q^{-r_{1}}}x_{2}q^{-r_{2}}
\end{align*}
and this belongs to $\U^{+}(0)\U^{+}(<;\b)$ 
by the induction hypothesis and Proposition \ref{prop-convex-2}.
Thus, the lemma is reduced to the case where 
$x=E_{n\d-\a_{1}}$ or $x=E_{2m\d-\a_{0}}$ with $m\geq 1$, 
but it follows from Proposition \ref{p-2} (2) and Corollary \ref{c-3} (1). 
\epf

\blem\label{lem-copro-psP}
Let $n\geq 0$. Then,
\bit
\item[(1)]
$\D(\ps_{n})\equiv (q_{1}-q_{1}^{-1})
\sum_{i=0}^{n}c^{n-i}\ps_{i}\otimes\ps_{n-i}$,
\item[(2)]
$\D(P_{n})\equiv \sum_{i=0}^{n}c^{n-i}P_{i}\otimes P_{n-i}$
\eit
where $\equiv$ means the congruence modulo 
$\U^{0}\U^{+}(0)\U^{+}(<)_{\leq -1}\otimes\U^{+}$.
\elem
\bpf
(1)
We can assume that $n\geq 1$. Note that in Lemma \ref{lem-r-E}, 
the truncated part can be written as 
$\sum_{\xi,\mu\in Q^{+};\ \xi +\mu =n\d -\a_{1}}
k_{\mu}x_{\xi}\otimes y_{\mu}$ 
with $x_{\xi}\in \U^{+}(<)_{\leq -2}\cap \U^{+}_{\xi}, 
\ y_{\mu}\in \U^{+}_{\mu}$,
for which we have 
$[k_{\mu}x_{\xi},E_{\a_{1}}]_{q^{-1}}
=k_{\mu}[x_{\xi},E_{\a_{1}}]_{q^{\h (\xi)}}$ and this belongs to
$\U^{0}\U^{+}(0)\U^{+}(<)_{\leq -1}$ by Lemma \ref{lem-convex-E}.
Hence, using Definition \ref{defn-ps} and Lemma \ref{lem-r-E}, we have
\bal
\D(\ps_{n})
&=[\D(E_{n\d -\a_{1}}),\D(E_{\a_{1}})]_{q^{-1}}\\
&\equiv \Big[ c^{n}k_{1}^{-1}\otimes E_{n\d -\a_{1}} 
+(q_{1}-q_{1}^{-1})\sum_{i=1}^{n}c^{n-i}E_{i\d -\a_{1}}\otimes \ps_{n-i},\\
&\ k_{1}\otimes E_{\a_{1}}+E_{\a_{1}}\otimes 1\Big]_{q^{-1}}\\
&\equiv c^{n}\otimes \ps_{n}
+c^{n}(k_{1}^{-1}E_{\a_{1}}-q^{-1}E_{\a_{1}}k_{1}^{-1})
\otimes E_{n\d -\a_{1}}\\
&\ +(q_{1}-q_{1}^{-1})\sum_{i=1}^{n}c^{n-i}\ps_{i}\otimes\ps_{n-i}.
\end{align*}
The second term vanishes and we obtain (1).

(2)
We argue by the induction on $n$. 
We can assume that $n\geq 1$.
By Definition \ref{defn-P}, the induction hypothesis, (1), 
and Proposition \ref{prop-convex-2}, we have
\bal
\D(P_{n})[2n]_{1}
&=\sum_{k=0}^{n-1}\D(P_{k}\ps_{n-k})q^{-k}\\
&\equiv\sum_{k=0}^{n-1}\Big(\sum_{i=0}^{k}c^{k-i}P_{i}\otimes P_{k-i}\Big)
\Big(\sum_{j=0}^{n-k}c^{n-k-j}\ps_{j}\otimes \ps_{n-k-j}\Big)
q^{-k}(q_{1}-q_{1}^{-1})\\
&=\sum_{k=0}^{n-1}\sum_{i=0}^{k}\sum_{j=0}^{n-k}
c^{n-i-j}P_{i}\ps_{j}\otimes P_{k-i}\ps_{n-k-j}q^{-k}(q_{1}-q_{1}^{-1}).
\end{align*}
Putting $i+j=i',\ i=j',\ k-i=k'$, we see that $j=i'-j',\ k=j'+k'$
and that $0\leq k\leq n-1,\ 0\leq i\leq k,\ 0\leq j\leq n-k$
is equivalent to $0\leq i'\leq n,\ 0\leq j'\leq i',
\ 0\leq k'\leq\min (n-i',n-j'-1)$. Hence,
\bal
&\D(P_{n})[2n]_{1}\\
&\equiv\sum_{i=0}^{n}\sum_{j=0}^{i}\sum_{k=0}^{\min (n-i,n-j-1)}
c^{n-i}P_{j}\ps_{i-j}\otimes P_{k}\ps_{n-i-k}q^{-j-k}(q_{1}-q_{1}^{-1})\\
&=\sum_{i=0}^{n}\Big( \sum_{k=0}^{n-i-1}
c^{n-i}P_{i}\ps_{0}\otimes P_{k}\ps_{n-i-k}q^{-i-k}\\
&\ +\th(i\geq 1)\sum_{j=0}^{i-1}c^{n-i}P_{j}\ps_{i-j}q^{-j}
\otimes (P_{n-i}\ps_{0}q^{-n+i}+P_{n-i}[2n-2i]_{1})\Big) (q_{1}-q_{1}^{-1})\\
&=\sum_{i=0}^{n}c^{n-i}P_{i}\otimes P_{n-i}
\big( [2n-2i]_{1}q^{-i}
+[2i]_{1}(q^{-n+i}+[2n-2i]_{1}(q_{1}-q_{1}^{-1}))\big).
\end{align*}
Thus, the induction proceeds.
\epf

\bcor\label{cor-r-psP}
Let $n\geq 0$. Then,
\bit
\item[(1)]
$r(\ps_{n})\equiv\sum_{i=0}^{n}\ps_{i}\otimes\ps_{n-i}
\mod\U^{+}(0)\U^{+}(<)_{\leq -1}\otimes\U^{+}$,
\item[(2)]
$r(P_{n})\equiv \sum_{i=0}^{n}P_{i}\otimes P_{n-i}
\mod\U^{+}(0)\U^{+}(<)_{\leq -1}\otimes\U^{+}$.
\eit
\ecor
\bpf
This follows from Lemma \ref{lem-copro-psP}.
\epf

\blem\label{lem-inner-psP}
For $n\geq 0$, we have
\bal
(\ps_{n},P_{n})=(q_{1}-q_{1}^{-1})^{-1}[2n+1]_{1}.
\end{align*}
\elem
\bpf
We argue by the induction on $n$. 
We can assume that $n\geq 1$.
Using Corollary \ref{cor-r-psP} (1) and Definition \ref{defn-P}, 
in view of Definition \ref{defn-inner} (6), we have
\bal
(\ps_{n},P_{n})[2n]_{1}
&=\Big( \sum_{i=0}^{n}\ps_{i}\otimes\ps_{n-i},
\sum_{j=0}^{n-1}P_{j}\otimes\ps_{n-j}q^{-j}\Big) (q_{1}-q_{1}^{-1})\\
&=\sum_{i=0}^{n-1}(\ps_{i},P_{i})(\ps_{n-i},\ps_{n-i})
q^{-i}(q_{1}-q_{1}^{-1}).
\end{align*}
Applying the induction hypothesis and Lemma \ref{lem-inner-ps}, we have
\bal
(\ps_{n},P_{n})[2n]_{1}
=\sum_{i=0}^{n-1}b_{2(n-i)}[2i+1]_{1}q^{-i}(1+q^{-1})(q_{1}-q_{1}^{-1})^{-2}.
\end{align*}
By Lemma \ref{b-4} (1), the induction proceeds.
\epf

\blem\label{lem-inner-P}
Let $n\geq 0$. Then,
\bal
(P_{n},P_{n})\equiv 1\mod\ q_{1}^{-1}\A.
\end{align*}
\bpf
We can assume that $n\geq 1$.
We argue by the induction on $n$.
By Corollary \ref{cor-r-psP} (2) and Definition \ref{defn-P}, 
and in view of Definition \ref{defn-inner} (6), we have 
\bal
(P_{n},P_{n})[2n]_{1}
&=\Big( \sum_{i=0}^{n}P_{i}\otimes P_{n-i},
\sum_{j=0}^{n-1}P_{j}\otimes\ps_{n-j}q^{-j}\Big)\\
&=\sum_{i=1}^{n}(P_{n-i},P_{n-i})(P_{i},\ps_{i})q^{-n+i}.
\end{align*}
Applying the induction hypothesis and Lemma \ref{lem-inner-psP}, we have
\bal
(P_{n},P_{n})
&=\sum_{i=1}^{n}(1+\a_{i})[2n]_{1}^{-1}
(q_{1}-q_{1}^{-1})^{-1}[2i+1]_{1}q^{-n+i}\\
&=\sum_{i=1}^{n}(1+\a_{i})(q^{n}-q^{-n})^{-1}
(q^{i}+q^{i-1}+\cdots +q^{-i})q^{-n+i}\\
&=\sum_{i=1}^{n}(1+\a_{i})(1+q^{-1}+\cdots +q^{-2i})
(1-q^{-2n})^{-1}q^{-2n+2i}
\end{align*}
for some $\a_{i}\in q_{1}^{-1}\A$.
Thus, $(P_{n},P_{n})\equiv 1\mod q_{1}^{-1}\A$ and the induction proceeds.
\epf
\elem

\bdefn\label{defn-M}\cite[Chap.1]{M}
Let us recall the $\Z$-algebra $\Lambda$ of symmetric functions.
\bit
\item[(1)]
For $n\geq 1$, let $h_{n}$ be 
the complete symmetric functions,
which are algebraically independent and generate $\Lambda$. 
We also set $h_{0}=1$.
\item[(2)]
We define the grading on $\Lambda$ by $\deg(h_{n})=n$ for $n\geq 0$.
\item[(3)]
Let $\D '$ be the algebra homomorphism from $\Lambda$ to
$\Lambda\otimes\Lambda$ given by
$\D'(h_{n})=\sum_{i=0}^{n}h_{i}\otimes h_{n-i}$.
\item[(4)]
Let $(\cdot,\cdot )_{\Lambda}$ be the $\Z$-valued positive definite 
symmetric bilinear form determined by
the following properties$:$ $(h_{n},h_{n})_{\Lambda}=1$ for $n\geq 0;$
$(f,g)_{\Lambda}=0$ for homogeneous $f,g\in\Lambda$ with $deg(f)\neq deg(g);$
and $(f,gh)_{\Lambda}=(\D '(f),g\otimes h)_{\Lambda\otimes\Lambda}$ 
for $f,g,h\in\Lambda$. Here, $(\cdot,\cdot )_{\Lambda\otimes\Lambda}$ is 
the symmetric bilinear form on $\Lambda\otimes\Lambda$ such that 
$(f_{1}\otimes f_{2},g_{1}\otimes g_{2})_{\Lambda\otimes\Lambda}
=(f_{1},g_{1})_{\Lambda}(f_{2},g_{2})_{\Lambda}$ for 
$f_{1},f_{2},g_{1},g_{2}\in\Lambda$.
\item[(5)]
Let $s_{\lambda}$ be the Schur functions
where $\lambda$ runs through the set of partitions.
The $s_{\lambda}$ form an orthonormal basis of $\Lambda$ 
with respect to $(\cdot,\cdot )_{\Lambda}$.
\eit
\edefn

\bdefn\label{defn-Svarphi}
Let $S$ be the $\Z$-subalgebra of $\U^{+}(0)$ 
generated by $P_{n}$ for $n\geq 1$.
Let $\varphi$ be the surjective 
$\Z$-algebra homomorphism from 
$\Z[h_{1},h_{2},\cdots]$ to $S$
given by $\varphi(h_{n})=P_{n}$ for $n\geq 1$.
\edefn

\blem\label{lem-com-varphi}
Let $f\in\Lambda$. Then, 
\bal
r(\varphi(f))\equiv(\varphi\otimes\varphi)(\D'(f))
\mod\U^{+}(0)\U^{+}(<)_{\leq -1}\otimes\U^{+}.
\end{align*}
\elem
\bpf
Let $f,g\in\Lambda$. By Lemma \ref{lem-copro-psP} (2) 
and Proposition \ref{prop-convex-2}, we have
\bal
\D(\varphi(f))\equiv\sum_{i}c^{n_{i}}x_{i}\otimes y_{i}
\mod\U^{0}\U^{+}(0)\U^{+}(<)_{\leq -1}\otimes\U^{+},\\
\D(\varphi(g))\equiv\sum_{j}c^{n_{j}'}x_{j}'\otimes y_{j}'
\mod\U^{0}\U^{+}(0)\U^{+}(<)_{\leq -1}\otimes\U^{+}
\end{align*}
for some $x_{i},y_{i},x_{j}',y_{j}'\in S^{h}$ with 
$n_{i}=\ih(y_{i}),n_{j}'=\ih(y_{j}')$. 
Then, 
\bal
r(\varphi(f))\equiv\sum_{i}x_{i}\otimes y_{i}
\mod\U^{+}(0)\U^{+}(<)_{\leq -1}\otimes\U^{+},\\
r(\varphi(g))\equiv\sum_{i}x_{i}'\otimes y_{i}'
\mod\U^{+}(0)\U^{+}(<)_{\leq -1}\otimes\U^{+}.
\end{align*}
By Proposition \ref{prop-convex-2}, we have 
\bal
\D(\varphi(fg))&=\D(\varphi(f)\varphi(g))=\D(\varphi(f))\D(\varphi(g))\\
&\equiv
\sum_{i,j}c^{n_{i}+n_{j}'}x_{i}x_{j}'\otimes y_{i}y_{j}'
\mod\U^{0}\U^{+}(0)\U^{+}(<)_{\leq -1}\otimes\U^{+}
\end{align*}
and
\bal
r(\varphi(f))r(\varphi(g))
&\equiv\sum_{i,j}x_{i}x_{j}'\otimes y_{i}y_{j}'
\mod\U^{+}(0)\U^{+}(<)_{\leq -1}\otimes\U^{+}.
\end{align*}
Thus, 
\bal
r(\varphi(fg))
&\equiv
r(\varphi(f))r(\varphi(g))
\mod\U^{+}(0)\U^{+}(<)_{\leq -1}\otimes\U^{+}.
\end{align*}
On the other hand, we have
\bal
(\varphi\otimes\varphi)(\D'(fg))
=(\varphi\otimes\varphi)(\D'(f))(\varphi\otimes\varphi)(\D'(g)).
\end{align*}
Thus, the lemma is reduced to the case 
where $f=h_{n}$ for $n\geq 0$, 
but it follows from Corollary \ref{cor-r-psP} (2) and Definition 
\ref{defn-inner} (4).
\epf

\blem\label{lem-inner-varphi}
Let $f,g\in\Lambda$. Then,
\bit
\item[(1)]
$(\varphi(f),\varphi(g))\in\A$,
\item[(2)]
$(\varphi(f),\varphi(g))\equiv(f,g)_{\Lambda}\mod q_{1}^{-1}\A$.
\eit
\elem
\bpf
We can assume that $f,g\in\Lambda$ are homogeneous. 
We can also assume that $\deg(f)=\deg(g)$; otherwise, both sides vanish 
by Definition \ref{defn-inner} (6) and Definition \ref{defn-M} (4).
We prove (1) and (2) at once by the induction on $\deg(f)$. 
First, assume that $g=g_{1}g_{2}$ for homogeneous $g_{i}\in\Lambda$ with 
$\deg(g_{i})\geq 1$.
Then, by Lemma \ref{lem-inner-varphi} 
and Definition \ref{defn-inner} (6), we have
\bal
(\varphi(f),\varphi(g))
&=(\varphi(f),\varphi(g_{1})\varphi(g_{2}))\\
&=(r(\varphi(f)),\varphi(g_{1})\otimes\varphi(g_{2}))\\
&=\big((\varphi\otimes\varphi)(\D'(f)),
\varphi(g_{1})\otimes\varphi(g_{2})\big).
\end{align*}
Hence, by the induction hypothesis, 
we obtain (1): $(\varphi(f),\varphi(g))\in\A$; and 
\bal
(\varphi(f),\varphi(g))
&\equiv(\D'(f),g_{1}\otimes g_{2})_{\Lambda\otimes\Lambda}\mod q_{1}^{-1}\A.
\end{align*}
Hence, in view of Definition \ref{defn-M} (4), we obtain 
(2): $(\varphi(f),\varphi(g))\equiv(f,g)_{\Lambda}\mod q_{1}^{-1}\A$, 
and the induction proceeds.
Thus, the lemma is reduced to the case where $g=h_{n}$ for $n\geq 0$.
Similarly, we can also assume that $f=h_{n}$ for $n\geq 0$.
Thus, the lemma follows from Lemma \ref{lem-inner-P} 
and Definition \ref{defn-M} (4).
\epf

\bdefn\label{defn-Schur}
For a partition $\lambda$, we set 
${\bf S}_{\lambda}=\varphi(s_{\lambda})$.
\edefn

\bprop\label{prop-alg-indep}
We have 
\bit
\item[(1)]
$\varphi$ is an isomorphism of $\Z$-algebras from $\Lambda$ to $S$,
\item[(2)]
$({\bf S}_{\lambda},{\bf S}_{\mu})\equiv\d_{\lambda,\mu}\mod q_{1}^{-1}\A$,
\item[(3)]
the ${\bf S}_{\lambda}$ form a $\Q(q_{1})$-basis of $\U^{+}(0)$,
\item[(4)]
both of $\{ P_{n}|\ n\geq 1\}$ and $\{ \ps_{n}|\ n\geq 1\}$ are 
algebraically independent over $\Q (q_{1})$.
\eit
\eprop
\bpf
By virtue of Lemma \ref{lem-inner-varphi}, if we define 
the $\Z$-valued symmetric bilinear form $(\cdot,\cdot )'$ on $S$ 
as the composition of $(\cdot,\cdot )$ and
the canonical projection from $\A$ to $\A /q_{1}^{-1}\A=\Z$, 
then we have 
\bal
(\varphi(f),\varphi(g))'=(f,g)_{\Lambda}\ \ \text{for}\ f,g\in\Lambda.
\end{align*}
Since $(\cdot,\cdot )_{\Lambda}$ is positive definite, we see
that $\varphi$ is injective; hence, we obtain (1).
Thus,
the ${\bf S}_{\lambda}$'s form a $\Z$-basis of $S$ orthonormal with respect to
$(\cdot,\cdot )'$; hence, they span $\U^{+}(0)$ 
over $\Q(q_{1})$ and satisfy (2).
Now, let $\sum_{\lambda} c_{\lambda}{\bf S}_{\lambda}=0$ 
with $c_{\lambda}\in\Z [q_{1}]$.
Assume that there exists a nonzero $c_{\lambda}$. 
For each of such $c_{\lambda}$, 
take the smallest $n_{\lambda}\geq 0$ such that 
$c_{\lambda}q_{1}^{-n_{\lambda}}\in\Z +q_{1}^{-1}\Z[q_{1}^{-1}]$.
Let $n_{\lambda_{0}}$ be the largest of $n_{\lambda}$. 
Then, $c_{\lambda}q_{1}^{-n_{\lambda_{0}}}\in\A$ and 
$c_{\lambda_{0}}q_{1}^{-n_{\lambda_{0}}}\notin q_{1}^{-1}\A$.
Then, $0=(q_{1}^{-n_{\lambda_{0}}}\sum c_{\lambda}{\bf S}_{\lambda},
{\bf S}_{\lambda_{0}})\equiv c_{\lambda_{0}}q_{1}^{-n_{\lambda_{0}}}
\mod q_{1}^{-1}\A$. 
A contradiction is deduced.
Hence, the ${\bf S}_{\lambda}$'s are linearly independent 
over $\Z [q_{1}]$; thus, 
they are linearly independent over $\Q(q_{1})$.
(3) is proved.

It follows from (3) that for $n\geq 1$, the dimension of 
$\U^{+}(0)\cap\U^{+}_{n\d}$ over $\Q(q_{1})$ 
is equal to the partition number of $n$. (4) follows.
\epf

\bthm\label{thm-q-basis}
Each of the following is a $\Q(q_{1})$-basis of $\U^{+}:$
\bal
\{ {\bf E_{c_{+}}} {\bf E_{c_{0}}'} {\bf E_{c_{-}}}|
\ {\bf c_{+}}\in\oplus_{i\in R_{re}^{+}(>)}\Z_{\geq 0}^{(i)},
\ {\bf c_{0}}\in\oplus_{n\geq 1}\Z_{\geq 0}^{(n)},
\ {\bf c_{-}}\in\oplus_{i\in R_{re}^{+}(<)}\Z_{\geq 0}^{(i)}\},\tag{1}\\
\{ {\bf E_{c_{+}}} {\bf E_{c_{0}}} {\bf E_{c_{-}}}|
\ {\bf c_{+}}\in\oplus_{i\in R_{re}^{+}(>)}\Z_{\geq 0}^{(i)},
\ {\bf c_{0}}\in\oplus_{n\geq 1}\Z_{\geq 0}^{(n)},
\ {\bf c_{-}}\in\oplus_{i\in R_{re}^{+}(<)}\Z_{\geq 0}^{(i)}\},\tag{2}\\
\{ {\bf E_{c_{+}}} {\bf S}_{\lambda} {\bf E_{c_{-}}}|
\ {\bf c_{+}}\in\oplus_{i\in R_{re}^{+}(>)}\Z_{\geq 0}^{(i)},
\ \lambda\ \text{is a partition},
\ {\bf c_{-}}\in\oplus_{i\in R_{re}^{+}(<)}\Z_{\geq 0}^{(i)}\}.\tag{3}
\end{align*}
\ethm
\bpf
This follows from Proposition \ref{prop-isom} and 
Proposition \ref{prop-alg-indep}.
\epf

\brem\label{rem-convex}
$\rm{Corollary}\ \ref{cor-convex-1}\ (3), 
\ \rm{Proposition}\ \ref{prop-convex-2}\ ((2),\ (4))$, and
$\rm{Proposition}\ \ref{prop-convex-3}\ (2)$
are referred to as the convexity of the bases in $(1)$ and $(2)$ 
of $\rm{Theorem}\ \ref{thm-q-basis}$.
\erem

\blem\label{lem-quasi}
The basis in $\rm{Proposition}\ \ref{thm-q-basis}\ (3)$ 
is quasi-orthonormal with respect to the inner product$:$ namely, 
for ${\bf c_{+}},{\bf c_{+}'}\in\oplus_{i\in R_{re}^{+}(>)}\Z_{\geq 0}^{(i)}$,
for ${\bf c_{-}},{\bf c_{-}'}\in\oplus_{i\in R_{re}^{+}(<)}\Z_{\geq 0}^{(i)}$, 
and for partitions $\lambda$ and $\mu$, we have
\bal
({\bf E_{c_{+}}}{\bf S}_{\lambda} {\bf E_{c_{-}}},
{\bf E_{c_{+}'}}{\bf S}_{\mu} {\bf E_{c_{-}'}})
\equiv \d_{{\bf c_{+}},{\bf c_{+}'}}\d_{{\bf c_{-}},{\bf c_{-}'}}
\d_{\lambda,\mu}\mod q_{1}^{-1}\A.
\end{align*}
\elem
\bpf
This follows from Proposition \ref{prop-inner-L} and 
Proposition \ref{prop-alg-indep} (2).
\epf

\bdefn
Let $n\geq 0$. 
In view of $\rm{Proposition}\ \ref{prop-alg-indep}\ (4)$, we set
\bal
&\U_{\Z}^{+}(0)=\Z [q_{1},q_{1}^{-1}][P_{1},P_{2},\cdots ]
\subset\U^{+}(0),\tag{1}\\
&\U_{\Z}^{+}(0;n)=\Z [q_{1},q_{1}^{-1}][P_{1},\cdots ,P_{n}]
\subset\U^{+}(0;n).\tag{2}
\end{align*}
We understand that $\U_{\Z}^{+}(0;0)=\Z [q_{1},q_{1}^{-1}]$.
\edefn

It follows from Lemma \ref{lem-det} that $\ps_{n}\in\U_{\Z}^{+}(0;n)$ 
for $n\geq 1$. 

\blem\label{lem-z-Schur}
The ${\bf S}_{\lambda}$ form a $\Z [q_{1},q_{1}^{-1}]$-basis 
of $\U^{+}_{\Z}(0)$.
\elem
\bpf
In the proof of Proposition \ref{prop-alg-indep}, 
it is shown that the ${\bf S}_{\lambda}$ form a $\Z$-basis of $S$; 
thus, they span $\U^{+}_{\Z}(0)$ over $\Z [q_{1},q_{1}^{-1}]$.
Their linearly independence over $\Z [q_{1},q_{1}^{-1}]$ 
follows from Proposition \ref{prop-alg-indep} (3).
\epf

\bdefn
We set $V_{\Z}=\U^{+}_{\Z}(>)\U^{+}_{\Z}(0)\U^{+}_{\Z}(<)\subset\U^{+}$.
\edefn

It is clear that 
the basis in Theorem \ref{thm-q-basis} (2) is 
a $\Z [q_{1},q_{1}^{-1}]$-basis of $V_{\Z}$.
By Lemma \ref{lem-z-Schur},  
the basis in Theorem \ref{thm-q-basis} (3) is also 
a $\Z [q_{1},q_{1}^{-1}]$-basis of $V_{\Z}$.
We are going to prove that 
$\U_{\Z}^{+}(0)\subset\U_{\Z}^{+}$ (hence $V_{\Z}\subset\U_{\Z}^{+}$) 
and that $V_{\Z}$ is closed under multiplication (Section 8). 
As a result of these, we conclude that $V_{\Z}=\U_{\Z}^{+}$, 
since $V_{\Z}$ contains the generators of 
the $\Z [q_{1},q_{1}^{-1}]$-algebra $\U_{\Z}^{+}$; 
thus, we obtain $\Z [q_{1},q_{1}^{-1}]$-bases of $\U_{\Z}^{+}$.


\section{Key Proposition}

In this section, $\equiv$ means the congruence  
modulo $\U^{0}\U^{+}(<)_{\leq -2}\otimes\U^{+}$,
unless otherwise stated.

\bdefn\label{defn-D(-)}
For $s,t\geq 0$, we define the elements $D^{-}_{s\a_{0}+t\a_{1}}$ 
of $\U^{+}_{\Z}(<)\cap\U^{+}_{s\a_{0}+t\a_{1}}$ 
by the induction on $s$ as follows$:$
\bal
D^{-}_{s\a_{0}+t\a_{1}}
&=\sum_{k=[(t+2)/2]}^{\min(t-1,s)}\sum_{l=t-k+1}^{\min(3(t-k),k)}
T_{\pa}(D^{-}_{(t-k)\a_{0}+(3(t-k)-l)\a_{1}})
E_{\d -\a_{1}}^{(k-l)}E_{\a_{0}}^{(s-k)}\\
&\ \times q_{1}^{-(3(t-k)-l)(k-l)}q^{-2(t-k)(s-k)}
\th(s\geq 2)\th(3\leq t\leq 2s-1)\\
&\ +E_{\d -\a_{1}}^{(t)}E_{\a_{0}}^{(s-t)}\th (s\geq t).
\end{align*}
\edefn

Note that if $s\geq 2,\ 3\leq t\leq 2s-1,\ [(t+2)/2]\leq k\leq \min(t-1,s),
\ t-k+1\leq l\leq \min(3(t-k),k)$, 
then we have $0\leq t-k\leq s-1,\ 3(t-k)-l\geq 0$; 
hence, the defining process by the induction on $s$ works.

\blem \label{lem-D(-)}
Let $s,t\geq 0$. Then,
\bit
\item[(1)]
$D^{-}_{s\a_{0}}=E_{\a_{0}}^{(s)}$,
\item[(2)]
$D^{-}_{t\a_{1}}=\d_{t,0}$,
\item[(3)]
$D^{-}_{s\a_{0}+t\a_{1}}=0$ if $s\geq 1,\ t\geq 2s$,
\item[(4)]
$D^{-}_{s\a_{0}+\a_{1}}=E_{\d -\a_{1}}E_{\a_{0}}^{(s-1)}$ if $s\geq 1$.
\item[(5)]
$D^{-}_{s\a_{0}+(2s-1)\a_{1}}=E_{s\d -\a_{1}}$ if $s\geq 1$.
\eit
\elem
\bpf
This follows from Definition \ref{defn-D(-)}:
(1)--(4) are clear, and (5) is checked by the induction on $s$.
\epf

\bdefn
For $s\geq 1,\ 1\leq t \leq 2s-1,\ 1\leq p\leq [(t+1)/2]$, we define 
the elements $d_{s\a_{0}+t\a_{1},p}$ of $\U^{+}_{\Z}(0)\U^{+}_{\Z}(<)
\cap\U^{+}_{(s-p)\a_{0}+(t-2p+1)\a_{1}}$ by 
\begin{align*}
d_{s\a_{0}+t\a_{1},p}=&\sum_{i=0}^{[(t-2p+1)/2]}
\ps_{i}D^{-}_{(s-p-i)\a_{0}+(t-2p-2i+1)\a_{1}}
 q_{1}^{(-2s+t)(2p+2i)+t+1}(q_{1}-q_{1}^{-1}).
\end{align*}
\edefn

\bex
We have $d_{\a_{0}+\a_{1},1}=1$.
\eex

Let us study some properties of $D^{-}_{s\a_{0}+t\a_{1}}$.

\bprop\label{prop-D(-)}
Let $s\geq 1$ and let $0\leq t\leq 2s-1$. Then,
\bal
&D^{-}_{s\a_{0}+t\a_{1}}[t]_{1}
=\sum_{i=1}^{[(t+1)/2]}E_{i\d -\a_{1}}
D^{-}_{(s-i)\a_{0}+(t-2i+1)\a_{1}}[2i-1]_{1}q^{(-2s+t+1)(i-1)}\th(t\geq 1),
\tag{1}\\
&[D^{-}_{s\a_{0}+t\a_{1}},E_{\a_{1}}]_{q^{-2s+t}}\tag{2}\\
&=D^{-}_{s\a_{0}+(t+1)\a_{1}}[t+1]_{1}
+\sum_{i=1}^{[(t+1)/2]}\ps_{i}D^{-}_{(s-i)\a_{0}+(t-2i+1)\a_{1}}
q_{1}^{(-2s+t)(2i+1)+2s+1}\th (t\geq 1),\\
&\D (D^{-}_{s\a_{0}+t\a_{1}})\tag{3}\\
&\equiv 
c^{s}k_{1}^{-2s+t}\otimes D^{-}_{s\a_{0}+t\a_{1}}
+\sum_{p=1}^{[(t+1)/2]}c^{s-p}k_{1}^{-2s+t+1}E_{p\d -\a_{1}}\otimes 
d_{s\a_{0}+t\a_{1},p}\th(t\geq 1),\\
&\D(T_{\pa}(D^{-}_{s\a_{0}+t\a_{1}}))\tag{4}\\
&\equiv 
c^{3s-t}k_{1}^{-2s+t}\otimes T_{\pa}(D^{-}_{s\a_{0}+t\a_{1}})\\
&\ +\sum_{p=2}^{[(t+3)/2]}c^{3s-t-p}k_{1}^{-2s+t+1}E_{p\d -\a_{1}}
\otimes T_{\pa}(d_{s\a_{0}+t\a_{1},p-1})\th (t\geq 1)\\
&\ +c^{3s-t-1}k_{1}^{-2s+t+1}E_{\d -\a_{1}}\otimes T_{\pa}
([D^{-}_{s\a_{0}+t\a_{1}},E_{\a_{1}}]_{q^{-2s+t}})(q_{1}-q_{1}^{-1}),\\
&[E_{\d -\a_{1}},T_{\pa}(D^{-}_{s\a_{0}+t\a_{1}})]_{q^{2s-t}}\tag{5}\\
&=T_{\pa}(D^{-}_{s\a_{0}+(t-1)\a_{1}})[4s-t+1]_{1}\th (t\geq 1)\\
&\ -\sum_{i=1}^{[t/2]}E_{(i+1)\d -\a_{i}}
T_{\pa}(D^{-}_{(s-i)\a_{0}+(t-2i)\a_{1}})[2i+1]_{1}q^{(-2s+t)i}\th (t\geq 2).
\end{align*}
\eprop

\bpf
We denote by (a)$_{r}$ the statement (a) for $s=1,\ldots,r$.
We prove (1)$_{s}$--(5)$_{s}$ at once by the induction on $s$.
If $s=1$, we have (1) clearly, 
(2) and (5) by Corollary \ref{c-1}, and 
(3) and (4) by Lemma \ref{lem-copro-E} and Lemma \ref{lem-copro-EE}. 
Now, assume that $s\geq 2$.

(1)$_{s}$\ \ 
We use (5)$_{s-1}$.
We rewrite the the right hand side of (1), which is denoted by RHS.
By Definition \ref{defn-D(-)}, we have
\bal
\text{RHS}&=\sum_{i=1}^{[(t+1)/2]}E_{i\d -\a_{1}}
\Big( \sum _{k=[(t-2i+3)/2]}^{\min(t-2i,s-i)}
\sum_{l=t-2i-k+2}^{\min(3(t-2i-k+1),k)}\\
&\ T_{\pa}(D^{-}_{(t-2i-k+1)\a_{0}+(3(t-2i-k+1)-l)\a_{1}})
E_{\d -\a_{1}}^{(k-l)}E_{\a_{0}}^{(s-i-k)}\\
&\ \times q_{1}^{-(3(t-2i-k+1)-l)(k-l)}q^{-2(t-2i-k+1)(s-i-k)}\\
&\ \times\th(s-i\geq 2)\th(3\leq t-2i+1\leq 2(s-i)-1)\\
&\ +E_{\d -\a_{1}}^{(t-2i+1)}E_{\a_{0}}^{(s-t+i-1)}\th (s-t+i-1\geq 0)
\Big) [2i-1]_{1}q^{(-2s+t+1)(i-1)}.
\end{align*}
By $(5)_{s-1}$, the $i=1$ part is equal to
\bal 
&\sum_{k=[(t+1)/2]}^{\min(t-2,s-1)}\sum_{l=t-k}^{\min(3(t-k-1),k)}
\Big( T_{\pa}(D^{-}_{(t-k-1)\a_{0}+(3(t-k-1)-l)\a_{1}})
E_{\d -\a_{1}}q^{-t+k+l+1}\\
&\ +T_{\pa}(D^{-}_{(t-k-1)\a_{0}+(3(t-k-1)-l-1)\a_{1}})[t-k+l]_{1}
\th(3(t-k-1)-l\geq 1)\\
&\ -\sum_{i=1}^{[(3(t-k-1)-l)/2]}E_{(i+1)\d -\a_{1}}
T_{\pa}(D^{-}_{(t-k-i-1)\a_{0}+(3(t-k-1)-l-2i)\a_{1}})\\
&\ \times[2i+1]_{1}q^{(t-k-l-1)i}\th(3(t-k-1)-l\geq 2)\Big)
E_{\d -\a_{1}}^{(k-l)}E_{\a_{0}}^{(s-k-1)}\\
&\ \times q_{1}^{-(3(t-k-1)-l)(k-l)}q^{-2(t-k-1)(s-k-1)}
\th(s\geq 3)\th(4\leq t\leq 2s-2)\\
&\ +E_{\d -\a_{1}}^{(t)}E_{\a_{0}}^{(s-t)}[t]_{1}\th(s\geq t).
\end{align*}
Hence, if we set 
\bal
A=&\sum_{k=[(t+1)/2]}^{\min(t-2,s-1)}\sum_{l=t-k}^{\min(3(t-k-1),k)}
 T_{\pa}(D^{-}_{(t-k-1)\a_{0}+(3(t-k-1)-l)\a_{1}})\\
&\ \times E_{\d -\a_{1}}^{(k-l+1)}E_{\a_{0}}^{(s-k-1)}
[k-l+1]_{1}q^{-t+k+l+1} \\
&\ \times q_{1}^{-(3(t-k-1)-l)(k-l)}q^{-2(t-k-1)(s-k-1)}
\th(s\geq 3)\th(4\leq t\leq 2s-2),\\
B=&\sum_{k=[(t+1)/2]}^{\min(t-2,s-1)}\sum_{l=t-k}^{\min(3(t-k-1)-1,k)}
T_{\pa}(D^{-}_{(t-k-1)\a_{0}+(3(t-k-1)-l-1)\a_{1}})
E_{\d -\a_{1}}^{(k-l)}E_{\a_{0}}^{(s-k-1)}\\
&\ \times [t-k+l]_{1}q_{1}^{-(3(t-k-1)-l)(k-l)}q^{-2(t-k-1)(s-k-1)}
\th(s\geq 3)\th(4\leq t\leq 2s-2),\\
C=&-\sum_{k=[(t+1)/2]}^{\min(t-3,s-1)}\sum_{l=t-k}^{\min(3(t-k-1)-2,k)}
\sum_{i=1}^{[(3(t-k-1)-l)/2]}\\
&\ E_{(i+1)\d -\a_{1}}T_{\pa}(D^{-}_{(t-k-i-1)\a_{0}+(3(t-k-1)-l-2i)\a_{1}})
E_{\d -\a_{1}}^{(k-l)}E_{\a_{0}}^{(s-k-1)}\\
&\ \times [2i+1]_{1}q_{1}^{-(3(t-k-1)-l)(k-l)}
q^{-2(t-k-1)(s-k-1)}q^{(t-k-l-1)i}\\
&\ \times\th(s\geq 4)\th(6\leq t\leq 2s-2),\\
D&=\sum_{i=2}^{[(t+1)/2]}\sum_{k=[(t-2i+3)/2]}^{\min(t-2i,s-i)}
\sum_{l=t-2i-k+2}^{\min(3(t-2i+1-k),k)}\\
&\ E_{i\d -\a_{1}}T_{\pa}(D^{-}_{(t-2i-k+1)\a_{0}+(3(t-2i-k+1)-l)\a_{1}})
E_{\d -\a_{1}}^{(k-l)}E_{\a_{0}}^{(s-i-k)}\\
&\ \times [2i-1]_{1}q_{1}^{-(3(t-2i-k+1)-l)(k-l)}q^{-2(t-2i-k+1)(s-i-k)}
q^{(-2s+t+1)(i-1)},\\
E=&\sum_{i=2}^{[(t+1)/2]}E_{i\d -\a_{1}}E_{\d -\a_{1}}^{(t-2i+1)}
E_{\a_{0}}^{(s-t+i-1)}\\
&\ \times [2i-1]_{1}q^{(-2s+t+1)(i-1)}\th(s-t+i-1\geq 0)\th(t\geq 3),
\end{align*} 
then we have 
\bal
\text{RHS}=A+B+C+D+E+E_{\d -\a_{1}}^{(t)}E_{\a_{0}}^{(s-t)}\th(s\geq t)[t]_{1}.
\end{align*} 
We rewrite $A$. Putting $k-l+1=k'-l',s-k-1=s-k'$, we see that 
$l=l',k=k'-1$ and that 
$[(t+1)/2]\leq k\leq \min(t-2,s-1),t-k \leq l\leq \min(3(t-k-1),k)$
is equivalent to 
$[(t+3)/2]\leq k'\leq \min(t-1,s),t-k'+1 \leq l'\leq \min(3(t-k'),k'-1)$.
Hence, 
\bal
A=&\sum_{k=[(t+3)/2]}^{\min(t-1,s)}\sum_{l=t-k+1}^{\min(3(t-k),k-1)}
 T_{\pa}(D^{-}_{(t-k)\a_{0}+(3(t-k)-l)\a_{1}})
 E_{\d -\a_{1}}^{(k-l)}E_{\a_{0}}^{(s-k)}\\
&\ \times [k-l]_{1}q^{-t+k+l}q_{1}^{-(3(t-k)-l)(k-l-1)}q^{-2(t-k)(s-k)}
\th(s\geq 3)\th(4\leq t\leq 2s-2).
\end{align*} 
We rewrite $B$. Putting $k-l=k'-l',\ s-k-1=s-k'$, we see that 
$l=l'-1,k=k'-1$ and that 
$[(t+1)/2]\leq k\leq \min(t-2,s-1),\ t-k \leq l\leq\min(3(t-k-1)-1,k)$
is equivalent to 
$[(t+3)/2]\leq k'\leq \min(t-1,s),\ t-k'+2 \leq l'\leq\min(3(t-k'),k')$.
Hence,
\bal
B=&\sum_{k=[(t+3)/2]}^{\min(t-1,s)}\sum_{l=t-k+2}^{\min(3(t-k),k-1)}
 T_{\pa}(D^{-}_{(t-k)\a_{0}+(3(t-k)-l)\a_{1}})
E_{\d -\a_{1}}^{(k-l)}E_{\a_{0}}^{(s-k)}\\
&\ \times [t-k+l]_{1}q_{1}^{-(3(t-k)-l+1)(k-l)}q^{-2(t-k)(s-k)}
\th(s\geq 3)\th(4\leq t\leq 2s-2).
\end{align*} 
We rewrite $D$. Putting $i=i'+1,\ t-2i-k+1=t-k'-i'-1,\ 
3(t-2i-k+1)-l=3(t-k'-1)-l'-2i'$, we see that 
$k=k'-i',\ l=l'-i'$ and that 
$2\leq i\leq [(t-2)/2],\ [(t-2i+3)/2]\leq k\leq\min(t-2i,s-i),
\ t-2i-k+2 \leq l\leq\min(3(t-2i-k+1),k)$
is equivalent to 
$[(t+1)/2]\leq k'\leq\min(t-3,s-1),\ t-k'\leq l'\leq\min(3(t-k'-1)-2,k'),
\ 1\leq i'\leq [(3(t-k'-1)-l')/2]$.
Hence,
\bal
D&=\sum_{k=[(t+1)/2]}^{\min(t-3,s-1)}\sum_{l=t-k}^{\min(3(t-k-1)-2,k)}
\sum_{i=1}^{[(3(t-k-1)-l)/2]}\\
&\ E_{(i+1)\d -\a_{1}}T_{\pa}(D^{-}_{(t-k-i-1)\a_{0}+(3(t-k+1)-l-2i)\a_{1}})
E_{\d -\a_{1}}^{(k-l)}E_{\a_{0}}^{(s-k-1)}\\
&\ \times [2i+1]_{1}q_{1}^{-(3(t-k-i-1)-l+i)(k-l)}q^{-2(t-k-i-1)(s-k-1)}
q^{(-2s+t+1)i}\\
&\ \times \th(s\geq 4)\th(6\leq t\leq 2s-2),
\end{align*} 
which cancels out with $C$.
We rewrite $E$. Putting $t-2i+1=k'-l',\ s-t+i-1=s-k'$, we see that
$i=l'=t-k'+1$ and that $2\leq i\leq [(t+1)/2],\  s-t+i-1$ is equivalent to
$[(t+2)/2]\leq k'\leq\min(t-1,s)$. 
Hence,
\bal
E&=\sum_{k=[(t+2)/2]}^{\min(t-1,s)}
T_{\pa}(D^{-}_{(t-k)\a_{0}+(3(t-k)-l)\a_{1}})E_{\d -\a_{1}}^{(k-l)}
E_{\a_{0}}^{(s-k)}\\
&\ \times [2l-1]_{1}q^{(l-1)(-2s+t+1)}\th(k+l=t+1)
\th(s\geq 2)\th(3\leq t\leq 2s-1).
\end{align*} 
It is now easy to see that 
\bal
A+B+E&=\sum_{k=[(t+2)/2]}^{\min(t-1,s)}\sum_{l=t-k+1}^{\min(3(t-k),k)}
T_{\pa}(D^{-}_{(t-k)\a_{0}+(3(t-k)-l)\a_{1}})
E_{\d -\a_{1}}^{(k-l)}E_{\a_{0}}^{(s-k)}\\
&\ \times q_{1}^{-(3(t-k)-l)(k-l)}q^{-2(t-k)(s-k)}\th(s\geq 2)
\th(3\leq t\leq 2s-1)[t]_{1}. 
\end{align*}
Thus, RHS$=D^{-}_{s\a_{0}+t\a_{1}}[t]_{1}$. We obtain (1) for $s$.

(2)$_{s}$\ \ 
We use (1)$_{s}$.
The case where $t=0$ is directly checked by the induction on $s$ 
in view of Lemma \ref{lem-D(-)} (5). The case where $t=2s-1$ is clear.
We assume that $s\geq 2,\ 1\leq t\leq 2s-2$.
We rewrite the left hand side, which is denoted by LHS. 
By $(1)_{s}$, we have 
\bal
\text{LHS}[t]_{1}&=\sum_{l=1}^{[(t+1)/2]}
[E_{l\d -\a_{1}}D^{-}_{(s-l)\a_{0}+(t-2l+1)\a_{1}},E_{\a_{1}}]_{q^{-2s+t}}
[2l-1]_{1}q^{(-2s+t+1)(l-1)}.
\end{align*} 
For $1\leq l\leq [(t+1)/2]$,
using Lemma \ref{lem-bracket}, 
$(2)_{s-1}$, and Proposition \ref{p-2} ((2), (5)), we have 
\bal
&[E_{l\d -\a_{1}}D^{-}_{(s-l)\a_{0}+(t-2l+1)\a_{1}},E_{\a_{1}}]_{q^{-2s+t}}\\
&=E_{l\d -\a_{1}}
[D^{-}_{(s-l)\a_{0}+(t-2l+1)\a_{1}},E_{\a_{1}}]_{q^{-2s+t+1}}\\
&\ +[E_{l\d -\a_{1}},E_{\a_{1}}]_{q^{-1}}
D^{-}_{(s-l)\a_{0}+(t-2l+1)\a_{1}}q^{-2s+t+1}\\
&=\ps_{l}D^{-}_{(s-l)\a_{0}+(t-2l+1)\a_{1}}q^{-2s+t+1}
+E_{l\d -\a_{1}}D^{-}_{(s-l)\a_{0}+(t-2l+2)\a_{1}}[t-2l+2]_{1}\\
&\ +\sum_{i=1}^{[(t-2l+2)/2]}\Big(\ps_{i}E_{l\d -\a_{1}}
+(1+q^{-1})\sum_{j=0}^{i-1}\ps_{j}E_{(l+i-j)\d -\a_{1}}b_{2(i-j)}\Big)\\
&\ \times D^{-}_{(s-l-i)\a_{0}+(t-2l-2i+2)\a_{1}}
q_{1}^{(-2s+t+1)(2i+1)+2(s-l)+1}\th(t\geq 2l).
\end{align*} 
Hence, if we set 
\bal
A&=\sum_{l=1}^{[(t+1)/2]}
\ps_{l}D^{-}_{(s-l)\a_{0}+(t-2l+1)\a_{1}}
[2l-1]_{1}q^{(-2s+t+1)l},\\
B&=\sum_{l=1}^{[(t+1)/2]}
E_{l\d -\a_{1}}D^{-}_{(s-l)\a_{0}+(t-2l+2)\a_{1}}
[t-2l+2]_{1}[2l-1]_{1}q^{(-2s+t+1)(l-1)},\\
C&=\sum_{l=1}^{[t/2]}\sum_{i=1}^{[(t-2l+2)/2]}
\ps_{i}E_{l\d -\a_{1}}D^{-}_{(s-l-i)\a_{0}+(t-2l-2i+2)\a_{1}}\\
&\ \times [2l-1]_{1}q^{(-2s+t+1)(l+i-1)}q_{1}^{t-2l+2},\\
D&=\sum_{l=1}^{[t/2]}\sum_{i=1}^{[(t-2l+2)/2]}\sum_{j=0}^{i-1}
\ps_{j}E_{(l+i-j)\d -\a_{1}}D^{-}_{(s-l-i)\a_{0}+(t-2l-2i+2)\a_{1}}\\
&\ \times b_{2(i-j)}[2l-1]_{1}(1+q^{-1})q^{(-2s+t+1)(l+i-1)}q_{1}^{t-2l+2},
\end{align*} 
then we have 
\bal
\text{LHS}[t]_{1}=A+B+C+D.
\end{align*} 
In the above expression of $B$, on account of the factor $[t-2l+2]_{1}$, 
we can include the case where $l=[(t+2)/2]$.
We rewrite $C$. Since $1\leq l\leq [t/2],\ 1\leq i\leq [(t-2l+2)/2]$
is equivalent to $1\leq i\leq [t/2],\ 1\leq l\leq [(t-2i+2)/2]$, we have
\bal
C&=\sum_{i=1}^{[t/2]}\sum_{l=1}^{[(t-2i+2)/2]}
\ps_{i}E_{l\d -\a_{1}}D^{-}_{(s-l-i)\a_{0}+(t-2i-2l+2)\a_{1}}\\
&\ \times [2l-1]_{1}q^{(-2s+t+1)(i+l-1)}q_{1}^{t-2l+2}\th(t\geq 2).
\end{align*}  
We rewrite $D$. Putting $j=i',\ l+i-j=l',\ i-j=j'$, we see that 
$i=i'+j',\ l=l'-j'$ and that 
$1\leq l\leq [t/2],\ 1\leq i\leq [(t-2l+2)/2],\ 0\leq j\leq i-1$ 
is equivalent to
$0\leq i'\leq [(t-2)/2],\ 2\leq l'\leq [(t-2i'+2)/2],\ 1\leq j\leq l'-1$.
Hence, using Lemma \ref{b-4} (1), we have 
\bal
D&=\sum_{i=0}^{[(t-2)/2]}\sum_{l=2}^{[(t-2i+2)/2]}\sum_{j=1}^{l-1}
\ps_{i}E_{l\d -\a_{1}}D^{-}_{(s-l-i)\a_{0}+(t-2i-2l+2)\a_{1}}\\
&\ \times [2l-2j-1]_{1}(1+q^{-1})b_{2j}
q^{(-2s+t+1)(i+l-1)+j}q_{1}^{t-2l+2}\th(t\geq 2)\\
&=\sum_{i=0}^{[(t-2)/2]}\sum_{l=2}^{[(t-2i+2)/2]}
\ps_{i}E_{l\d -\a_{1}}D^{-}_{(s-l-i)\a_{0}+(t-2i-2l+2)\a_{1}}\\
&\ \times [2l-2]_{1}[2l-1]_{1}(q_{1}-q_{1}^{-1})
q^{(-2s+t+1)(i+l-1)}q_{1}^{t}\th(t\geq 2).
\end{align*}
Here, on account of the factor $[2l-2]_{1}$, 
we can include the case where $l=1$; 
hence, we can also include the case where $i=[t/2]$ 
and thus the case where $t=1$. Thus,
\bal
D&=\sum_{i=0}^{[t/2]}\sum_{l=1}^{[(t-2i+2)/2]}
\ps_{i}E_{l\d -\a_{1}}D^{-}_{(s-l-i)\a_{0}+(t-2i-2l+2)\a_{1}}\\
&\ \times [2l-1]_{1}
q^{(-2s+t+1)(i+l-1)}q_{1}^{t}(q_{1}^{2l-2}-q_{1}^{-2l+2}),
\end{align*}
the $i=0$ part of which is denoted by $E$.
Then,
\bal
C+(D-E)&=\sum_{i=1}^{[t/2]}\sum_{l=1}^{[(t-2i+2)/2]}
\ps_{i}E_{l\d -\a_{1}}D^{-}_{(s-l-i)\a_{0}+(t-2i-2l+2)\a_{1}}\\
&\ \times [2l-1]_{1}q^{(-2s+t+1)(i+l-1)}q_{1}^{t+2l-2}\th(t\geq 2).
\end{align*}
Applying $(1)_{s-1}$, we have 
\bal
C+(D-E)=\sum_{i=1}^{[t/2]}\ps_{i}D^{-}_{(s-i)\a_{0}+(t-2i+1)\a_{1}}
[t-2i+1]_{1}q_{1}^{t-2i},
\end{align*}
where the range of $i$ can be changed to $1\leq i\leq [(t+1)/2]$
on account of the factor $[t-2i+1]_{1}$.
Thus, 
\bal
A+C+(D-E)=\sum_{i=1}^{[(t+1)/2]}\ps_{i}D^{-}_{(s-i)\a_{0}+(t-2i+1)\a_{1}}
[t]_{1}q^{(-2s+t)i}q_{1}^{t+1}. 
\end{align*}
On the other hand, by $(1)_{s}$, we have 
\bal
B+E
&=\sum_{l=1}^{[(t+2)/2]}E_{l\d -\a_{1}}D^{-}_{(s-l)\a_{0}+(t-2l+2)\a_{1}}
[2l-1]_{1}[t]_{1}q^{(-2s+t+2)(l-1)}\\
&=D^{-}_{s\a_{0}+(t+1)\a_{1}}[t+1]_{1}[t]_{1}.
\end{align*}
Hence, we obtain (2) for $s$.


(3)$_{s}$\ \ 
We use (1)$_{s}$.
We can assume that $t\geq 1$. By (1)$_{s}$, we have 
\bal
\D(D^{-}_{s\a_{0}+t\a_{1}})[t]_{1}
&=\sum_{l=1}^{[(t+1)/2]}\D(E_{l\d -\a_{1}}D^{-}_{(s-l)\a_{0}+(t-2l+1)\a_{1}})\\
&\ \times[2l-1]_{1}q^{(-2s+t+1)l}(q_{1}-q_{1}^{-1}).
\end{align*}
For $1\leq l\leq [(t+1)/2]$, using Lemma \ref{lem-copro-E}, (3)$_{s-1}$, 
and Corollary \ref{cor-convex-1}, we have
\bal
&\D(E_{l\d -\a_{1}}D^{-}_{(s-l)\a_{0}+(t-2l+1)\a_{1}})\\
&\equiv 
\Big(c^{l}k_{1}^{-1}\otimes E_{l\d -\a_{1}}
+(q_{1}-q_{1}^{-1})\sum_{i=1}^{l}c^{l-i}E_{i\d -\a_{1}}\otimes\ps_{l-i}\Big)\\
&\ \times \Big( c^{s-l}k_{1}^{-2s+t+1}
\otimes D^{-}_{(s-l)\a_{0}+(t-2l+1)\a_{1}}\\
&\ +\sum_{p=1}^{[(t-2l+2)/2]}c^{s-l-p}k_{1}^{-2s+t+2}E_{p\d -\a_{1}}\otimes 
d_{(s-l)\a_{0}+(t-2l+1)\a_{1},p}\th(t\geq 2l)\Big)\\
&\equiv 
c^{s}k_{1}^{-2s+t}\otimes E_{l\d -\a_{1}}D^{-}_{(s-l)\a_{0}+(t-2l+1)\a_{1}}\\
&\ +\sum_{i=1}^{l}c^{s-i}k_{1}^{-2s+t+1}E_{i\d -\a_{1}}
\otimes\ps_{l-i}D^{-}_{(s-l)\a_{0}+(t-2l+1)\a_{1}}q^{-2s+t+1}
(q_{1}-q_{1}^{-1})\\
&\ +\sum_{p=1}^{[(t-2l+2)/2]}c^{s-p}k_{1}^{-2s+t+1}E_{p\d -\a_{1}}\otimes 
E_{l\d -\a_{1}}d_{(s-l)\a_{0}+(t-2l+1)\a_{1},p}\th(t\geq 2l).
\end{align*}
Hence, if we set 
\bal
A&=\sum_{l=1}^{[(t+1)/2]}\sum_{i=1}^{l}
c^{s-i}k_{1}^{-2s+t+1}E_{i\d -\a_{1}}\otimes 
\ps_{l-i}D^{-}_{(s-l)\a_{0}+(t-2l+1)\a_{1}}\\
&\ \times [2l-1]_{1}q^{(-2s+t+1)l}(q_{1}-q_{1}^{-1}),\\
B&=\sum_{l=1}^{[t/2]}\sum_{p=1}^{[(t-2l+2)/2]}
c^{s-p}k_{1}^{-2s+t+1}E_{p\d -\a_{1}}\otimes 
E_{l\d -\a_{1}}d_{(s-l)\a_{0}+(t-2l+1)\a_{1},p}\\
&\ \times [2l-1]_{1}q^{(-2s+t+1)(l-1)}\th(t\geq 2),
\end{align*}
then we have
\bal
\D(D^{-}_{s\a_{0}+t\a_{1}})[t]_{1}\equiv 
c^{s}k_{1}^{-2s+t}\otimes D^{-}_{s\a_{0}+t\a_{1}}[t]_{1}+A+B.
\end{align*}
We rewrite $A$. Putting $i=p,\ l-i=i'$, we see that $l=p+i'$ and that 
$1\leq l\leq [(t+1)/2],\ 1\leq i\leq l$ is equivalent to
$1\leq p\leq [(t+1)/2],\ 0\leq i\leq [(t-2p+1)/2]$. 
Hence, 
\bal
A&=\sum_{p=1}^{[(t+1)/2]}c^{s-p}k_{1}^{-2s+t+1}E_{p\d -\a_{1}}\otimes A_{p}
\end{align*}
where we set
\bal
A_{p}&= \sum_{i=0}^{[(t-2p+1)/2]}\ps_{i}
D^{-}_{(s-p-i)\a_{0}+(t-2p-2i+1)\a_{1}}
[2p+2i-1]_{1}q^{(-2s+t+1)(p+i)}(q_{1}-q_{1}^{-1}).
\end{align*}
We rewrite $B$. 
Since $1\leq l\leq [t/2],\ 1\leq p\leq [(t-2l+2)/2]$ 
is equivalent to $1\leq p\leq [t/2],\ 1\leq l\leq [(t-2p+2)/2]$, we have 
\bal
B=\sum_{p=1}^{[t/2]}
c^{s-p}k_{1}^{-2s+t+1}E_{p\d -\a_{1}}\otimes B_{p}
\end{align*}
where we set
\bal
B_{p}&=\sum_{l=1}^{[(t-2p+2)/2]}E_{l\d -\a_{1}}
d_{(s-l)\a_{0}+(t-2l+1)\a_{1},p}
[2l-1]_{1}q^{(-2s+t+1)(l-1)}\th(t\geq 2).
\end{align*}
Note that if $t$ is odd and $p=(t+1)/2$, then
$A_{p}=D^{-}_{(s-p)\a_{0}}[t]_{1}q^{(-2s+t+1)p}=d_{s\a_{0}+t\a_{1},p}[t]_{1}$.
Hence, the proof of (3) for $s$ is reduced
to showing that $A_{p}+B_{p}=d_{s\a_{0}+t\a_{1}}[t]_{1}$ 
for $t\geq 2,\ 1\leq p\leq [t/2]$. 
We rewrite $B_{p}$. 
Using Proposition \ref{p-2} (5), we have
\bal
B_{p}
&=\sum_{l=1}^{[(t-2p+2)/2]}\sum_{i=0}^{[(t-2p-2l+2)/2]}
E_{l\d -\a_{1}}\ps_{i}D^{-}_{(s-l-p-i)\a_{0}+(t-2p-2l-2i+2)\a_{1}}\\
&\ \times [2l-1]_{1}q_{1}^{(-2s+t+1)(2p+2i)+t-2l+2}
q^{(-2s+t+1)(l-1)}(q_{1}-q_{1}^{-1})\\
&=\sum_{l=1}^{[(t-2p+2)/2]}\sum_{i=0}^{[(t-2p-2l+2)/2]}
\Big(\ps_{i}E_{l\d -\a_{1}}\\
&\ +(1+q^{-1})\sum_{j=0}^{i-1}b_{2(i-j)}
\ps_{j}E_{(l+i-j)\d -\a_{1}}\th(i\geq 1)\Big)\\
&\ \times D^{-}_{(s-p-l-i)\a_{0}+(t-2p-2l-2i+2)\a_{1}}[2l-1]_{1}
q^{(-2s+t+1)(p+l+i-1)}q_{1}^{t-2l+2}(q_{1}-q_{1}^{-1}).
\end{align*}
In the double summation, 
$1\leq l\leq [(t-2p+2)/2],\ 0\leq i\leq [(t-2p-2l+2)/2]$ is equivalent to
$0\leq i\leq [(t-2p)/2],\ 1\leq l\leq [(t-2p-2i+2)/2]$.
In the triple one, putting $j=i',\ l+i-j=l',\ i-j=j'$, 
we see that $i=i'+j',\ l=l'-j'$ and that
$1\leq l\leq [(t-2p+2)/2],\ 1\leq i\leq [(t-2l-2p+2)/2],\ 0\leq j\leq i-1$ 
is equivalent to 
$0\leq i'\leq [(t-2p-2)/2],\ 2\leq l'\leq [(t-2p-2i'+2)/2],\ 1\leq j\leq l'-1$.
Thus,
\bal
B_{p}
&=\sum_{i=0}^{[(t-2p)/2]}\sum_{l=1}^{[(t-2p-2i+2)/2]}
\ps_{i}E_{l\d -\a_{1}}D^{-}_{(s-l-p-i)\a_{0}+(t-2l-2p-2i+2)\a_{1}}\\
&\ \times [2l-1]_{1}q^{(-2s+t+1)(p+l+i-1)}q_{1}^{t-2l+2}(q_{1}-q_{1}^{-1})\\
&\ +\sum_{i=0}^{[(t-2p-2)/2]}
\sum_{l=2}^{[(t-2p-2i)/2]}\sum_{j=1}^{l-1}
\ps_{i}E_{l\d -\a_{1}}D^{-}_{(s-p-l-i)\a_{0}+(t-2p-2l-2i+2)\a_{1}}\\
&\ \times [2l-2j-1]_{1}b_{2j}q^{(-2s+t+1)(p+i+l-1)+j}q_{1}^{t-2l+2}
(1+q^{-1})(q_{1}-q_{1}^{-1}).
\end{align*} 
By Lemma \ref{b-4}, the latter term is equal to 
\bal
&\sum_{i=0}^{[(t-2p-2)/2]}\sum_{l=2}^{[(t-2p-2i)/2]}
\ps_{i}E_{l\d -\a_{1}}D^{-}_{(s-p-l-i)\a_{0}+(t-2p-2l-2i+2)\a_{1}}\\
&\ \times [2l-1]_{1}(q^{l-1}-q^{-l+1})q^{l-1}
q^{(-2s+t+1)(p+i+l-1)}q_{1}^{t-2l+2}(q_{1}-q_{1}^{-1}).
\end{align*}
Here, on account of the factor $(q^{l-1}-q^{-l+1})$, we can include the case 
where $l=1$; hence, we can also include the case where 
$i=[(t-2p-2)/2]+1=[(t-2p)/2]$. Thus, 
\bal
B_{p}&=\sum_{i=0}^{[(t-2p)/2]}\sum_{l=1}^{[(t-2p-2i+2)/2]}
\ps_{i}E_{l\d -\a_{1}}D^{-}_{(s-p-l-i)\a_{0}+(t-2p-2l-2i+2)\a_{1}}\\
&\ \times [2l-1]_{1}q^{(-2s+t+1)(p+i+l-1)}q_{1}^{t+2l-2}(q_{1}-q_{1}^{-1}).
\end{align*}
Applying (1)$_{s-1}$, we have 
\bal
B_{p}&=\sum_{i=0}^{[(t-2p)/2]}
\ps_{i}D^{-}_{(s-p-i)\a_{0}+(t-2p-2i+1)\a_{1}}\\
&\ \times [t-2p-2i+1]_{1}q^{(-2s+t+1)(p+i)}q_{1}^{t}(q_{1}-q_{1}^{-1}),
\end{align*}
where the range of $i$ can be changed to $0\leq i\leq [(t-2p+1)/2]$ 
on account of the factor $[t-2p-2i+1]_{1}$.
Hence, 
\bal
A_{p}+B_{p}&=
\sum_{i=0}^{[(t-2p+1)/2]}
\ps_{i}D^{-}_{(s-p-i)\a_{0}+(t-2p-2i+1)\a_{1}}\\
&\ \times ([2p+2i-1]_{1}+[t-2p-2i+1]_{1}q_{1}^{t})q^{(-2s+t+1)(p+i)}
(q_{1}-q_{1}^{-1}),
\end{align*}
which is equal to $d_{s\a_{0}+t\a_{1},p}[t]_{1}$.
We obtain (3) for $s$.

(4)$_{s}$\ \ 
We use (3)$_{s}$.
By Lemma \ref{lem-copro-braid} and (3)$_{s}$, we have
\bal
&\D(T_{\pa}(D^{-}_{s\a_{0}+t\a_{1}}))\\
&\equiv  
T_{\pa}(c^{s}k_{1}^{-2s+t})\otimes T_{\pa}(D^{-}_{s\a_{0}+t\a_{1}})\\
&\ +\sum_{p=1}^{[(t+1)/2]}
T_{\pa}(c^{s-p}k_{1}^{-2s+t+1}E_{p\d -\a_{1}})
\otimes T_{\pa}(d_{s\a_{0}+t\a_{1},p})\th(t\geq 1)\\
&\ +(q_{1}-q_{1}^{-1})\big[T_{\pa}(c^{s}k_{1}^{-2s+t})\otimes
T_{\pa}(D^{-}_{s\a_{0}+t\a_{1}}),
c^{-1}k_{1}E_{\d -\a_{1}}\otimes T_{\pa}(E_{\a_{1}})\big]\\
&\equiv  
c^{3s-t}k_{1}^{-2s+t}\otimes T_{\pa}(D^{-}_{s\a_{0}+t\a_{1}})\\
&\ +\sum_{p=2}^{[(t+3)/2]}
c^{3s-t-p}k_{1}^{-2s+t+1}E_{p\d -\a_{1}}
\otimes T_{\pa}(d_{s\a_{0}+t\a_{1},p-1})\th(t\geq 1)\\
&\ +(q_{1}-q_{1}^{-1})c^{3s-t-1}k_{1}^{-2s+t+1}E_{\d -\a_{1}}
\otimes T_{\pa}([D^{-}_{s\a_{0}+t\a_{1}},E_{\a_{1}}]_{q^{-2s+t}})
\end{align*}
where $\equiv$ means the congruence modulo 
$\U^{0}\U^{+}(<)_{\leq -2}\otimes\U$. 
Noting that both sides belong to $\U^{0}\U^{+}(<)\otimes\U^{+}$ 
by Corollary \ref{cor-copro-E(<)} and Lemma \ref{lem-convex-E}, 
we obtain (4) for $s$.

(5) for $s$ with $t=2s-1$\ \ 
We use (1)$_{s}$.
We rewrite the right hand side of (5) for $s$ with $t=2s-1$, 
which is denoted by RHS. Using Lemma \ref{lem-D(-)} (2), we have 
\bal
\text{RHS}&=T_{\pa}(D^{-}_{s\a_{0}+(2s-2)\a_{1}})[2s+2]_{1}
-\sum_{l=1}^{s-1}E_{(l+1)\d -\a_{1}}E_{(s-l+1)\d -\a_{1}}[2l+1]_{1}q^{-l}.
\end{align*} 
By $(1)_{s}$ and Lemma \ref{lem-D(-)} (2), we have 
\bal
T_{\pa}(D^{-}_{s\a_{0}+(2s-2)\a_{1}})[2s-2]_{1}
=\sum_{l=1}^{s-1}E_{(l+1)\d -\a_{1}}E_{(s-l+1)\d -\a_{1}}[2l-1]_{1}q^{-l+1}.
\end{align*} 
Hence,
\bal
\text{RHS}[2s-2]_{1}
=\sum_{l=1}^{s-1}E_{(l+1)\d -\a_{1}}E_{(s-l+1)\d -\a_{1}}c_{s,l}
\end{align*}
where we set
\bal
c_{s,l}=([2s+2]_{1}[2l-1]_{1}q-[2s-2]_{1}[2l+1]_{1})q^{-l}
\end{align*}
for $s,l\in\Z$. 
Using Lemma \ref{b-4}, for $s\in\Z,\ l\geq 1$, we have 
\bal
\sum_{i=1}^{l-1}b_{2i}c_{s,l-i}\th(l\geq 2)
&=[2s+2]_{1}[2l-1]_{1}[l-1](q-1)\\
&\ -[2s-2]_{1}[2l+1]_{1}[l](q-1)+b_{2l}[2s-2]_{1}.
\end{align*}
Hence, for $s\in\Z,\ l\geq 1$, we have 
\bal
c_{s,s-l}+c_{s,l}q+\sum_{i=1}^{l-1}b_{2i}c_{s,l-i}(1+q)\th(l\geq 2)
=(1+q)b_{2l}[2s-2]_{1}\tag{$\star$}.
\end{align*}
Using Lemma \ref{b-4}, we also have 
\bal
\sum_{i=1}^{m}b_{2i-1}c_{2m+1,m-i+1}=b_{2m+1}[4m]_{1}
\ \ \text{for}\ \ m\geq 1\tag{$\star\star$}.
\end{align*}
Now, we consider the case where $s=2m$ with $m\geq 1$.
Using Corollary \ref{c-2} (2), we have
\bal
\text{RHS}[4m-2]_{1}
&=\sum_{l=m}^{2m-1}E_{(l+1)\d -\a_{1}}E_{(2m-l+1)\d -\a_{1}}c_{2m,l}\\
&\ +\sum_{l=1}^{m-1}
\Big( E_{(2m-l+1)\d -\a_{1}}E_{(l+1)\d -\a_{1}}q\\
&\ +\sum_{i=1}^{m-l-1}E_{(2m-l-i+1)\d -\a_{1}}E_{(l+i+1)\d -\a_{1}}
b_{2i}(1+q)\th(l\leq m-2)\\
&\ +E_{(m+1)\d -\a_{1}}^{2}b_{2(m-l)}\Big) c_{2m,l}\th(m\geq 2).
\end{align*}
In the first term, putting $l=2m-l'$, we see that 
$m\leq l\leq 2m-1$ is equivalent to $1\leq l'\leq m$.
In the double summation, putting $l=l'-i$, we see that 
$1\leq l\leq m-2,\ 1\leq i\leq m-l-1$ is equivalent to 
$2\leq l'\leq m-1,\ 1\leq i\leq l'-1$. Hence,
\bal
\text{RHS}[4m-2]_{1}
&=\th(m\geq 2)\sum_{l=1}^{m-1}E_{(2m-l+1)\d -\a_{1}}E_{(l+1)\d -\a_{1}}\\
&\ \times\Big( c_{2m,2m-l}+c_{2m,l}q
+\sum_{i=1}^{l-1}b_{2i}c_{2m,l-i}(1+q)\th(l\geq 2)\Big)\\
&\ +E_{(m+1)\d -\a_{1}}^{2}
\Big( c_{2m,m}+\sum_{i=1}^{m-1}b_{2i}c_{2m,m-i}\th(m\geq 2)\Big).
\end{align*}
Applying \thetag{$\star$}, we have 
\bal
\text{RHS}
=(1+q)\sum_{l=1}^{m-1}b_{2l}E_{(2m-l+1)\d -\a_{1}}E_{(l+1)\d -\a_{1}}
+b_{2m}E_{(m+1)\d -\a_{1}}^{2},
\end{align*}
which is equal to $[E_{\d -\a_{1}},E_{(2m+1)\d -\a_{1}}]_{q}$ 
by Corollary \ref{c-2} (2).
Next, we consider the case where $s=2m+1$ with $m\geq 1$.
Using Corollary \ref{c-2} (2), we have
\bal
\text{RHS}[4m]_{1}
&=\sum_{l=m+1}^{2m}E_{(l+1)\d -\a_{1}}E_{(2m-l+2)\d -\a_{1}}c_{2m+1,l}\\
&\ +\sum_{l=1}^{m}
\Big( E_{(2m-l+2)\d -\a_{1}}E_{(l+1)\d -\a_{1}}q\\
&\ +\sum_{i=1}^{m-l}E_{(2m-l-i+2)\d -\a_{1}}E_{(l+i+1)\d -\a_{1}}
b_{2i}(1+q)\th(l\leq m-1)\\
&\ +E_{(2m+2)\d -\a_{0}}b_{2(m-l)+1}[4]_{1}\Big) c_{2m+1,l}.
\end{align*}
In the first term, putting $l=2m-l'+1$, we see that 
$m+1\leq l\leq 2m$ is equivalent to $1\leq l'\leq m$.
In the double summation, putting $l=l'-i$, we see that 
$1\leq l\leq m-1,\ 1\leq i\leq m-l$ is equivalent to 
$2\leq l'\leq m,\ 1\leq l'-1$. Hence, 
\bal
\text{RHS}[4m]_{1}
&=\sum_{l=1}^{m}E_{(2m-l+2)\d -\a_{1}}E_{(l+1)\d -\a_{1}}\\
&\ \times \Big( c_{2m+1,2m-l+1}+c_{2m+1,l}q
+\sum_{i=1}^{l-1}b_{2i}c_{2m+1,l-i}(1+q)\th(l\geq 2)\Big)\\
&\ +E_{(2m+2)\d -\a_{0}}\sum_{i=1}^{m}b_{2i-1}c_{2m+1,m-i+1}[4]_{1}.
\end{align*}
Applying \thetag{$\star$} and \thetag{$\star\star$}, we have 
\bal
\text{RHS}=(1+q)\sum_{l=1}^{m}b_{2l}E_{(2m-l+2)\d -\a_{1}}E_{(l+1)\d -\a_{1}}
+[4]_{1}b_{2m+1}E_{(2m+2)\d+\a_{0}},
\end{align*}
which is equal to
$[E_{\d -\a_{1}},E_{(2m+2)\d -\a_{1}}]_{q}$ by Corollary \ref{c-2} (2).
Thus, we obtain (5) for $s$ with $t=2s-1$.


(5)$_{s}$\ \ 
We use (2)$_{s}$ and (4)$_{s}$.
The case where $t=2s-1$ is proved above.
The case where $t=0$ is clear from Lemma \ref{l-6}.
We assume that $s\geq 2,\ 1\leq t\leq 2s-2$ and
argue by the descending induction on $t$.
The left and right hand sides of (5) (for $s$ and $t$) are
denoted by LHS and RHS respectively. 
By Corollary \ref{cor-convex-1}, both of LHS and RHS 
belong to $\U^{+}(<;\d-\a_{1})$: 
thus, we can use Corollary \ref{cor-copro} to prove (5) for $s$. 
Using Lemma \ref{lem-copro-E} with $n=1$ and (4)$_{s}$, we have
\bal
\D(\text{LHS})
&=[\D(E_{\d-\a_{1}}),\D(T_{\pa}(D^{-}_{s\a_{0}+t\a_{1}}))]_{q^{2s-t}}\\
&\equiv 
\Big[ ck_{1}^{-1}\otimes E_{\d -\a_{1}}+E_{\d -\a_{1}}\otimes 1,
\ c^{3s-t}k_{1}^{-2s+t}\otimes T_{\pa}(D^{-}_{s\a_{0}+t\a_{1}})\\
&\ +\sum_{p=2}^{[(t+3)/2]}c^{3s-t-p}k_{1}^{-2s+t+1}E_{p\d -\a_{1}}
\otimes T_{\pa}(d_{s\a_{0}+t\a_{1},p-1})\th (t\geq 1)\\
&\ +c^{3s-t-1}k_{1}^{-2s+t+1}E_{\d -\a_{1}}\otimes T_{\pa}
([D^{-}_{s\a_{0}+t\a_{1}},E_{\a_{1}}]_{q^{-2s+t}})
(q_{1}-q_{1}^{-1})\Big]_{q^{2s-t}}\\
&\equiv 
c^{3s-t+1}k_{1}^{-2s+t-1}\otimes \text{LHS}
+\sum_{p=1}^{[(t+3)/2]}c^{3s-t-p+1}k_{1}^{-2s+t}E_{p\d -\a_{1}}\otimes A_{p}
\end{align*}
where we set
\bal
A_{p}=[E_{\d -\a_{1}},T_{\pa}(d_{s\a_{0}+t\a_{1},p-1})]_{q^{2s-t-1}}
\end{align*}
for $2\leq p\leq [(t+3)/2]$, and
\bal
A_{1}
&=[E_{\d -\a_{1}},T_{\pa}([D^{-}_{s\a_{0}+t\a_{1}},E_{\a_{1}}]_{q^{-2s+t}})]
_{q^{2s-t-1}}(q_{1}-q_{1}^{-1})\\
&\ +T_{\pa}(D^{-}_{s\a_{0}+t\a_{1}})(q^{-2s+t}-q^{2s-t}).
\end{align*}
Note that if $t$ is odd and $p=(t+3)/2$, then $A_{p}=0$,
so that we have $A_{p}$ for $1\leq p\leq [(t+2)/2]$.
Next, we calculate the coproduct of RHS.
By (4)$_{s}$, we have
\bal
\D (\text{RHS})
&=\D\Big( T_{\pa}(D^{-}_{s\a_{0}+(t-1)\a_{1}})[4s-t+1]_{1}\\
&\ -\sum_{i=1}^{[t/2]}E_{(i+1)\d -\a_{i}}
T_{\pa}(D^{-}_{(s-i)\a_{0}+(t-2i)\a_{1}})
[2i+1]_{1}q^{(-2s+t)i}\th (t\geq 2)\Big).
\end{align*}
For $t\geq 2,\ 1\leq i\leq [t/2]$, 
using Lemma \ref{lem-copro-E} and (4)$_{s-1}$, we have 
\bal
&\D\big( E_{(i+1)\d -\a_{1}}T_{\pa}(D^{-}_{(s-i)\a_{0}+(t-2i)\a_{1}})\big)\\
&\equiv
\Big( c^{i+1}k_{1}^{-1}\otimes E_{(i+1)\d -\a_{1}} 
+(q_{1}-q_{1}^{-1})\sum_{p=1}^{i+1}c^{i-p+1}E_{p\d -\a_{1}}
\otimes\ps_{i-p+1}\Big)\\
&\times
\Big( c^{3s-t-i}k_{1}^{-2s+t}\otimes 
T_{\pa}(D^{-}_{(s-i)\a_{0}+(t-2i)\a_{1}})\\
&\ +\sum_{p=2}^{[(t-2i+3)/2]}c^{3s-t-i-p}k_{1}^{-2s+t+1}E_{p\d -\a_{1}}
\otimes T_{\pa}(d_{(s-i)\a_{0}+(t-2i)\a_{1},p-1})\th(t-2i\geq 1)\\
&\ +c^{3s-t-i-1}k_{1}^{-2s+t+1}E_{\d-\a_{1}}\otimes T_{\pa}
([D^{-}_{(s-i)\a_{0}+(t-2i)\a_{1}},E_{\a_{1}}]_{q^{-2s+t}})
(q_{1}-q_{1}^{-1})\Big)\\
&\equiv 
c^{3s-t+1}k_{1}^{-2s+t-1}\otimes 
E_{(i+1)\d -\a_{1}}T_{\pa}(D^{-}_{(s-i)\a_{0}+(t-2i)\a_{1}})\\
&\ +\sum_{p=2}^{[(t-2i+3)/2]}c^{3s-t-p+1}k_{1}^{-2s+t}E_{p\d -\a_{1}}\otimes
E_{(i+1)\d -\a_{1}}T_{\pa}(d_{(s-i)\a_{0}+(t-2i)\a_{1},p-1})\\
&\ +E_{\d -\a_{1}}c^{3s-t}k_{1}^{-2s+t}\otimes 
E_{(i+1)\d -\a_{1}}
T_{\pa}([D^{-}_{(s-i)\a_{0}+(t-2i)\a_{1}},E_{\a_{1}}]_{q^{-2s+t}})
(q_{1}-q_{1}^{-1})q^{2s-t}\\
&\  +\sum_{p=1}^{i+1}k_{1}^{-2s+t}E_{p\d -\a_{1}}\otimes
\ps_{i-p+1}T_{\pa}(D^{-}_{(s-i)\a_{0}+(t-2i)\a_{1}})
q^{-2s+t}(q_{1}-q_{1}^{-1}).
\end{align*}
We sum up this for $1\leq i\leq [t/2]$ after multiplying 
$[2i+1]_{1}q^{(-2s+t)i}$.
Note that $1\leq i\leq [t/2],\ 2\leq p\leq [(t-2i+3)/2]$ is equivalent to
$2\leq p\leq [(t+1)/2],\ 1\leq i\leq [(t-2p+3)/2]$
and that $1\leq i\leq [t/2],\ 1\leq p\leq i+1$ is equivalent to 
$1\leq p\leq [(t+2)/2],\ \max(1,p-1)\leq i\leq [t/2]$. 
Thus, for $t\geq 2$, we have
\bal
&\sum_{i=1}^{[t/2]}
\D\big( E_{(i+1)\d -\a_{1}}T_{\pa}(D^{-}_{(s-i)\a_{0}+(t-2i)\a_{1}})\big)
[2i+1]_{1}q^{(-2s+t)i}\\
&\equiv
c^{3s-t+1}k_{1}^{-2s+t-1}\otimes \sum_{i=1}^{[t/2]}
E_{(i+1)\d -\a_{1}}T_{\pa}(D^{-}_{(s-i)\a_{0}+(t-2i)\a_{1}})
[2i+1]_{1}q^{(-2s+t)i}\\
&\ +\sum_{p=2}^{[(t+1)/2]}
c^{3s-t-p+1}k_{1}^{-2s+t}E_{p\d -\a_{1}}\\
&\ \otimes \sum_{i=1}^{[(t-2p+3)/2]}
E_{(i+1)\d -\a_{1}}T_{\pa}(d_{(s-i)\a_{0}+(t-2i)\a_{1},p-1})
[2i+1]_{1}q^{(-2s+t)i}\\
&\ +E_{\d -\a_{1}}c^{3s-t}k_{1}^{-2s+t}\otimes 
\sum_{i=1}^{[t/2]} E_{(i+1)\d -\a_{1}}
T_{\pa}([D^{-}_{(s-i)\a_{0}+(t-2i)\a_{1}},E_{\a_{1}}]_{q^{-2s+t}})\\
&\ \times [2i+1]_{1}q^{(-2s+t)(i-1)}(q_{1}-q_{1}^{-1})\\
&\ +\sum_{p=1}^{[(t+2)/2]}
k_{1}^{-2s+t}E_{p\d -\a_{1}}\otimes \sum_{i=\max(1,p-1)}^{[t/2]}
\ps_{i-p+1}T_{\pa}(D^{-}_{(s-i)\a_{0}+(t-2i)\a_{1}})\\
&\ \times [2i+1]_{1}q^{(-2s+t)(i+1)}(q_{1}-q_{1}^{-1}).
\end{align*}
Hence,
\bal
\D (\text{RHS})
&\equiv 
c^{3s-t+1}k_{1}^{-2s+t-1}\otimes \text{RHS}
+\sum_{p=1}^{[(t+2)/2]}c^{3s-t-p+1}k_{1}^{-2s+t}E_{p\d -\a_{1}}\otimes B_{p}
\end{align*}
where we set
\bal
B_{p}&=T_{\pa}(d_{s\a_{0}+(t-1)\a_{1},p-1})[4s-t+1]_{1}\\
&\ -\sum_{i=1}^{[(t-2p+3)/2]}E_{(i+1)\d -\a_{1}}
T_{\pa}(d_{(s-i)\a_{0}+(t-2i)\a_{1},p-1})\\
&\ \times [2i+1]_{1}q^{(-2s+t)i}\th(p\leq [(t+1)/2])\\
&\ -\sum_{i=p-1}^{[t/2]}\ps_{i-p+1}T_{\pa}(D^{-}_{(s-i)\a_{0}+(t-2i)\a_{1}})
[2i+1]_{1}q^{(-2s+t)(i+1)}(q_{1}-q_{1}^{-1})
\end{align*}
for $t\geq 2,\ 2\leq p\leq [(t+2)/2]$ and
\bal
B_{1}&=T_{\pa}([D^{-}_{s\a_{0}+(t-1)\a_{1}},E_{\a_{1}}]_{q^{-2s+t-1}})
[4s-t+1]_{1}(q_{1}-q_{1}^{-1})\\
&\ -\sum_{i=1}^{[t/2]}E_{(i+1)\d -\a_{1}}
T_{\pa}([D^{-}_{(s-i)\a_{0}+(t-2i)\a_{1}},E_{\a_{1}}]_{q^{-2s+t-1}})\\
&\ \times [2i+1]_{1}(q_{1}-q_{1}^{-1})q^{(-2s+t)i}\th(t\geq 2)\\
&\ -\sum_{i=1}^{[t/2]}\ps_{i}T_{\pa}(D^{-}_{(s-i)\a_{0}+(t-2i)\a_{1}})
[2i+1]_{1}q^{(-2s+t)(i+1)}(q_{1}-q_{1}^{-1})\th(t\geq 2).
\end{align*}
By virtue of Corollary \ref{cor-copro},
the proof of (5) for $s$ and $t$ is reduced to showing that
$A_{p}=B_{p}$ for $1\leq p\leq [(t+2)/2]$.
First, assuming that $t\geq 2$, we shall prove that 
$A_{p}=B_{p}$ for $2\leq p\leq [(t+2)/2]$.
We rewrite $A_{p}$. Using Proposition \ref{p-2} (5), we have
\bal
A_{p}
&=[E_{\d -\a_{1}},T_{\pa}(d_{s\a_{0}+t\a_{1},p-1})]_{q^{2s-t-1}}\\
&=\Big[ E_{\d -\a_{1}},\sum_{i=0}^{[(t-2p+3)/2]}
\ps_{i}T_{\pa}(D^{-}_{(s-p-i)\a_{0}+(t-2p-2i+3)\a_{1}})\Big]_{q^{2s-t-1}}\\
&\ \times q_{1}^{(-2s+t)(2p+2i-2)+t+1}(q_{1}-q_{1}^{-1})\\
&=\sum_{i=0}^{[(t-2p+3)/2]}\left(\ps_{i}[E_{\d -\a_{1}},
T_{\pa}(D^{-}_{(s-p-i+1)\a_{0}+(t-2p-2i+3)\a_{1}})]_{q^{2s-t-1}}\right.\\
&\ +(1+q^{-1})\sum_{j=0}^{i-1}b_{2(i-j)}\left.\ps_{j}
E_{(i-j+1)\d-\a_{1}}T_{\pa}(D^{-}_{(s-p-i)\a_{0}+(t-2p-2i+3)\a_{1}})\right)\\
&\ \times q_{1}^{(-2s+t)(2p+2i-2)+t+1}(q_{1}-q_{1}^{-1})\th(p\leq [(t+1)/2]).
\end{align*}
We apply (5)$_{s-1}$ to the first term. In the second one, putting
$j=i',\ i-j=j'$, we see that $i=i'+j'$ and that
$1\leq i\leq [(t-2p+3)/2],\ 0\leq j\leq i-1$ is equivalent to
$0\leq i'\leq [(t-2p+1)/2],\ 1\leq j\leq [(t-2p-2i+3)/2]$.
Hence,
\bal
A_{p}
&=\sum_{i=0}^{[(t-2p+3)/2]}\ps_{i}
\Big( T_{\pa}(D^{-}_{(s-p-i+1)\a_{0}+(t-2p-2i+2)\a_{1}})\\
&\ \times [4s-t-2p-2i+2]_{1}\th(t-2p-2i+2\geq 0)\\
&\ -\sum_{j=1}^{[(t-2p-2i+3)/2]}E_{(j+1)\d -\a_{i}}
T_{\pa}(D^{-}_{(s-p-i-j+1)\a_{0}+(t-2p-2i-2j+3)\a_{1}})
[2j+1]_{1}\\
&\ \times q^{(-2s+t+1)j}\th(t-2p-2i+1\geq 0)\Big)
q_{1}^{(-2s+t)(2p+2i-2)+t+1}(q_{1}-q_{1}^{-1})\\
&\ +\sum_{i=0}^{[(t-2p+1)/2]}\sum_{j=1}^{[(t-2p-2i+3)/2]}\ps_{i}
E_{(j+1)\d-\a_{1}}T_{\pa}(D^{-}_{(s-p-i-j+1)\a_{0}+(t-2p-2i-2j+3)\a_{1}})\\
&\ \times b_{2j}q_{1}^{(-2s+t)(2p+2i+2j-1)+2s+1}(1+q^{-1})(q_{1}-q_{1}^{-1})
\th(p\leq [(t+1)/2])\\
&=\sum_{i=0}^{[(t-2p+2)/2]}\ps_{i}
T_{\pa}(D^{-}_{(s-p-i+1)\a_{0}+(t-2p-2i+2)\a_{1}})[4s-t-2p-2i+2]_{1}\\
&\ \times q_{1}^{(-2s+t)(2p+2i-1)+2s+1}(q_{1}-q_{1}^{-1})\\
&\ +\sum_{i=0}^{[(t-2p+1)/2]}\sum_{j=1}^{[(t-2p-2i+3)/2]}\ps_{i}
E_{(j+1)\d-\a_{1}}T_{\pa}(D^{-}_{(s-p-i-j+1)\a_{0}+(t-2p-2i-2j+3)\a_{1}})\\
&\ \times((1+q^{-1})b_{2j}-[2j+1]_{1}q^{j})
q_{1}^{(-2s+t)(2p+2i+2j-1)+2s+1}\\
&\ \times (q_{1}-q_{1}^{-1})\th(p\leq [(t+1)/2]).
\end{align*}
We rewrite $B_{p}$. Assuming that $p\leq [(t+1)/2]$ and 
using Proposition \ref{p-2} (5), we have
\bal
&\sum_{i=1}^{[(t-2p+3)/2]}
E_{(i+1)\d -\a_{1}}T_{\pa}(d_{(s-i)\a_{0}+(t-2i)\a_{1},p-1})
[2i+1]_{1}q^{(-2s+t)i}\th(p\leq [(t+1)/2])\\
&=\sum_{i=1}^{[(t-2p+3)/2]}\sum_{j=0}^{[(t-2p-2i+3)/2]}
E_{(i+1)\d -\a_{1}}\ps_{j}
T_{\pa}(D^{-}_{(s-p-i-j+1)\a_{0}+(t-2p-2i-2j+3)\a_{1}})\\
&\ \times q_{1}^{(-2s+t)(2p+2j-1)+2s-2i+1}(q_{1}-q_{1}^{-1})\\
&=\sum_{i=1}^{[(t-2p+3)/2]}\sum_{j=0}^{[(t-2p-2i+3)/2]}
\Big( \ps_{j}E_{(i+1)\d -\a_{1}}\\
&\ +(1+q^{-1})\sum_{k=0}^{j-1}b_{2(j-k)}
\ps_{k}E_{(i+j-k+1)\d -\a_{1}}\th(j\geq 1)\Big)\\
&\ \times T_{\pa}(D^{-}_{(s-p-i-j+1)\a_{0}+(t-2p-2i-2j+3)\a_{1}})
q_{1}^{(-2s+t)(2p+2j-1)+2s-2i+1}(q_{1}-q_{1}^{-1}).
\end{align*}
In the double summation, we see that 
$1\leq i\leq [(t-2p+3)/2],\ 0\leq j\leq [(t-2p-2i+3)/2]$ is equivalent to
$0\leq j\leq [(t-2p+1)/2],\ 1\leq i\leq [(t-2p-2j+3)/2]$.
In the triple one, putting $k=i',\ i+j-k=j',\ j-k=k'$, we see that
$i=j'-k',\ j=i'+k'$ and that
$1\leq i\leq [(t-2p+3)/2],\ 1\leq j\leq [(t-2p-2i+3)/2],\ 0\leq k\leq j-1$
is equivalent to
$0\leq i'\leq [(t-2p+1)/2],\ 2\leq j'\leq [(t-2p-2i'+3)/2],
\ 0\leq k'\leq j'-1$.
Hence,
\bal
B_{p}&=\sum_{i=0}^{[(t-2p+2)/2]}
\ps_{i}T_{\pa}(D^{-}_{(s-p-i+1)\a_{0}+(t-2p-2i+2)\a_{1}})
c_{i}(q_{1}-q_{1}^{-1})\\
&\ -\sum_{i=0}^{[(t-2p+1)/2]}\sum_{j=1}^{[(t-2p-2i+3)/2]}\ps_{i}
E_{(j+1)\d -\a_{1}}T_{\pa}(D^{-}_{(s-p-i-j+1)\a_{0}+(t-2p-2i-2j+2)\a_{1}})\\
&\ \times c_{ij}(q_{1}-q_{1}^{-1})\th(p\leq [(t+1)/2])
\end{align*}
where we set
\bal
c_{i}&=[4s-t+1]_{1}q_{1}^{(-2s+t)(2p+2i-1)+2s-2p-2i+2}
-[2p+2i-1]_{1}q^{(-2s+t)(p+i)}
\end{align*}
and
\bal
c_{ij}&=[2j+1]_{1}q_{1}^{(-2s+t)(2p+2i+2j-1)+2s-2p-2j+1}\\
&\ +(1+q^{-1})\sum_{k=0}^{j-1}b_{2k}[2j-2k+1]_{1}
q_{1}^{(-2s+t)(2p+2i+2j-1)+2s-2j+2k+1}\th(j\geq 1).
\end{align*}
It directly follows that
\bal
c_{i}&=[4s-t-2p-2i+2]_{1}q_{1}^{(-2s+t)(2p+2i-1)+2s+1}.
\end{align*}
Using Lemma \ref{b-4} (1), we have
\bal
c_{ij}=([2j+1]_{1}q^{j}-b_{2j}(1+q^{-1}))q_{1}^{(-2s+t)(2p+2i+2j-1)+2s+1}.
\end{align*}
Hence we have $A_{p}=B_{p}$ for $2\leq p\leq [(t+2)/2]$.
It remains to show that $A_{1}=B_{1}$.
We rewrite $A_{1}$. Using (2)$_{s}$, Lemma \ref{lem-bracket}, 
and Proposition \ref{p-2} (5), we have
\bal
&\big[ E_{\d-\a_{1}},T_{\pa}([D^{-}_{s\a_{0}+t\a_{1}},E_{\a_{1}}]_{q^{-2s+t}})
\big]_{q^{2s-t-1}}\\
&=\Big[ E_{\d-\a_{1}},T_{\pa}(D^{-}_{s\a_{0}+(t+1)\a_{1}})[t+1]_{1}\\
&\ +\sum_{i=1}^{[(t+1)/2]}\ps_{i}T_{\pa}(D^{-}_{(s-i)\a_{0}+(t-2i+1)\a_{1}})
q_{1}^{(-2s+t)(2i+1)+2s+1}\Big]_{q^{2s-t-1}}\\
&=[E_{\d-\a_{1}},T_{\pa}(D^{-}_{s\a_{0}+(t+1)\a_{1}})]_{q^{2s-t-1}}[t+1]_{1}\\
&\ +\sum_{i=1}^{[(t+1)/2]}\ps_{i}
[E_{\d-\a_{1}},T_{\pa}(D^{-}_{(s-i)\a_{0}+(t-2i+1)\a_{1}})]_{q^{2s-t-1}}
q_{1}^{(-2s+t)(2i+1)+2s+1}\\
&\ +\sum_{i=1}^{[(t+1)/2]}\sum_{j=0}^{i-1}\ps_{j}
E_{(i-j+1)\d-\a_{1}}T_{\pa}(D^{-}_{(s-i)\a_{0}+(t-2i+1)\a_{1}})\\
&\ \times b_{2(i-j)}q_{1}^{(-2s+t)(2i+1)+2s+1}(1+q^{-1}).
\end{align*}
We apply the induction hypothesis ((5) for $s$ and $t+1$) to the first term 
and (5)$_{s-1}$ to the second one. 
In the third one, putting $j=i',\ i-j=j'$, we see that $i=i'+j'$
and that $1\leq i\leq [(t+1)/2],\ 0\leq j\leq i-1$ is equivalent to
$0\leq i'\leq [(t-1)/2],\ 1\leq j'\leq [(t-2i'+1)/2]$.
Hence,
\bal
&\big[ E_{\d-\a_{1}},T_{\pa}([D^{-}_{s\a_{0}+t\a_{1}},E_{\a_{1}}]_{q^{-2s+t}})
\big]_{q^{2s-t-1}}\\
&=T_{\pa}(D^{-}_{s\a_{0}+t\a_{1}})[4s-t]_{1}[t+1]_{1}\\
&\ -\sum_{i=1}^{[(t+1)/2]}E_{(i+1)\d -\a_{1}}
T_{\pa}(D^{-}_{(s-i)\a_{0}+(t-2i+1)\a_{1}})
[2i+1]_{1}q^{(-2s+t+1)i}[t+1]_{1}\\
&\ +\sum_{i=1}^{[(t+1)/2]}\ps_{i}
\Big( T_{\pa}(D^{-}_{(s-i)\a_{0}+(t-2i)\a_{1}})[4s-t-2i]_{1}\th(t-2i\geq 0)\\
&\ -\sum_{j=1}^{[(t-2i+1)/2]}E_{(j+1)\d -\a_{1}}
T_{\pa}(D^{-}_{(s-i-j)\a_{0}+(t-2i-2j+1)\a_{1}})\\
&\ \times [2j+1]_{1}q^{(-2s+t+1)j}\th(t-2i\geq 1)\Big)
q_{1}^{(-2s+t)(2i+1)+2s+1}\\
&\ +\sum_{i=0}^{[(t-1)/2]}\sum_{j=1}^{[(t-2i+1)/2]}\ps_{i}
E_{(j+1)\d-\a_{1}}T_{\pa}(D^{-}_{(s-i-j)\a_{0}+(t-2i-2j+1)\a_{1}})\\
&\ \times b_{2j}q_{1}^{(-2s+t)(2i+2j+1)+2s+1}(1+q^{-1}).
\end{align*}
Thus,
\bal
A_{1}&=T_{\pa}(D^{-}_{s\a_{0}+t\a_{1}})
([4s-t]_{1}[t+1]_{1}(q_{1}-q_{1}^{-1})+(q^{-2s+t}-q^{2s-t}))\\
&\ +\sum_{i=1}^{[t/2]}\ps_{i}T_{\pa}(D^{-}_{(s-i)\a_{0}+(t-2i)\a_{1}})
[4s-t-2i]_{1}q_{1}^{(-2s+t)(2i+1)+2s+1}(q_{1}-q_{1}^{-1})\\
&\ -\sum_{j=1}^{[(t+1)/2]}E_{(j+1)\d -\a_{1}}
T_{\pa}(D^{-}_{(s-j)\a_{0}+(t-2j+1)\a_{1}})\\
&\ \times [2j+1]_{1}[t+1]_{1}q^{(-2s+t+1)j}(q_{1}-q_{1}^{-1})\\
&\ -\sum_{i=0}^{[(t-1)/2]}\sum_{j=1}^{[(t-2i+1)/2]}\ps_{i}E_{(j+1)\d -\a_{1}}
T_{\pa}(D^{-}_{(s-i-j)\a_{0}+(t-2i-2j+1)\a_{1}})\\
&\ \times [2j+1]_{1}[t+1]_{1}q^{(-2s+t+1)j}q_{1}^{(-2s+t)(2i+1)+2s+1}
(q_{1}-q_{1}^{-1})\\
&\ -\sum_{i=0}^{[(t-1)/2]}\sum_{j=1}^{[(t-2i+1)/2]}\ps_{i}E_{(j+1)\d -\a_{1}}
T_{\pa}(D^{-}_{(s-i-j)\a_{0}+(t-2i-2j+1)\a_{1}})\\
&\ \times b_{2j}q^{(-2s+t+1)j}q_{1}^{(-2s+t)(2i+2j+1)+2s+1}(1+q^{-1})
(q_{1}-q_{1}^{-1})\\
&=T_{\pa}(D^{-}_{s\a_{0}+t\a_{1}})[4s-t+1]_{1}[t]_{1}(q_{1}-q_{1}^{-1})\\
&\ +\sum_{i=1}^{[t/2]}\ps_{i}T_{\pa}(D^{-}_{(s-i)\a_{0}+(t-2i)\a_{1}})
[4s-t-2i]_{1}\\
&\ \times q_{1}^{(-2s+t)(2i+1)+2s+1}(q_{1}-q_{1}^{-1})\th(t\geq 2)\\
&\ +\sum_{j=1}^{[(t+1)/2]}E_{(j+1)\d -\a_{1}}
T_{\pa}(D^{-}_{(s-j)\a_{0}+(t-2j+1)\a_{1}})\\
&\ \times\big((1+q^{-1})b_{2j}q_{1}^{(-2s+t)(2j+1)+2s+1}
-[2j+1]_{1}[t+1]_{1}q^{(-2s+t+1)j}(q_{1}-q_{1}^{-1})\big)\\
&\ +\sum_{i=1}^{[(t-1)/2]}\sum_{j=1}^{[(t-2i+1)/2]}\ps_{i}E_{(j+1)\d -\a_{1}}
T_{\pa}(D^{-}_{(s-i-j)\a_{0}+(t-2i-2j+1)\a_{1}})\\
&\ \times ((1+q^{-1})b_{2j}-[2j+1]_{1}q^{j})q_{1}^{(-2s+t)(2i+2j+1)+2s+1}
(q_{1}-q_{1}^{-1})\th(t\geq 3).
\end{align*}
We rewrite $B_{1}$. Using (2)$_{s-1}$ and Proposition \ref{p-2} (5), we have
\bal
&\sum_{i=1}^{[t/2]}E_{(i+1)\d -\a_{1}}
T_{\pa}([D^{-}_{(s-i)\a_{0}+(t-2i)\a_{1}},E_{\a_{1}}]_{q^{-2s+t-1}})
[2i+1]_{1}q^{(-2s+t)i}\th(t\geq 2)\\
&=\sum_{i=1}^{[t/2]}E_{(i+1)\d -\a_{1}}
\Big( T_{\pa}(D^{-}_{(s-i)\a_{0}+(t-2i+1)\a_{1}})[t-2i+1]_{1}\\
&\ +\sum_{i=1}^{[(t-2i+1)/2]}\ps_{i}
T_{\pa}(D^{-}_{(s-i-j)\a_{0}+(t-2i-2j+1)\a_{1}})
q_{1}^{(-2s+t)(2j+1)+2(s-i)+1}\th (t-2i\geq 1)\Big) \\
&\ \times [2i+1]_{1}q^{(-2s+t)i}\th(t\geq 2)\\
&=\sum_{i=1}^{[t/2]}E_{(i+1)\d -\a_{1}}
T_{\pa}(D^{-}_{(s-i)\a_{0}+(t-2i+1)\a_{1}})[t-2i+1]_{1}
[2i+1]_{1}q^{(-2s+t)i}\th(t\geq 2)\\
&\ +\sum_{i=1}^{[t/2]}\sum_{j=1}^{[(t-2i+1)/2]}
\Big(\ps_{j}E_{(i+1)\d -\a_{1}}+(1+q^{-1})\sum_{k=0}^{j-1}b_{2(j-k)}
\ps_{k}E_{(i+j-k+1)\d -\a_{1}}\Big)\\
&\ \times T_{\pa}(D^{-}_{(s-i-j)\a_{0}+(t-2i-2j+1)\a_{1}})
[2i+1]_{1}\\
&\ \times q_{1}^{(-2s+t)(2i+2j+1)+2(s-i)+1}\th(t-2i\geq 1)\th(t\geq 3).
\end{align*}
In the first term, on account of the factor $[t-2i+1]_{1}$,
the range of $i$ can be changed to $1\leq i\leq [(t+1)/2]$,
so that the case where $t=1$ can be included.
In the double summation, we see that 
$1\leq i\leq [(t-1)/2],\ 1\leq j\leq [(t-2i+1)/2]$ is equivalent to
$1\leq j\leq [(t-1)/2],\ 1\leq i\leq [(t-2j+1)/2]$.
In the triple one, putting $k=i',\ i+j-k=j',\ j-k=k'$, we see that
$i=j'-k',\ j=i'+k'$ and that 
$1\leq i\leq [(t-1)/2],\ 1\leq j\leq [(t-2i+1)/2],\ 0\leq k\leq j-1$
is equivalent to
$0\leq i'\leq [(t-3)/2],\ 2\leq j'\leq [(t-2i'+1)/2],\ 1\leq k'\leq j'-1$.
Hence,
\bal
&\sum_{i=1}^{[t/2]}E_{(i+1)\d -\a_{1}}
T_{\pa}([D^{-}_{(s-i)\a_{0}+(t-2i)\a_{1}},E_{\a_{1}}]_{q^{-2s+t-1}})
[2i+1]_{1}q^{(-2s+t)i}\th(t\geq 2)\\
&=\sum_{i=1}^{[(t+1)/2]}E_{(i+1)\d -\a_{1}}
T_{\pa}(D^{-}_{(s-i)\a_{0}+(t-2i+1)\a_{1}})[t-2i+1]_{1}
[2i+1]_{1}q^{(-2s+t)i}\\
&\ +\sum_{i=1}^{[(t-1)/2]}\sum_{j=1}^{[(t-2j+1)/2]}\ps_{i}E_{(j+1)\d -\a_{1}}
T_{\pa}(D^{-}_{(s-i-j)\a_{0}+(t-2i-2j+1)\a_{1}})\\
&\ \times [2j+1]_{1}q_{1}^{(-2s+t)(2i+2j+1)+2(s-j)+1}\th(t\geq 3)\\
&\ +\sum_{i=0}^{[(t-3)/2]}\sum_{j=2}^{[(t-2i+1)/2]}\sum_{k=1}^{j-1}
\ps_{i}E_{(j+1)\d -\a_{1}}T_{\pa}(D^{-}_{(s-i-j)\a_{0}+(t-2i-2j+1)\a_{1}})\\
&\ \times b_{2k}[2j-2k+1]_{1}q_{1}^{(-2s+t)(2i+2j+1)+2(s-j+k)+1}
(1+q^{-1})\th(t\geq 3).
\end{align*}
Using Lemma \ref{b-4} (1), we see that the last term is equal to
\bal
&\sum_{i=0}^{[(t-3)/2]}\sum_{j=2}^{[(t-2i+1)/2]}
\ps_{i}E_{(j+1)\d -\a_{1}}T_{\pa}(D^{-}_{(s-i-j)\a_{0}+(t-2i-2j+1)\a_{1}})\\
&\ \times ((q^{j}-q^{-j})[2j+1]_{1}-b_{2j}(1+q^{-1}))
q_{1}^{(-2s+t)(2i+2j+1)+2s+1}\th(t\geq 3).
\end{align*}
Here, we can include the case where $j=1$, 
in which the factor $((q^{j}-q^{-j})[2j+1]_{1}-b_{2j}(1+q^{-1}))$ vanishes,
and the case where $i=[(t-1)/2](=[(t-3)/2]+1)$, 
in which we have $j=1$ since $[(t+1)/2]-[(t-1)/2]=1$, 
so that we can also include the case where $t=1,2$.
Hence, applying (2)$_{s}$ to the definition of $B_{1}$, we have
\bal
B_{1}&=T_{\pa}(D^{-}_{s\a_{0}+t\a_{1}})[4s-t+1]_{1}[t]_{1}(q_{1}-q_{1}^{-1})\\
&\ +(q_{1}-q_{1}^{-1})\sum_{i=1}^{[t/2]}
\ps_{i}T_{\pa}(D^{-}_{(s-i)\a_{0}+(t-2i)\a_{1}})\\
&\ \times 
\big([4s-t+1]_{1}q_{1}^{(-2s+t-1)(2i+1)+2s+1}-[2i+1]_{1}q^{(-2s+t)(i+1)}\big)
\th(t\geq 2)\\
&\ -\sum_{j=1}^{[(t+1)/2]}E_{(j+1)\d -\a_{1}}
T_{\pa}(D^{-}_{(s-j)\a_{0}+(t-2j+1)\a_{1}})\\
&\ \times [t-2j+1]_{1}[2j+1]_{1}q^{(-2s+t)j}(q_{1}-q_{1}^{-1})\\
&\ -\sum_{i=1}^{[(t-1)/2]}\sum_{j=1}^{[(t-2i+1)/2]}
\ps_{i}E_{(j+1)\d -\a_{1}}T_{\pa}(D^{-}_{(s-i-j)\a_{0}+(t-2i-2j+1)\a_{1}})\\
&\ \times [2j+1]_{1}q_{1}^{(-2s+t)(2i+2j+1)+2(s-j)+1}(q_{1}-q_{1}^{-1})
\th(t\geq 3)\\
&\ -\sum_{i=0}^{[(t-1)/2]}\sum_{j=1}^{[(t-2i+1)/2]}
\ps_{i}E_{(j+1)\d -\a_{1}}T_{\pa}(D^{-}_{(s-i-j)\a_{0}+(t-2i-2j+1)\a_{1}})\\
&\ \times ((q^{j}-q^{-j})[2j+1]_{1}-b_{2j}(1+q^{-1}))
q_{1}^{(-2s+t)(2i+2j+1)+2s+1}(q_{1}-q_{1}^{-1}).
\end{align*}
It is now easy to see that $A_{1}=B_{1}$. We obtain (5) for $s$.

The proposition is proved.
\epf


\bdefn\label{defn-D(+)}
For $n,r\geq 0$, we define the elements $D^{+}_{n\d +k\a_{1}}$ of 
$\U_{\Z}^{+}(>)\cap\U_{n\d+k\a_{1}}^{+}$ by
$D^{+}_{n\d}=\d_{n,0},\ D^{+}_{k\a_{1}}=E_{\a_{1}}^{(k)}$, and 
\bal 
D^{+}_{n\d +r\a_{1}}=\sum_{i=1}^{\min(r,2n)}E_{\a_{1}}^{(r-i)}
(T_{1}^{-1}\ast )(D^{-}_{n\a_{0}+(2n-i)\a_{1}})q_{1}^{(2n-i)(-k+i)}
\ \ \ \text{for}\ n,r \geq 1.
\end{align*} 
\edefn

\bex
Let $n\geq 0$. Then, $D^{+}_{n\d +\a_{1}}=E_{n\d +\a_{1}}$.
\eex

The following is a consequence of Proposition \ref{prop-D(-)} (5)
and will be used in the next section.

\bcor\label{cor-D(+)E}
Let $n,r\geq 0$. Then, 
\bal
\sum_{i=0}^{n}D^{+}_{(n-i)\d +r\a_{1}}E_{i\d+\a_{1}}[2i+1]_{1}q^{-ir}
=D^{+}_{n\d +(r+1)\a_{1}}[2n+r+1]_{1}.
\end{align*} 
\ecor

\bpf
We can assume that $n\geq 1$.
We rewrite the left hand side, which is denoted by LHS. 
By Definition \ref{defn-D(+)}, we have 
\bal 
\text{LHS}
&=E_{\a_{1}}^{(r)}E_{n\d+\a_{1}}[2n+1]_{1}q^{-nr}\\
&\ +\sum_{i=0}^{n-1}\sum_{j=1}^{\min(r,2n-2i)}E_{\a_{1}}^{(r-j)}
(T_{1}^{-1}\ast )(D^{-}_{(n-i)\a_{0} +(2n-2i-j)\a_{1}})
E_{i\d+\a_{1}}\\
&\ \times [2i+1]_{1}q_{1}^{(-2n+2i+j)(r-j)}q^{-ir}.
\end{align*} 
We rewrite the $l=0$ part. Using Proposition \ref{prop-D(-)} (5), we have 
\bal
(&T_{1}^{-1}\ast )(D^{-}_{n\a_{0} +(2n-j)\a_{1}})E_{\a_{1}}\\
&=(T_{1}^{-1}\ast )((T_{1}^{-1}\ast )
(E_{\a_{1}})D^{-}_{n\a_{0} +(2n-j)\a_{1}})\\
&=(T_{1}^{-1}\ast T_{\pa}^{-1}T_{\pa})
((T_{1}^{-1}\ast )(E_{\a_{1}})D^{-}_{n\a_{0} +(2n-j)\a_{1}})\\
&=(T_{1}^{-1}\ast T_{\pa}^{-1})
(E_{\d -\a_{1}}T_{\pa}(D^{-}_{n\a_{0} +(2n-j)\a_{1}}))\\
&=(T_{1}^{-1}\ast T_{\pa}^{-1})
\Big( T_{\pa}(D^{-}_{n\a_{0} +(2n-j)\a_{1}})E_{\d -\a_{1}}q^{j} \\
&\ +T_{\pa}(D^{-}_{n\a_{0} +(2n-j-1)\a_{1}})[2n+j+1]_{1}\th (j\leq 2n-1)\\
& \ -\sum_{k=1}^{[(2n-j)/2]}
E_{(k+1)\d -\a_{1}}T_{\pa}(D^{-}_{(n-k)\a_{0} +(2n-j-2k)\a_{1}})
[2k+1]_{1}q^{-jk}\th(j\leq 2n-2)\Big) \\
&=E_{\a_{1}}(T_{1}^{-1}\ast )(D^{-}_{n\a_{0} +(2n-j)\a_{1}})q^{j}\\
&\ +(T_{1}^{-1}\ast )(D^{-}_{n\a_{0} +(2n-j-1)\a_{1}})[2n+j+1]_{1}
\th (j\leq 2n-1)\\
&\ -\sum_{k=1}^{[(2n-j)/2]}(T_{1}^{-1}\ast )
(E_{k\d -\a_{1}}D^{-}_{(n-k)\a_{0} +(2n-j-2k)\a_{1}})
[2k+1]_{1}q^{-jk}\th(j\leq 2n-2).
\end{align*} 
Hence,
\bal
\text{LHS}&=E_{\a_{1}}^{(r)}E_{n\d+\a_{1}}[2n+1]_{1}q^{-nr}\\
&\ +\sum_{i=1}^{n-1}\sum_{j=1}^{\min(r,2n-2i)}E_{\a_{1}}^{(r-j)}
(T_{1}^{-1}\ast ) (E_{i\d -\a_{1}} D^{-}_{(n-i)\a_{0} +(2n-2i-j)\a_{1}})\\
&\ \times [2i+1]_{1}q_{1}^{(-2n+2i+j)(r-j)}q^{-ir}\\
&\ +\sum_{j=1}^{\min(r,2n)}E_{\a_{1}}^{(r-j)}
E_{\a_{1}}(T_{1}^{-1}\ast )(D^{-}_{n\a_{0} +(2n-j)\a_{1}})
q^{j}q_{1}^{(-2n+j)(r-j)}\\
&\ +\sum_{j=1}^{\min(r,2n-1)}E_{\a_{1}}^{(r-j)}(T_{1}^{-1}\ast )
(D^{-}_{n\a_{0} +(2n-j-1)\a_{1}})[2n+j+1]_{1}q_{1}^{(-2n+j)(k-j)}\\
&\ -\sum_{j=1}^{\min(r,2n-2)}\sum_{k=1}^{[(2n-j)/2]}E_{\a_{1}}^{(r-j)}
(T_{1}^{-1}\ast )(E_{k\d -\a_{1}}D^{-}_{(n-k)\a_{0} +(2n-j-2k)\a_{1}})\\
&\ \times [2k+1]_{1}q^{-jk}q_{1}^{(-2n+j)(r-j)}.
\end{align*} 
By Lemma \ref{lem-D(-)} (4), the first term is equal to
$E_{\a_{1}}^{(r)}(T_{1}^{-1}\ast)(D^{-}_{n\a_{0} +(2n-1)\a_{1}})
[2n+1]_{1}q_{1}^{-2nk}$,
which can be added to the fourth one for $i=0$.
In the second one, we see that $1\leq i\leq n-1,\ 1\leq j\leq \min(r,2n-2i)$
is equivalent to $1\leq j\leq \min(r,2n-2),\ 1\leq i\leq [(2n-j)/2]$,
so that it cancels out with the fifth one.
Hence, replacing $j$ by $j+1$ in the third one, we have 
\bal
\text{LHS}&=\sum_{j=0}^{\min(r,2n-1)}E_{\a_{1}}^{(r-j)}
(T_{1}^{-1}\ast )(D^{-}_{n\a_{0} +(2n-j-1)\a_{1}})\\
&\ \times ([r-j]_{1}q^{j+1}q_{1}^{(-2n+j+1)(r-j-1)}
+[2n+j+1]_{1}q_{1}^{(-2n+j)(r-j)})\\
&=\sum_{i=0}^{\min(r,2n-1)}E_{\a_{1}}^{(r-i)}
(T_{1}^{-1}\ast )(D^{-}_{n\a_{0} +(2n-j-1)\a_{1}})
[2n+r+1]_{1}q_{1}^{(-2n+j+1)(r-j)}\\
&=\sum_{j=1}^{\min(r+1,2n)}E_{\a_{1}}^{((r+1)-j)}
(T_{1}^{-1}\ast )(D^{-}_{n\a_{0} +(2n-j)\a_{1}})
q_{1}^{(-2n+j)((r+1)-j)}[2n+r+1]_{1}\\
&=D^{+}_{n\d +(r+1)\a_{1}}[2n+r+1]_{1}.
\end{align*} 
The corollary is proved.
\epf


\section{Construction of Basis II}

\bdefn\label{defn-D(+0)}
For $n,r\geq 0$, we define the elements $D^{\geq 0}_{n\d +r\a_{1}}$
of $\U_{\Z}^{+}(>)\U_{\Z}^{+}(0)\cap\U_{n\d+r\a_{1}}^{+}$ by 
\bal 
D^{\geq 0}_{n\d +r\a_{1}}
=\sum_{i=0}^{n}D^{+}_{(n-i)\d +r\a_{1}}P_{i}q^{-ir}.
\end{align*} 
\edefn

\bex
Let $n,k\geq 0$. Then, 
$D^{\geq 0}_{n\d}=P_{n},\ D^{\geq 0}_{k\a_{1}}=E_{\a_{1}}^{(k)}$.
\eex

The following is a consequence of Corollary \ref{cor-D(+)E}.

\blem\label{lem-D(+0)E}
Let $n,r\geq 0$. Then, 
\bal
D^{\geq 0}_{n\d +r\a_{1}}E_{\a_{1}}
=\sum_{i=0}^{n}D^{+}_{(n-i)\d +(r+1)\a_{1}}P_{i}[2n-2i+r+1]_{1}q^{-ir}.
\end{align*} 
\elem

\bpf
By Definition \ref{defn-D(+0)} and Lemma \ref{lem-com-PE}, 
the left hand side is equal to 
\bal
\sum_{i=0}^{n}D^{+}_{(n-i)\d +r\a_{1}}P_{i}E_{\a_{1}}q^{-ir}
=\sum_{i=0}^{n}\sum_{j=0}^{i}
D^{+}_{(n-i)\d +r\a_{1}}E_{(i-j)\d+\a_{1}}P_{j}[2i-2j+1]_{1}q^{-ir}.
\end{align*} 
Putting $j=i',\ i-j=j'$, we see that $i=i'+j'$ and that
$0\leq i\leq n,\ 0\leq j\leq i$ is equivalent to 
$0\leq i'\leq n,\ 0\leq j'\leq n-i'$. Thus, this is equal to 
\bal
\sum_{i=0}^{n}\sum_{j=0}^{n-i}
D^{+}_{(n-i-j)\d +r\a_{1}}E_{j\d+\a_{1}}P_{i}[2j+1]_{1}q^{-(i+j)r}.
\end{align*}
Applying Corollary \ref{cor-D(+)E}, we obtain the lemma.
\epf

\blem\label{lem-D(+0)ps}
Let $n\geq 1,\ r\geq 0$. Then, 
\bal
\sum_{i=1}^{n}D^{\geq 0}_{(n-i)\d +r\a_{1}}\ps_{i}q^{i(-r+1)}
=\sum_{i=1}^{n}D^{+}_{(n-i)\d +r\a_{1}}P_{i}[2i]_{1}q^{i(-r+1)}.
\end{align*} 
\elem

\bpf
By Definition \ref{defn-D(+0)}, the left hand side is equal to 
\bal
\sum_{i=1}^{n}\sum_{j=0}^{n-i}
D^{+}_{(n-i-j)\d +r\a_{1}}P_{j}\ps_{i}q^{-jr}q^{i(-r+1)}.
\end{align*}
Putting $i=i'-j$, we see that $1\leq i\leq n,\ 0\leq j\leq i-1$ is 
equivalent to $1\leq i'\leq n,\ 0\leq j\leq i'-1$. Thus, this is equal to 
\bal
\sum_{i=1}^{n}\sum_{j=0}^{i-1}
D^{+}_{(n-i)\d +r\a_{1}}P_{j}\ps_{i-j}q^{-j}q^{i(-r+1)}.
\end{align*}
Applying Definition \ref{defn-P}, we obtain the lemma.
\epf

Using Proposition \ref{prop-D(-)} (2) and Lemma \ref{lem-D(+0)E},
we prove the following.

\bthm\lab{thm-com}
Let $r,s\geq 0$. Then,
\bal 
E_{\a_{0}}^{(s)}E_{\a_{1}}^{(r)}&=\sum_{t=0}^{\min(2s-1,r)}\sum_{i=0}^{[t/2]}
D^{\geq 0}_{i\d +(r-t)\a_{1}}D^{-}_{(s-i)\a_{0} +(t-2i)\a_{1}}
q^{i(r-2s)}q_{1}^{(r-t)(-4s+t)}\th(s\geq 1)\\
&\ +D^{\geq 0}_{s\d +(r-2s)\a_{1}}\th(r\geq 2s).
\end{align*} 
Note that if $s\geq 1$, then 
$0\leq t\leq \min(2s-1,r),\ 0\leq i\leq [t/2]$ is equivalent to
$0\leq i\leq\min(s-1,[r/2]),\ 2i\leq t\leq\min(2s-1,r)$.
\ethm

\bpf
We argue by the induction on $r$. 
The case where $r=0$ is clear.
Assuming the case for $r$, we shall prove the case for $r+1$.
We can assume that $s\geq 1$. 
Using the induction hypothesis and Proposition \ref{prop-D(-)} (2), we have 
\bal
E_{\a_{0}}^{(s)}E_{\a_{1}}^{(r)}E_{\a_{1}}
&=\sum_{t=0}^{\min(2s-1,r)}\sum_{i=0}^{[t/2]}
D^{\geq 0}_{i\d +(r-t)\a_{1}}
\Big(E_{\a_{1}}D^{-}_{(s-i)\a_{0}+(t-2i)\a_{1}}q^{-2s+t}\\
&\ +D^{-}_{(s-i)\a_{0} +(t-2i+1)\a_{1}}[t-2i+1]_{1}\\
&\ +\sum_{j=1}^{[(t-2i+1)/2]}
 \ps_{j}D^{-}_{(s-i-j)\a_{0}+(t-2i-2j+1)\a_{1}}
 q_{1}^{(-2s+t)(2j+1)+2(s-i)+1}\Big)\\
&\ \times q^{i(r-2s)}q_{1}^{(r-t)(-4s+t)}
+D^{\geq 0}_{s\a_{0} +(r-2s)\a_{1}}E_{\a_{1}}\th(r\geq 2s).
\end{align*} 
Hence, if we set 
\bal
A&=\sum_{t=0}^{\min(2s-1,r)}\sum_{i=0}^{[t/2]}\sum_{j=0}^{i}
D^{+}_{(i-j)\d +(r-t+1)\a_{1}}P_{j}D^{-}_{(s-i)\a_{0}+(t-2i)\a_{1}}\\
&\ \times [2i-2j+r-t+1]_{1}
q^{-j(r-t)}q^{-2s+t}q^{i(r-2s)}q_{1}^{(r-t)(-4s+t)},\\
B&=\sum_{t=0}^{\min(2s-1,r)}\sum_{i=0}^{[t/2]}
D^{\geq 0}_{i\d +(r-t)\a_{1}}D^{-}_{(s-i)\a_{0}+(t-2i+1)\a_{1}}\\
&\ \times [t-2i+1]_{1}q^{i(r-2s)}q_{1}^{(r-t)(-4s+t)},\\
C&=\sum_{t=0}^{\min(2s-1,r)}\sum_{i=0}^{[t/2]}\sum_{j=1}^{[(t-2i+1)/2]}
D^{\geq 0}_{i\d +(r-t)\a_{1}}\ps_{j}D^{-}_{(s-i-j)\a_{0}+(t-2i-2j+1)\a_{1}}\\
&\ \times q_{1}^{(-2s+t)(2j+1)+2(s-i)+1}q^{i(r-2s)}q_{1}^{(r-t)(-4s+t)},\\
D&=\sum_{j=0}^{s}D^{+}_{(s-j)\d +(r-2s+1)\a_{1}}P_{j}
[r-2j+1]_{1}q^{-j(r-2s)}\th(r\geq 2s),
\end{align*}
then we have 
\bal
E_{\a_{0}}^{(s)}E_{\a_{1}}^{(r+1)}[r+1]_{1}=A+B+C+D.
\end{align*}
Note that in the above expression of $A$, 
on account of the factor $[2i-2j+r-t+1]_{1}$, 
the range of $t$ can be changed to $0\leq t\leq\min(2s-1,r+1)$.
Also, in the above expression of $D$,  
on account of the factor $[r-2j+1]_{1}$,
the factor $\th(r\geq 2s)$ can be changed to $\th(r+1\geq 2s)$.
As for $B$, by Definition \ref{defn-D(+0)}, we have 
\bal
B&=\sum_{t=0}^{\min(2s-1,r)}\sum_{i=0}^{[t/2]}\sum_{j=0}^{i}
D^{+}_{(i-j)\d +(r-t)\a_{1}}P_{j}D^{-}_{(s-i)\a_{0}+(t-2i+1)\a_{1}}\\
&\ \times [t-2i+1]_{1}q^{-j(r-t)}q^{i(r-2s)}q_{1}^{(r-t)(-4s+t)}.
\end{align*}
Here, on account of the factor $[t-2i+1]_{1}$, 
we can change the range of $i$ to $0\leq i\leq [(t+1)/2]$.
We rewrite $C$. 
Putting $i+j=i'$, we see that $i=i'-j$ and that 
$0\leq i\leq [t/2],\ 1\leq j\leq [(t-2i+1)/2]$ is equivalent to
$0\leq i'\leq [(t+1)/2],\ 1\leq j\leq i'$. Thus,
\bal
C&=\sum_{t=0}^{\min(2s-1,r)}\sum_{i=0}^{[(t+1)/2]}\sum_{j=1}^{i}
D^{\geq 0}_{(i-j)\d +(r-t)\a_{1}}\ps_{j}D^{-}_{(s-i)\a_{0}+(t-2i+1)\a_{1}}\\
&\ \times q_{1}^{(-2s+t)(2j+1)+2(s-i+j)+1}q^{(i-j)(r-2s)}q_{1}^{(r-t)(-4s+t)}.
\end{align*}
Applying Lemma \ref{lem-D(+0)ps}, we have  
\bal
C&=\sum_{t=0}^{\min(2s-1,r)}\sum_{i=0}^{[(t+1)/2]}\sum_{j=1}^{i}
D^{+}_{(i-j)\d +(r-t)\a_{1}}P_{j}D^{-}_{(s-i)\a_{0}+(t-2i+1)\a_{1}}\\
&\ \times [2j]_{1}q^{j(-r+t+1)}q_{1}^{t-2i+1}q^{i(r-2s)}q_{1}^{(r-t)(-4s+t)}.
\end{align*}
Thus, 
\bal
B+C
&=\sum_{t=0}^{\min(2s-1,r)}\sum_{i=0}^{[(t+1)/2]}\sum_{j=0}^{i}
D^{+}_{(i-j)\d +(r-t)\a_{1}}P_{j}D^{-}_{(s-i)\a_{0}+(t-2i+1)\a_{1}}\\
&\ \times [t-2i+2j+1]_{1}q^{j(-r+t+1)}q^{i(r-2s)}q_{1}^{(r-t)(-4s+t)}.
\end{align*}
Replacing $t$ by $t-1$, we have
\bal
B+C=&\sum_{t=1}^{\min(2s,r+1)}\sum_{i=0}^{[t/2]}\sum_{j=0}^{i}
D^{+}_{(i-j)\d +(r-t+1)\a_{1}}P_{j}D^{-}_{(s-i)\a_{0}+(t-2i)\a_{1}}\\
&\ \times [t-2i+2j]_{1}q^{j(-r+t)}q^{i(r-2s)}q_{1}^{(r-t+1)(-4s+t-1)}.
\end{align*}
Here, on account of the factor $[t-2i+2j]_{1}$, 
we can include the case where $t=0$.
Hence, 
\bal
&A+B+C\\
&=\sum_{t=0}^{\min(2s-1,r+1)}\sum_{i=0}^{[t/2]}\sum_{j=0}^{i}
D^{+}_{(i-j)\d +(r-t+1)\a_{1}}P_{j}D^{-}_{(s-i)\a_{0}+(t-2i)\a_{1}}\\
&\ \times [r+1]_{1}q^{j(-r+t-1)}q^{i(r-2s+1)}q_{1}^{(r-t+1)(-4s+t)}\\
&\ +\sum_{j=0}^{s}D^{+}_{(s-j)\d +(r-2s+1)\a_{1}}P_{j}
[2j]_{1}q^{(s-j)(r-2s)}q_{1}^{-(r-2s+1)(2s+1)}\th(r+1\geq 2s).
\end{align*}
By Definition \ref{defn-D(+0)}, we see that the former term is equal to 
\bal
\sum_{t=0}^{\min(2s-1,r+1)}\sum_{i=0}^{[t/2]}
D^{\geq 0}_{i\d +(r-t+1)\a_{1}}D^{-}_{(s-i)\a_{0}+(t-2i)\a_{1}}
[r+1]_{1}q^{i(r-2s+1)}q_{1}^{(r-t+1)(-4s+t)}
\end{align*}
and that 
the latter plus $D$ is equal to
\bal
&\sum_{j=0}^{s}D^{+}_{(s-j)\d +(r-2s+1)\a_{1}}P_{j}[r+1]_{1}
q^{-j(r-2s+1)}\th(r+1\geq 2s)\\
&=D^{\geq 0}_{s\a_{0} +(r-2s+1)\a_{1}}[r+1]_{1}\th(r+1\geq 2s).
\end{align*}
Thus, we obtain the case for $r+1$. The theorem is proved.
\epf

\bcor\label{cor-z}
We have $\U_{\Z}^{+}(0)\subset\U_{\Z}^{+}$.                     
\ecor
\bpf
It is enough to show that $P_{s}\in\U_{\Z}^{+}$ for $s\geq 0$, 
but it follows from Theorem \ref{thm-com} with $r=2s$ 
by the induction on $s$.
\epf


We shall prove that $V_{\Z}$ is closed under multiplication 
(Proposition \ref{prop-z-convex-3}).
First, we shall prove that both of $\U_{\Z}^{+}(>)\U_{\Z}^{+}(0)$
and $\U_{\Z}^{+}(0)\U_{\Z}^{+}(<)$ are 
closed under multiplication (Corollary \ref{cor-z-convex-2}).

\bprop\label{prop-z}
Let $s\geq 0$. Then,
\bit
\item[(1)]
$D^{-}_{s\a_{0} +r\a_{1}}E_{\a_{1}}^{(k)}\in V_{\Z}$ for $r,k\geq 0$,
\item[(2)]
$P_{s}E_{n\d +\a_{1}}^{(k)}\in
\U^{+}_{\Z}(>;n\d +\a_{1})\U^{+}_{\Z}(0;s)$
 for $n,k\geq 0$,
\item[(3)]
$E_{\a_{0}}^{(k)}D^{\geq 0}_{s\d +r\a_{1}}\in V_{\Z}$ for $k,r\geq 0$,
\item[(4)]
$P_{s}E_{(2n+2)\d -\a_{0}}^{(k)}\in
\U^{+}_{\Z}(>;(2n+2)\d -\a_{0})\U^{+}_{\Z}(0;s)$
 for $n,k\geq 0$.
\eit
\eprop

\bpf
We denote by (a)$_{r}$ the statement (a) for $s=0,1,\ldots ,r$.
First, we prove the following.

\bslem\label{sub}
Let $s\geq 1$ and  
assume that $(2)_{s-1}$ and $(4)_{s-1}$ are valid.
Let $\a\in R^{+}_{re}(>),\ x\in\U^{+,h}_{\Z}(>;\a)$ and let 
$y\in\U^{+,h}_{\Z}(0;s-1)$. Then,
$yx=\sum_{i}x_{i}y_{i}$ for some
$x_{i}\in\U^{+,h}_{\Z}(>;\a),\ y_{i}\in\U^{+,h}_{\Z}(0;s-1)$
with $\h(x_{i})=\h(x),\ \ih(y_{i})\leq\ih(y)$.
\eslem
\bpf
By the same argument in the proof of Proposition \ref{prop-convex-2},
the sublemma is reduced to the case where
$x=E_{n\d+\a_{1}}^{(k)}$ with $n\d+\a_{1}\geq\a$
or $x=E_{(2n+2)\d-\a_{0}}^{(k)}$ with $(2n+2)\d-\a_{0}\geq\a$,
and $y=P_{i}$ for $1\leq i\leq s-1$, 
but it follows from $(2)_{s-1}$ and $(4)_{s-1}$.
\epf
We shall prove (1)$_{s}$--(4)$_{s}$ at once by the induction on $s$.
If $s=0$, then (3) follows from Theorem \ref{thm-com},
and the others are clear. Now, assume that $s\geq 1$.

(1)$_{s}$\ \ 
We argue by the induction on $r$.
The case where $r=0$ follows from Theorem \ref{thm-com}.
We can assume that $1\leq r\leq 2s-1$. 
By Theorem \ref{thm-com}, we have
\bal
&E_{\a_{0}}^{(s)}E_{\a_{1}}^{(r)}E_{\a_{1}}^{(k)}\\
&=\sum_{i=1}^{[r/2]}\sum_{t=2i}^{r}D^{\geq 0}_{i\d+(r-t)\a_{1}}
D^{-}_{(s-i)\a_{0}+(t-2i)\a_{1}}E_{\a_{1}}^{(k)}
q^{i(r-2s)}q_{1}^{(r-t)(-4s+t)}\th(r\geq 2)\\
&\ +\sum_{t=0}^{r-1}E_{\a_{1}}^{(r-t)}
D^{-}_{s\a_{0}+t\a_{1}}E_{\a_{1}}^{(k)}q_{1}^{(r-t)(-4s+t)}
+D^{-}_{s\a_{0}+r\a_{1}}E_{\a_{1}}^{(k)}.
\end{align*} 
The left hand side is equal to $E_{\a_{0}}^{(s)}E_{\a_{1}}^{(r+k)}
\genfrac{[}{]}
{0pt}
{}{r+k}{k}_{1}$, which belongs to $V_{\Z}$ by Theorem \ref{thm-com}.
The first term on the right hand side belongs to $V_{\Z}$, since
for $1\leq i\leq [r/2],\ 2i\leq t\leq r$, we have
$D^{-}_{(s-i)\a_{0}+(t-2i)\a_{1}}E_{\a_{1}}^{(k)}\in V_{\Z}$
by (1)$_{s-1}$; hence, 
$D^{\geq 0}_{i\d+(r-t)\a_{1}}
D^{-}_{(s-i)\a_{0}+(t-2i)\a_{1}}E_{\a_{1}}^{(k)}\in V_{\Z}$
by Sublemma \ref{sub} and Proposition \ref{prop-convex-2} (2).
The second one on the right hand side belongs to $V_{\Z}$
by the induction hypothesis.
Thus, $D^{-}_{s\a_{0}+r\a_{1}}E_{\a_{1}}^{(k)}\in V_{\Z}$.
We obtain (1) for $s$.

(2)$_{s}$\ \ 
By Theorem \ref{thm-com}, we have
\bal
&E_{\a_{0}}^{(s)}E_{\a_{1}}^{(2s)}E_{\a_{1}}^{(k)}\\
&=\sum_{t=0}^{2s-1}\sum_{i=0}^{[t/2]}
D^{\geq 0}_{i\d +(2s-t)\a_{1}}D^{-}_{(s-i)\a_{0} +(t-2i)\a_{1}}E_{\a_{1}}^{(k)}
q_{1}^{(2s-t)(-4s+t)}+P_{s}E_{\a_{1}}^{(k)}.
\end{align*}
The left hand side is equal to 
$E_{\a_{0}}^{(s)}E_{\a_{1}}^{(2s+k)}
\genfrac{[}{]}
{0pt}
{}{2s+k}{k}_{1}$,
which belongs to $V_{\Z}$
by Theorem \ref{thm-com}.
The first term on the right hand side belongs to $V_{\Z}$, since
for $0\leq t\leq 2s-1,\ 0\leq i\leq [t/2]$, we have
$D^{-}_{(s-i)\a_{0} +(t-2i)\a_{1}}E_{\a_{1}}^{(k)}\in V_{\Z}$
by virtue of (1)$_{s}$; hence,
$D^{\geq 0}_{i\d +(2s-t)\a_{1}}D^{-}_{(s-i)\a_{0} +(t-2i)\a_{1}}
E_{\a_{1}}^{(k)}\in V_{\Z}$ by Sublemma \ref{sub} and 
Corollary \ref{cor-z-convex-1}.
Hence, by Proposition \ref{prop-convex-2}, we have 
$P_{s}E_{n\d +\a_{1}}^{(k)}\in V_{\Z}
\cap\U^{+}(>;n\d +\a_{1})\U^{+}(0;s)
=\U^{+}_{\Z}(>;n\d +\a_{1})\U^{+}_{\Z}(0;s)$.
We obtain (2) for $s$.

(3)$_{s}$\ \ 
By Theorem \ref{thm-com}, we have
\bal
&E_{\a_{0}}^{(k)}E_{\a_{0}}^{(s)}E_{\a_{1}}^{(r+2s)}\\
&=\sum_{i=0}^{s-1}\sum_{t=2i}^{2s-1}E_{\a_{0}}^{(k)}
D^{\geq 0}_{i\d +(r+2s-t)\a_{1}}D^{-}_{(s-i)\a_{0} +(t-2i)\a_{1}}
q^{ir}q_{1}^{(-4s+t)(r+2s-t)}
+E_{\a_{0}}^{(k)}D^{\geq 0}_{s\d +r\a_{1}}.
\end{align*}
The left hand side is equal to $E_{\a_{0}}^{(k+s)}E_{\a_{1}}^{(r+2s)}
\genfrac{[}{]}
{0pt}
{}{k+s}{k}_{0}$,
which belongs to $V_{\Z}$ by Theorem \ref{thm-com}.
The first term on the right hand side belongs to
$V_{\Z}$, since for $0\leq i\leq s-1,\ 2i\leq t\leq 2s-1$, we have
$E_{\a_{0}}^{(k)}D^{\geq 0}_{i\d +(r+2s-t)\a_{1}}\in V_{\Z}$ 
by (3)$_{s-1}$; thus, $E_{\a_{0}}^{(k)}
D^{\geq 0}_{i\d +(r+2s-t)\a_{1}}D^{-}_{(s-i)\a_{0} +(t-2i)\a_{1}}\in V_{\Z}$ 
by Corollary \ref{cor-convex-1} (4).
Hence, $E_{\a_{0}}^{(k)}D^{\geq 0}_{s\d +r\a_{1}}\in V_{\Z}$.
We obtain (3) for $s$.

(4)$_{s}$\ \ 
By Theorem \ref{thm-com}, we have
\bal
&E_{\a_{0}}^{(k)}E_{\a_{0}}^{(s)}E_{\a_{1}}^{(2s)}\\
&=\sum_{t=0}^{2s-1}\sum_{i=0}^{[t/2]}E_{\a_{0}}^{(k)}
D^{\geq 0}_{i\d +(2s-t)\a_{1}}D^{-}_{(s-i)\a_{0} +(t-2i)\a_{1}}
q_{1}^{(-4s+t)(2s-t)}
+E_{\a_{0}}^{(k)}P_{s}.
\end{align*}
The left hand side is equal to $E_{\a_{0}}^{(k+s)}E_{\a_{1}}^{(2s)}
\genfrac{[}{]}
{0pt}
{}{k+s}{k}_{0}$,
which belongs to $V_{\Z}$ by Theorem \ref{thm-com}.
The first term on the right hand side belongs to $V_{\Z}$, since
for $0\leq t\leq 2s-1,\ 0\leq i\leq [t/2]$, we have
$E_{\a_{0}}^{(k)}D^{\geq 0}_{i\d +(2s-t)\a_{1}}\in V_{\Z}$ 
by (3)$_{s-1}$; hence, $E_{\a_{0}}^{(k)}
D^{\geq 0}_{i\d +(2s-t)\a_{1}}D^{-}_{(s-i)\a_{0} +(t-2i)\a_{1}}\in V_{\Z}$
by Corollary \ref{cor-z-convex-1} (4).
Hence, by Proposition \ref{prop-convex-2}, we have 
$E_{\a_{0}}^{(k)}P_{s}\in V_{\Z}\cap\U^{+}(0;s)\U^{+}(<)
=\U^{+}_{\Z}(0;s)\U^{+}_{\Z}(<)$.
Applying $T_{\pa}^{n}$, we have 
$E_{2n\d+\a_{0}}^{(k)}P_{s}\in\U^{+}_{\Z}(0;s)\U^{+}_{\Z}(<;2n\d+\a_{0})$ 
for $n\geq 0$.
Applying $T_{1}^{-1}\ast$, we have 
$P_{s}E_{(2n+2)\d-\a_{0}}^{(k)}
\in\U^{+}_{\Z}(>;(2n+2)\d-\a_{0})\U^{+}_{\Z}(0;s)$ for $n\geq 0$.
We obtain (4) for $s$.

The proposition is proved.
\epf

\bcor\label{cor-z-convex-2}
Let $\a\in R^{+}_{re}(>),\ \b\in R^{+}_{re}(<)$ and let $n\geq 0$.
\bit
\item[(1)]
Let $x\in\U^{+,h}_{\Z}(>;\a),\ y\in\U^{+,h}_{\Z}(0;n)$. Then, 
$yx=\sum_{i}x_{i}y_{i}$ for some
$x_{i}\in\U^{+,h}_{\Z}(>;\a),\ y_{i}\in\U^{+,h}_{\Z}(0;n)$
with $\h(x_{i})=\h(x),\ \ih(y_{i})\leq\ih(y)$.
\item[(2)]
$\U^{+}_{\Z}(>;\a)\U^{+}_{\Z}(0;n)$ is closed under multiplication.
\item[(3)]
Let $y\in\U^{+,h}_{\Z}(0;n),\ z\in\U^{+,h}_{\Z}(<;\b)$. Then, 
$zy=\sum_{i}y_{i}z_{i}$ for some
$y_{i}\in\U^{+,h}_{\Z}(0;n),\ z_{i}\in\U^{+,h}_{\Z}(<;\b)$
with $\ih(y_{i})\leq\ih(y),\ \h(z_{i})=\h(z)$.
\item[(4)]
$\U^{+}_{\Z}(0;n)\U^{+}_{\Z}(<;\b)$ is closed under multiplication.
\eit
\ecor
\bpf
(1) is nothing but Sublemma \ref{sub}.
(2) follows from (1) by the same argument in the proof of Proposition 
\ref{prop-convex-2} (2).
Applying $T_{1}^{-1}\ast$ to (1) and (2), we obtain (3) and (4) respectively.
\epf

\bcor\label{cor-z-ED(+)}
Let $s,t,k\geq 0$. Then,
$E_{\a_{0}}^{(k)}D^{+}_{s\a_{0}+t\a_{1}}\in V_{\Z}$.
\ecor
\bpf
By Corollary \ref{cor-z-convex-2} (4), 
we have $V_{\Z}P_{n}\subset V_{\Z}$ for $n\geq 0$. 
Thus,  
the corollary follows from Proposition \ref{prop-z} (3) and  
Definition \ref{defn-D(+0)} by the induction on $s$.
\epf

\blem\label{lem-z-D(-)}
Let $n,k\geq 0$. Then, 
\bit
\item[(1)]
$D^{-}_{(n+1)k\a_{0}+(2n+1)k\a_{1}}=E_{(n+1)\d -\a_{1}}^{(k)}
+\sum_{i}z_{i1}z_{i2}$ 
for some $z_{ij}\in\U_{\Z}^{+,h}(<)$ with $\h (z_{ij})\leq -1$,
\item[(2)]
$D^{-}_{(2n+1)k\a_{0}+4nk\a_{1}}=E_{2n\d +\a_{0}}^{(k)}
+\sum_{i}z_{i1}z_{i2}$ 
for some $z_{ij}\in\U_{\Z}^{+,h}(<)$ with $\h (z_{ij})\leq -1$.
\eit
\elem
\bpf
This follows from Definition \ref{defn-D(-)} by the induction on $n$.
\epf

\blem\label{lem-z-D(+)}
Let $n,k\geq 0$. Then,
\bit
\item[(1)]
$D^{+}_{nk\d +k\a_{1}}=E_{n\d +\a_{1}}^{(k)}
+\sum_{i}x_{i1}x_{i2}$ 
for some $x_{ij}\in\U_{\Z}^{+,h}(>)$ with $\h (x_{ij})\geq 1$,
\item[(2)]
$D^{+}_{(2n+1)k\d +2k\a_{1}}=E_{2n\d +\a_{0}}^{(k)}
+\sum_{i}x_{i1}x_{i2}$ 
for some $x_{ij}\in\U_{\Z}^{+,h}(>)$ with $\h (x_{ij})\geq 1$.
\eit
\elem
\bpf
In view of Definition \ref{defn-D(+)}, we obtain the lemma 
by applying $T_{1}^{-1}\ast$ to Lemma \ref{lem-z-D(-)}.
\epf

\bprop\label{prop-z-convex-3}
Let $\a\in R^{+}_{re}(>),\ \b\in R^{+}_{re}(<)$ 
and let $x\in\U_{\Z}^{+,h}(>;\a),\ z\in\U_{\Z}^{+,h}(<;\b)$. Then, 
\bit
\item[(1)]
$zx=\sum_{i}x_{i}y_{i}z_{i}$ for some
$x_{i}\in\U_{\Z}^{+,h}(>;\a),\ y_{i}\in\U_{\Z}^{+,h}(0),
\ z_{i}\in\U_{\Z}^{+,h}(<;\b)$ 
with $\h (x_{i})\leq\h (x),\ \h (z_{i})\geq\h (z)$,
\item[(2)]
$\U_{\Z}^{+}(>;\a)\U_{\Z}^{+}(0)\U_{\Z}^{+}(<;\b)$ 
is closed under multiplication. 
\eit
\eprop
\bpf
We prove (1) by the induction on $\h (x)-\h (z)$.
We can assume that $\h (x)\geq 1$ and $\h (z)\leq -1$.
By the same argument in the proof of Proposition \ref{prop-convex-3} (1), 
we see that (1) is reduced to the case where
$x=E_{n\d+\a_{1}}^{(k)}$ with $n\d+\a_{1}\geq\a$ or 
$x=E_{(2n+2)\d-\a_{0}}^{(k)}$ with $(2n+2)\d-\a_{0}\geq\a$, and
$z=E_{m\d-\a_{1}}^{(l)}$ with $m\d-\a_{1}\leq\b$ or 
$z=E_{2m\d+\a_{0}}^{(l)}$ with $2m\d+\a_{0}\leq\b$.
By virtue of Proposition \ref{prop-convex-3}, we only have to show that 
$zx\in V_{\Z}$ for these cases.
The case where 
$x=E_{n\d+\a_{1}}^{(k)}$ is reduced to 
the case where 
$x=E_{\a_{1}}^{(k)}$ by using $T_{\pa}^{n}$, 
but it follows from Proposition \ref{prop-z} (1), Lemma \ref{lem-z-D(-)}, 
and the induction hypothesis
by the same argument in the proof of Proposition \ref{prop-convex-3} (1). 
The case where 
$x=E_{(2n+2)\d-\a_{0}}^{(k)}$ is reduced to 
the case where 
$x=E_{2\d-\a_{0}}^{(k)}$ and 
$z=E_{m\d-\a_{1}}^{(l)}$ or $E_{2m\d+\a_{0}}^{(l)}$
by using $T_{\pa}^{n}$, 
which in turn is reduced to the case where 
$z=E_{\a_{0}}^{(k)}$ and 
$x=E_{n\d+\a_{1}}^{(l)}$ or $E_{(2n+2)\d-\a_{0}}^{(l)}$ 
by using $T_{1}^{-1}\ast$,  
but it follows from 
Corollary \ref{cor-z-ED(+)}, Lemma \ref{lem-z-D(+)}, and 
the induction hypothesis 
by the same argument in the proof of Proposition \ref{prop-convex-3} (1). 
(1) is proved.
By the same argument in the proof of Proposition \ref{prop-convex-3} (2), 
we obtain (2).
\epf

\bthm\label{thm-z-basis}
Both of the following are $\Z [q_{1},q_{1}^{-1}]$-bases of $\U_{\Z}^{+}:$
\bal
\{ {\bf E_{c_{+}}} {\bf E_{c_{0}}} {\bf E_{c_{-}}}|
\ {\bf c_{+}}\in\oplus_{i\in R_{re}^{+}(>)}\Z_{\geq 0}^{(i)}, 
\ {\bf c_{0}}\in\oplus_{n\geq 1}\Z_{\geq 0}^{(n)},
\ {\bf c_{-}}\in\oplus_{i\in R_{re}^{+}(<)}\Z_{\geq 0}^{(i)}\},\tag{1}
\end{align*}
\bal
\{ {\bf E_{c_{+}}}{\bf S}_{\lambda}{\bf E_{c_{-}}}| 
\ {\bf c_{+}}\in\oplus_{i\in R_{re}^{+}(>)}\Z_{\geq 0}^{(i)},
\ \lambda\ \text{is a partition}, 
\ {\bf c_{-}}\in\oplus_{i\in R_{re}^{+}(<)}\Z_{\geq 0}^{(i)}\}.\tag{2}
\end{align*}
The basis in $(1)$ has the convexity 
$($see $\rm{Remark}\ \ref{rem-convex})$
and the one in $(2)$  has the quasi-orthonormality 
$($see $\rm{Lemma}\ \ref{lem-quasi})$. 
\ethm

\bpf
As is discussed at the end of Section 6, 
this follows from Proposition \ref{prop-z-convex-3} 
and Corollary \ref{cor-z}.
\epf


We shall prove that the basis in Theorem \ref{thm-z-basis} (2)
is an integral crystal basis of $\U^{+}$.
First, let us recall some properties of the canonical basis.
We set $L=\{x\in\U_{\Z}^{+}|\ (x,x)\in\A\}$, 
which is a $\Z[q_{1}^{-1}]$-submodule of $\U_{\Z}^{+}$. 
Then we have $q_{1}^{-1}L=\{x\in\U_{\Z}^{+}|\ (x,x)\in q_{1}^{-1}\A\}$. 
Let $\pi$ be the canonical projection from $L$ to ${L}/q_{1}^{-1}L$. 

\bprop\label{prop-canonical}$($see \cite{L}$)$
Let {\bf B} be the canonical basis of $\U^{+}$.
\bit
\item[(1)]
{\bf B} is quasi-orthonormal with respect to the inner product on $\U^{+}$, 
that is, $(b,b')\equiv\d_{b,b'}\mod q_{1}^{-1}\A$ for $b,b'\in {\bf B}$.
\item[(2)]
Each element of {\bf B} is homogeneous$;$
thus, ${\bf B}=\sqcup_{\a\in Q^{+}}({\bf B}\cap\U^{+}_{\a})$.
\item[(3)]
Each element of {\bf B} is invariant under $-$.
\item[(4)]
For $i\in I,\ k\geq 0$, we have $E_{\a_{i}}^{(k)}\in{\bf B}$.
\item[(5)]
Let $i\in I$ and let $x\in{\bf B}+q_{1}^{-1}{L}$. 
If $_{i}r(x)=0$, 
then $E_{\a_{i}}^{(k)}x\in{\bf B}+q_{1}^{-1}{L}$ for $k\geq 0$.
\item[(6)]
{\bf B} is a $\Q(q_{1})$-basis of $\U^{+}$ and 
a $\Z[q_{1},q_{1}^{-1}]$-basis of $\U_{\Z}^{+}$.
\item[(7)]
{\bf B} is a $\Z[q_{1}^{-1}]$-basis of $L$ and 
a $\Z[q_{1}]$-basis of $\overline{L}$.
\item[(8)]
{\bf B} is a $\Z$-basis of ${L}\cap\overline{L}$ and 
the restriction of $\pi$ gives an isomorphism of $\Z$-modules 
from ${L}\cap\overline{L}$ to ${L}/q_{1}^{-1}{L}$.
\eit
\eprop

\bdefn\label{defn-crystal}
An integral crystal basis of $\U^{+}$ is
a $\Z [q_{1}^{-1}]$-basis of $L$ that coincides with 
the canonical basis {\bf B} of $\U^{+}$ modulo $q_{1}^{-1}L$. 
\edefn

\bthm\label{thm-crystal}
Let $B$ be the $\Z [q_{1},q_{1}^{-1}]$-basis of $\U_{\Z}^{+}$ 
in $\rm{Theorem}\ \ref{thm-z-basis}\ (2)$. 
Then, $B$ is an integral crystal basis of $\U^{+}$.
\ethm
\bpf
By the quasi-orthonormality, 
$B$ is a $\Z [q_{1}^{-1}]$-basis of $L$; hence, 
$\pi (B)$ is a $\Z$-basis of ${L}/q_{1}^{-1}{L}$; thus, 
the transformation matrix with coefficients in $\Z$ 
between $\pi(B)$ and $\pi({\bf B})$ can be taken as 
the identity up to signs.
We have to show that all the signs are plus. 
For that purpose, by virtue of \cite[8.3]{L4}, it is enough to show that 
any ${\bf S}_{\lambda}$ belongs to ${\bf B}+q_{1}^{-1}L$. 
We argue by the induction on the length of $\lambda$,
which is denoted by $l(\lambda)$. 
First, assume that $l(\lambda)=0,1$.
Then, ${\bf S}_{\lambda}=P_{s}$ for some $s\geq 0$.
By a close look at Theorem \ref{thm-com} with $r=2s$ and 
by Proposition \ref{prop-inner-L}
together with Lemma \ref{lem-inner-P},  
we have $E_{\a_{0}}^{(s)}E_{\a_{1}}^{(2s)}\equiv P_{s} \mod q_{1}^{-1}{L}$. 
By Proposition \ref{prop-canonical} ((4), (5)), we have 
$E_{\a_{0}}^{(s)}E_{\a_{1}}^{(2s)}\in{\bf B}+q_{1}^{-1}{L}$; hence, 
$P_{s}\in {\bf B}+q_{1}^{-1}{L}$. Now, assume that $l(\lambda)\geq 2$. 
By the same argument in the proof of \cite[Lem.4.2]{BCP}, 
the induction proceeds.
\epf


\appendix
\section{Additional Commutation Relations}
We shall study the commutation relations 
between the real root vectors 
of height $2$ and $-2$ (Corollary \ref{cor-appA}).
\blem\label{lem-appA}
Let $k,l\geq 0$. Then,
\bal
&[E_{\a_{0}},E_{k\d +\a_{1}}E_{l\d +\a_{1}}]_{q^{-4}}\\
&=(q_{1}-q_{1}^{-1})\sum_{i=0}^{l}q^{-2}b_{2(l-i)+1}
E_{k\d +\a_{1}}\ps_{i}E_{(l-i+1)\d -\a_{1}}\\
&\ +(q_{1}-q_{1}^{-1})\sum_{i=0}^{k}q^{-1}b_{2(k-i)+1}
E_{l\d +\a_{1}}\ps_{i}E_{(k-i+1)\d -\a_{1}}\\
&\ +(q_{1}-q_{1}^{-1})(1+q^{-1})\sum_{i=0}^{k-1}\sum_{j=1}^{k-i}
q^{-1}b_{2i-1}b_{2(k-i-j+1)}
E_{(k+l-i-j+1)\d +\a_{1}}\ps_{i}E_{j\d -\a_{1}}\\
&\ +(q_{1}-q_{1}^{-1})\sum_{i=1}^{k+1}b_{2i-1}\ps_{k-i+1}\ps_{l+1}.
\end{align*}
\elem

\bpf
Using Proposition \ref{p-2} ((2), (4), (7)), we have
\bal
&E_{\a_{0}}E_{k\d +\a_{1}}E_{l\d +\a_{1}}\\
&=\Big(E_{k\d +\a_{1}}E_{\a_{0}}q^{-2}
+(q_{1}-q_{1}^{-1})\sum_{i=0}^{k}b_{2(k-i)+1}
\ps_{i}E_{(k-i+1)\d -\a_{1}}\Big)E_{l\d +\a_{1}}\\
&=E_{k\d +\a_{1}}\Big(E_{l\d +\a_{1}}E_{\a_{0}}q^{-2}
+(q_{1}-q_{1}^{-1})\sum_{i=0}^{l}b_{2(l-i)+1}
\ps_{i}E_{(l-i+1)\d -\a_{1}}\Big)q^{-2}\\
&\ +(q_{1}-q_{1}^{-1})\sum_{i=0}^{k}b_{2(k-i)+1}\ps_{i}
(E_{l\d +\a_{1}}E_{(k-i+1)\d -\a_{1}}q^{-1}+\ps_{k+l-i+1})\\
&=E_{k\d +\a_{1}}E_{l\d +\a_{1}}E_{\a_{0}}q^{-4}
+(q_{1}-q_{1}^{-1})\sum_{i=0}^{l}q^{-2}b_{2(l-i)+1}
E_{k\d +\a_{1}}\ps_{i}E_{(l-i+1)\d -\a_{1}}\\
&\ +(q_{1}-q_{1}^{-1})\sum_{i=0}^{k}q^{-1}b_{2(k-i)+1}
\Big( E_{l\d +\a_{1}}\ps_{i} +(1+q^{-1})\sum_{j=0}^{i-1}b_{2(i-j)}\\
&\ \times E_{(l+i-j)\d +\a_{1}}\ps_{j}\th(i\geq 1)\Big) E_{(k-i+1)\d -\a_{1}}
+(q_{1}-q_{1}^{-1})\sum_{i=0}^{k}b_{2(k-i)+1}\ps_{i}\ps_{k+l-i+1}.
\end{align*}
Rewriting the double summation, we obtain the lemma.
\epf

\bprop\label{prop-appA}
Let $n\geq 1$. Then,
\bal
[E_{\a_{0}},E_{2n\d -\a_{0}}]_{q^{-4}}[4]_{1}
&=\sum_{i=0}^{n-1}\sum_{j=1}^{n-i}
E_{(n+j-1)\d +\a_{1}}\ps_{i}E_{(n-i-j+1)\d -\a_{1}}\\
&\ \times b_{2j-1}b_{2(n-j)+1}(q-q^{-3})(q_{1}-q_{1}^{-1})\\
&\ +\Big(\ps_{n}^{2}+\sum_{i=1}^{n}\ps_{n-i}\ps_{n+i}(b_{2i+1}-qb_{2i-1})\Big)
(q_{1}-q_{1}^{-1}).
\end{align*}
\eprop

\bpf
By Corollary \ref{c-1}, we have
\bal
[E_{\a_{0}},E_{2n\d -\a_{0}}]_{q^{-4}}[4]_{1}
&=[E_{\a_{0}},E_{n\d +\a_{1}}E_{(n-1)\d +\a_{1}}]_{q^{-4}}\\
&\ -q[E_{\a_{0}},E_{n\d +\a_{1}}E_{(n-1)\d +\a_{1}}]_{q^{-4}}.
\end{align*}
Applying Lemma \ref{lem-appA}, we have
\bal
&[E_{\a_{0}},E_{2n\d -\a_{0}}]_{q^{-4}}[4]_{1}(q_{1}-q_{1}^{-1})^{-1}\\
&=\sum_{i=0}^{n-1}q^{-2}b_{2(n-i)-1}
E_{n\d +\a_{1}}\ps_{i}E_{(n-i)\d -\a_{1}}\\
&\ +\sum_{i=0}^{n}q^{-1}b_{2(n-i)+1}
E_{(n-1)\d +\a_{1}}\ps_{i}E_{(n-i+1)\d -\a_{1}}\\
&\ +(1+q^{-1})\sum_{i=0}^{n-1}\sum_{j=1}^{n-i}
q^{-1}b_{2j-1}b_{2(n-i-j+1)}
E_{(2n-i-j)\d +\a_{1}}\ps_{i}E_{j\d -\a_{1}}\\
&\ +\sum_{i=1}^{n+1}b_{2i-1}\ps_{n-i+1}\ps_{n+i-1}\\
& -q\Big( \sum_{i=0}^{n}q^{-2}b_{2(n-i)+1}
E_{(n-1)\d +\a_{1}}\ps_{i}E_{(n-i+1)\d -\a_{1}}\\
&\ +\sum_{i=0}^{n-1}q^{-1}b_{2(n-i)-1}
E_{n\d +\a_{1}}\ps_{i}E_{(n-i)\d -\a_{1}}\\
&\ +(1+q^{-1})\sum_{i=0}^{n-2}\sum_{j=1}^{n-i-1}
q^{-1}b_{2j-1}b_{2(n-i-j)}
E_{(2n-i-j)\d +\a_{1}}\ps_{i}E_{j\d -\a_{1}}\\
&\ +\sum_{i=1}^{n}b_{2i-1}\ps_{n-i}\ps_{n+i}\Big) .
\end{align*}
The second term and the fifth one cancel out.
Combining the first one and the sixth one, and then 
adding it to the seventh one
for $j=n-i,\ 1\leq i\leq n-1$, we have 
\bal
&[E_{\a_{0}},E_{2n\d -\a_{0}}]_{q^{-4}}[4]_{1}(q_{1}-q_{1}^{-1})^{-1}\\
&=(1+q^{-1})\sum_{i=0}^{n-1}\sum_{j=1}^{n-i}
b_{2j-1}(q^{-1}b_{2(n-i-j+1)}-b_{2(n-i-j)})
E_{(2n-i-j)\d +\a_{1}}\ps_{i}E_{j\d -\a_{1}}\\
&+\sum_{i=1}^{n+1}b_{2i-1}\ps_{n-i+1}\ps_{n+i-1}
-q\sum_{i=1}^{n}b_{2i-1}\ps_{n-i}\ps_{n+i}.
\end{align*}
Applying Lemma \ref{b-1} (2) to the first term, we obtain the proposition.
\epf

\bcor\label{cor-appA}
Let $m\geq 0,\ n\geq 1$. Then,
\bal
[E_{2m\d +\a_{0}}&,E_{2n\d -\a_{0}}]_{q^{-4}}[4]_{1}\\
&=\sum_{i=0}^{m+n-1}\sum_{j=1}^{m+n-i}
E_{(n+j-1)\d +\a_{1}}\ps_{i}E_{(2m+n-i-j+1)\d -\a_{1}}\\
&\ \times b_{2j-1}b_{2(m+n-j)+1}(q-q^{-3})(q_{1}-q_{1}^{-1})\\
&\ +\Big(\ps_{m+n}^{2}+\sum_{i=1}^{m+n}\ps_{m+n-i}\ps_{m+n+i}
(b_{2i+1}-qb_{2i-1})\Big) (q_{1}-q_{1}^{-1}).
\end{align*}
\ecor

\bpf
Applying $T_{\pa}^{m}$ to Proposition \ref{prop-appA} 
with $n$ replaced by $n+m$, we obtain the corollary.
\epf


\section{Connection with the Drinfeld Generators}

\bdefn\cite{Dr}
Let $\Dr$ be the $\Q(q_{1})$-algebra generated by
\bal
\{ x_{n}^{\pm},\ a_{k},\ K^{\pm 1},\ C^{\pm 1/2}|
\ n\in\Z,\ k\in\Z\backslash\{0\}\}
\end{align*}
with the following defining relations$:$
\bal
&C^{\pm 1/2}\ \text{are central},\ C^{1/2}C^{-1/2}=1,\
KK^{-1}=K^{-1}K=1,\tag{D1}\\
&[a_{k},K]=0,\
Kx_{n}^{\pm}K^{-1}=q^{\pm 1}x_{n}^{\pm},\tag{D2}\\
&[a_{k},x_{n}^{\pm}]=\pm\frac{[2k]_{1}}{k}
(q^{k}+q^{-k}+(-1)^{k+1})C^{\mp |k|/2}x_{n+k}^{\pm}\tag{D3},\\
&[a_{k},a_{l}]=\d_{k+l,0}\frac{[2k]_{1}}{k}(q^{k}+q^{-k}+(-1)^{k+1})
\frac{C^{k}-C^{-k}}{q_{1}-q_{1}^{-1}},\tag{D4}\\
&[x_{n}^{+},x_{m}^{-}]=
\frac{C^{(n-m)/2}\psd_{n+m}^{+}-C^{(m-n)/2}\psd_{n+m}^{-}}
{q_{1}-q_{1}^{-1}},\tag{D5}\\
&\sym_{n,m}(x_{n+2}^{\pm}x_{m}^{\pm}+(q^{\mp 1}-q^{\pm 2})
x_{n+1}^{\pm}x_{m+1}^{\pm}-q^{\pm 1}x_{n}^{\pm}x_{m+2}^{\pm})=0,\tag{D6}\\
&\sym_{k,l,m}(q_{1}^{3}x_{k\mp 1}^{\pm}x_{l}^{\pm}x_{m}^{\pm}
-(q_{1}+q_{1}^{-1})x_{k}^{\pm}x_{l\mp 1}^{\pm}x_{m}^{\pm}
+q_{1}^{-3}x_{k}^{\pm}x_{l}^{\pm}x_{m\mp 1}^{\pm})=0,\tag{D7}\\
&\sym_{k,l,m}(q_{1}^{-3}x_{k\pm 1}^{\pm}x_{l}^{\pm}x_{m}^{\pm}
-(q_{1}+q_{1}^{-1})x_{k}^{\pm}x_{l\pm 1}^{\pm}x_{m}^{\pm}
+q_{1}^{3}x_{k}^{\pm}x_{l}^{\pm}x_{m\pm 1}^{\pm})=0.\tag{D8}
\end{align*}
Here, each suffix runs through all the possible cases and we set
$\psd_{n}^{+}=0\ \text{for}\ n\leq -1,\ 
\psd_{n}^{-}=0\ \text{for}\ n\geq 1$, and
\bal
\sum_{i\geq 0}\psd _{\pm i}^{\pm}u^{i}
&=K^{\pm 1}\exp\Big(\pm (q_{1}-q_{1}^{-1})\sum_{j\geq 1}a_{\pm j}u^{j}\Big).
\end{align*}
\edefn

\bdefn
Let $\U_{\Dr}$ be the $\Q(q_{1})$-subalgebra of $\U$ generated by
$e_{i},\ f_{i},\ k_{i}^{\pm 1},\ c^{\pm 1/2}:=q^{\pm (h_{0}+2^{-1}h_{1})}$ 
for $i\in I$.
\edefn

\bprop
The following gives an isomorphism of $\Q(q_{1})$-algebras 
from $\Dr$ to $\U_{\Dr}:$
\bal
&C^{1/2}\mapsto c^{1/2},\ K\mapsto k_{1},\\
&x_{n}^{+}\mapsto T_{\pa}^{-n}(e_{1})\ \ \text{for}\ n\in\Z,\\
&x_{n}^{-}\mapsto T_{\pa}^{n}(f_{1})\ \ \text{for}\ n\in\Z,\\
&a_{k}\mapsto c^{-k/2}E_{k\d}\ \ \text{for}\ k\geq 1,\tag{$\star$}\\
&a_{-k}\mapsto c^{k/2}\Omega (E_{k\d})\ \ \text{for}\ k\geq 1,\\
&\psd^{+}_{n}\mapsto (q_{1}-q_{1}^{-1})c^{-n/2}k_{1}\ps_{n}
\ \ \text{for}\ n\geq 0,\\
&\psd^{-}_{-n}\mapsto -(q_{1}-q_{1}^{-1})c^{n/2}k_{1}^{-1}\Omega (\ps_{n})
\ \ \text{for}\ n\geq 0
\end{align*}
where we set
\bal
(q_{1}-q_{1}^{-1})\sum_{n\geq 1}E_{n\d}u^{n}
=\log\Big(1+(q_{1}-q_{1}^{-1})\sum_{k\geq 1}\ps_{k}u^{k}\Big),
\end{align*} 
or equivalently,
\bal
nE_{n\d}=n\ps_{n}-(q_{1}-q_{1}^{-1})\sum_{k\geq 1}kE_{k\d}\ps_{n-k}\th(n\geq 2)
\ \ \text{for}\ n\geq 1.
\end{align*} 
The inverse map is given by
\bal
&c^{1/2}\mapsto C^{1/2},\ k_{1}\mapsto K,\ k_{0}\mapsto CK^{-2},
\ e_{1}\mapsto x_{0}^{+},\ f_{1}\mapsto x_{0}^{-},\\
&e_{0}\mapsto [4]_{1}^{-1}q^{-1}CK^{-2}[x_{0}^{-},x_{1}^{-}]_{q},
\ f_{0}\mapsto [4]_{1}^{-1}qC^{-1}[x_{-1}^{+},x_{0}^{+}]_{q^{-1}}K^{2}.
\end{align*} 
\eprop

Sketch of Proof:  
Using the results in Section 3 and the $\Q$-linear 
anti-involution $\Omega$, we can directly check that
the correspondence \thetag{$\star$} is consistent with (D1)--(D6).
To check (D3), we use the following identity:
\bal
(1+q^{-1})nb_{2n}
&=(1+q^{-1})\sum_{i=1}^{n}b_{2(n-i)}(q^{i}-q^{-i})
(q^{i}+q^{-i}+(-1)^{i+1})\\
&\ +(q^{n}-q^{-n})(q^{n}+q^{-n}+(-1)^{n+1})q^{-2}\ \ \text{for}\ n\geq 1.
\end{align*}
We can also check (D7) and (D8), which will be written elsewhere. 
Thus, we obtain the surjective $\Q (q_{1})$-algebra homomorphism 
from $\Dr$ to $\U_{\Dr}$. 
Its injectivity follows from the 
specialization argument.
See \cite{B1} and \cite{J} for the related matter.


}




\begin{thebibliography}{99}

\bibitem[B1]{B1} J.~Beck, 
    {\em Braid group action and quantum affine algebras}, 
    Comm. Math. Phys. {\bf 165} (1994), 555-568.

\bibitem[B2]{B2} J.~Beck, 
    {\em Convex bases of PBW type for quantum affine algebras}, 
    Comm. Math. Phys. {\bf 165} (1994), 193-199.

\bibitem[BCP]{BCP} J.~Beck, V.~Chari, and A.~Pressley, 
    {\em  An algebraic characterization of the affine canonical basis}, 
    Duke Math.J. {\bf 99} (1999), 455-487.

\bibitem[CP]{CP} V.~Chari, and A.~Pressley, 
    {\em Quantum affine algebras at roots of unity}, 
    Represent. Theory {\bf 1} (1997), 280-328.

\bibitem[Da1]{Da1} I.~Damiani, 
    {\em A basis of type Poincar\'{e}-Birkhoff-Witt 
    for the quantum algebra of $\widehat{sl}(2)$}, 
    J.Algebra {\bf 161} (1993), 291-310.

\bibitem[Da2]{Da2} I.~Damiani, 
    {\em The R-matrix for $($twisted$)$ affine quantum algebras}, 
    Representations and Quantizations, China High Educ. Press, Beijing,
    2000, 89-144.

\bibitem[Dr]{Dr} V.~Drinfeld, 
    {\em A new realization of Yangians and quantized affine algebras}, 
    Soviet Math. Dokl. {\bf 36} (1988), 212-216.

\bibitem[GL]{GL} I.~Grojnowski, and G.~Lusztig, 
    {\em A comparison of bases of
    quantized enveloping algebras}, Linear algebraic groups and their
    representations, Contemp. Math. {\bf 153} (1993), 11-19.

\bibitem[J]{J} N.~Jing, 
    {\em On Drinfeld realization of quantum affine algebras},
    The Monster and Lie Algebras, 
    Ohio State Univ. Math. Res. Inst. Publ. {\bf 7},
    de Gruyter, Berlin, 1998, 195-206.

\bibitem[K1]{K1} M.~Kashiwara, 
    {\em On crystal bases of the $q$-analogue 
    of universal enveloping algebras}, 
    Duke Math. J. {\bf 63} (1991), 465-516.

\bibitem[K2]{K2} M.~Kashiwara, 
    {\em On level $0$ representations of quantized affine algebras}, 
    to appear.

\bibitem[KhT]{KhT} S.M.~Khoroshkin, and V.N.~Tolstoy, 
    {\em The uniqueness theorem for the universal R-matrix}, 
    Lett. Math. Phys. {\bf 24} (1992), 231-244.
                
\bibitem[L1]{L1} G.~Lusztig, {\em Finite dimensional Hopf algebras 
    arising from quantized universal 
    enveloping algebras}, J. Amer. Math. Soc.  {\bf 3} (1990), 257-296.  

\bibitem[L2]{L2} G.~Lusztig, {\em Canonical bases arising from quantized
    enveloping algebras}, J. Amer. Math. Soc.  {\bf 3} (1990), 447-498.  

\bibitem[L]{L} G.~Lusztig, 
    Introduction to Quantum Groups, 
    Birkh\"auser, Boston, 1993.

\bibitem[L4]{L4} G.~Lusztig, 
    {\em Braid group action and canonical bases},
    Adv. Math. {\bf 122} (1996), 237-261.

\bibitem[M]{M} I.G.~Macdonald,
    Symmetric Functions and Hall Polynomials,
    2nd ed., Oxford University Press, New York, 1995.

\end{thebibliography}
\end{document}